\RequirePackage{etex}
\documentclass[3p, sort&compress]{elsarticle}

\usepackage{lineno}
\usepackage[breaklinks,bookmarks = false]{hyperref}
\hypersetup{colorlinks, linkcolor = blue, citecolor = blue, urlcolor = blue, plainpages = false, pdfwindowui = false}

\usepackage{booktabs,multirow}
\usepackage[table,xcdraw]{xcolor}

\modulolinenumbers[5]

\journal{CAMWA}

\usepackage{amsmath, amsthm, amssymb, verbatim, comment, tikz, pgfplots, mathtools, mathrsfs, ifthen, import, float, subcaption, enumitem, todonotes, bbm}
\usepackage[nameinlink,noabbrev,capitalise]{cleveref}
\Crefname{section}{Section}{Sections}
\usetikzlibrary{cd}
\tikzcdset{every label/.append style = {font = \normalsize}}
\usepackage[ruled,vlined,linesnumbered]{algorithm2e}%for psuedo code
\pgfplotsset{compat=newest}

\usepackage{pdflscape}
\usepackage{afterpage}
\usepackage{capt-of}

\usepackage{adjustbox}
\usepackage{empheq}

\pgfplotsset{compat=newest}

\usepackage[ruled]{algorithm2e}
\SetKwComment{Comment}{/* }{ */}
\usepackage{pythonhighlight}
\usetikzlibrary{patterns}

\theoremstyle{definition}

\newtheorem*{example*}{Example}
\theoremstyle{definition}
\newtheorem{remark}{Remark}

\newcommand{\R}{\mathbb{R}}

\newcommand{\U}{\mathbb{U}}
\newcommand{\V}{\mathbb{V}}

\usepackage{xcolor}

\definecolor{cincinnati-red}{RGB}{190,0,0}
\definecolor{airforceblue}{rgb}{0.36, 0.54, 0.66}

\usepackage{tikz-network}
\usetikzlibrary{decorations.pathreplacing}
\usepackage{pgfplots}

\begin{document}

\begin{frontmatter}

\title{Optimizing Variational Physics-Informed Neural Networks\\ Using Least Squares}

% Group authors per affiliation:
\author[1]{Carlos Uriarte}
\ead{carlos.uriarte@ehu.eus}


\author[2]{Manuela Bastidas}
\ead{mbastidaso@unal.edu.co}

\author[1,3,4]{David Pardo}
\ead{david.pardo@ehu.eus}

\author[5]{Jamie M. Taylor}
\ead{jamie.taylor@cunef.edu}

\author[6]{Sergio Rojas}
\ead{sergio.rojas@monash.edu}

\address[1]{Universidad del País Vasco/Euskal Herriko Unibertsitatea (UPV/EHU), Leioa, Spain}
\address[2]{Universidad Nacional de Colombia -- Sede Medellín, Colombia}
\address[3]{Basque Center for Applied Mathematics (BCAM), Bilbao, Spain}
\address[4]{Basque Foundation for Science (Ikerbasque), Bilbao, Spain}
\address[5]{CUNEF Universidad, Madrid, Spain}
\address[6]{Monash University, Melbourne, Australia}

\begin{abstract}
Variational Physics-Informed Neural Networks often suffer from poor convergence when using stochastic gradient-descent-based optimizers. By introducing a least squares solver for the weights of the last layer of the neural network, we improve the convergence of the loss during training in most practical scenarios. This work analyzes the computational cost of the resulting hybrid least-squares/gradient-descent optimizer and explains how to implement it efficiently. In particular, we show that a traditional implementation based on backward-mode automatic differentiation leads to a prohibitively expensive algorithm. To remedy this, we propose using either forward-mode automatic differentiation or an ultraweak-type scheme that avoids the differentiation of trial functions in the discrete weak formulation. The proposed alternatives are up to one hundred times faster than the traditional one, recovering a computational cost-per-iteration similar to that of a conventional gradient-descent-based optimizer alone. To support our analysis, we derive computational estimates and conduct numerical experiments in one- and two-dimensional problems.
\end{abstract}


\end{frontmatter}

\section{Introduction}\label[section]{section1}

Several neural-network-based methods exist for solving partial differential equations (PDEs) \cite{raissi2019physics, karniadakis2021physics, blechschmidt2021three, huang2022partial, uriarte2024solving}. These methods have shown to be successful in multiple scenarios where traditional numerical methods face significant limitations, such as high-dimensional PDEs \cite{han2018solving, zang2020weak}, fractional PDEs \cite{pang2019fpinns, guo2022monte}, parametric PDEs \cite{uriarte2022finite, khoo2021solving}, and several non-linear applications \cite{brevis2024learning, brevis2022neural}.

However, when it comes to most linear, second-order models crucial in multiple engineering applications, neural networks typically fall short of being competitive against traditional finite-element- and finite-difference-based solvers. This lack of competitiveness occurs for diverse reasons. A notorious one is an excessive number of iterations required to converge, sometimes even without reaching accuracy thresholds that are usual in finite-element or finite-difference frameworks \cite{sacchetti2022neural, grossmann2024can}.

To partially overcome this problem, E. C. Cyr \textit{et al.} \cite{cyr2020robust} considered the case of a feed-forward neural network with a non-activated and unbiased output layer. This design allowed interpreting each neuron of the last hidden layer as a spanning function of a finite-dimensional vector space, while the output layer works as a linear combination of these spanning functions. Then, the authors proposed determining the optimal coefficients of the output layer using least squares (LS) when the loss function is quadratic with respect to these coefficients. As a result, they replaced the conventional gradient-descent-based (GD-based) optimizer during training with a hybrid LS/GD optimizer that, at each iteration, uses an LS solver in the last layer followed by a GD step over the weights in all hidden layers. As a result, the number of required iterations to converge when using LS/GD drastically diminished in most of their experiments, especially as increasing the number of neurons in the final hidden layer. The origins of this idea date back to the Extreme Learning Machine proposal of 2006 by G-B. Huang \textit{et al.} \cite{huang2006extreme, huang2014insight, huang2015extreme}, which in turn is inspired by the functioning of the (Multilayer) Perceptron of 1958 by F. Rosenblatt \cite{rosenblatt1958perceptron}.

Unfortunately, this hybrid LS/GD optimizer has important limitations when solving PDEs. First, while the error of the hybrid LS/GD optimizer largely outperforms a traditional GD-based optimizer in most practical scenarios, it may rapidly get trapped in some undesired local minima on some rare occasions (see Appendix A in \cite{cyr2020robust}). Additionally, increasing the number of neurons in the last hidden layer increases the LS computational cost. Indeed, to construct the LS matrix, we typically need to evaluate spatial derivatives of each spanning function that compose the last hidden layer of the neural network. This operation becomes computationally prohibitive when using automatic differentiation (AD) in backward mode, which is the one usually employed during the GD-based training of neural networks \cite{griewank2008evaluating, baydin2018automatic, margossian2019review}.

In this work, we realize there are two ways of overcoming the above computational limitations while retaining the gains of the hybrid LS/GD optimizer: (i) use forward-mode AD, or (ii) consider an ultraweak-type scheme that avoids differential computations in the spanning functions. By doing this, the cost of a single LS/GD iteration becomes similar to that of a conventional GD-based iteration. Additionally, while \cite{cyr2020robust} presented the construction of the LS system of linear equations using an overdetermined number of integration points, herein, we focus the discussion on a variational framework, following the spirit of the Petrov-Galerkin method \cite{ern2004theory, demkowicz2014overview}. Particularly, we consider Robust Variational Physics-Informed Neural Networks (RVPINNs) \cite{rojas2024robust}, which is a robust version of the original Variational Physics-Informed Neural Network (VPINN) method \cite{kharazmi2019variational,kharazmi2021hp}.

In the computational cost estimates and numerical results, we restrict for simplicity to an ideal case that uses orthonormal test functions, leading to the so-called \emph{Deep Fourier Residual (DFR)} method \cite{taylor2023deep, Manuela2023deep}. The computational cost analysis shows that using the LS/GD optimizer is comparable to the cost of employing a conventional GD-based optimizer when the network has a sufficiently large number of trainable parameters in the hidden layers. These results are confirmed by experiments measuring floating-point operations.

The remainder of this work is structured as follows. \Cref{section2} introduces RVPINNs and frames them within the context of the hybrid LS/GD optimizer. \Cref{section3} discusses the functioning of AD and \Cref{section4} analyzes the resulting computational cost estimates when applied to the LS/GD optimizer. Finally, \Cref{section5} is devoted to numerical experiments in one- and two-dimensional problems and \Cref{section7} summarizes and concludes the work.

\section{RVPINNs: A robust minimum-residual method using neural networks}\label[section]{section2}

Let $\U$ and $\V$ be two real Hilbert spaces equipped with inner products $(\cdot,\cdot)_\U$ and $(\cdot,\cdot)_\V$, respectively. Let $b:\U\times \V\longrightarrow\R$ be a bilinear form satisfying the hypotheses of the Babǔska-Lax-Milgram theorem \cite{babuvska1971error}\footnote{Equivalently, if it satisfies any of the conditions of the Banach–Nečas–Babǔska theorem \cite{saito2018notes}.}. Then, for any continuous linear functional $l:\V\longrightarrow\R$, the following variational problem admits a unique solution:
\begin{equation}\label{variational_problem}
\displaystyle{\left\{
\begin{tabular}{l}
Find $u^*\in \U$ such that\\
$b(u^*,v)=l(v), \; \forall v\in \V$,
\end{tabular}
\right.}
\end{equation} where spaces $\U$ and $\V$ are typically referred to as \emph{trial} and \emph{test} spaces, respectively. Below, we illustrate this formulation applied to Poisson's equation.
\newline

\begin{example*}[Poisson's Equation in weak form]\label{example0}
Let $\Omega\subset\mathbb{R}^d$ be a convex and Lipschitz domain, let $f\in L^2(\Omega)$, and let
\begin{subequations}
\begin{equation}
\begin{cases}
-\Delta u=f, &\\
u\vert_{\partial\Omega}=0. &
\end{cases}
\end{equation} We multiply the PDE by a test function and integrate by parts. Then, we obtain the following weak variational problem:
\begin{equation}\label{example02}
\displaystyle{\left\{
\begin{tabular}{l}
Find $u^*\in \U=H_0^1(\Omega)=\{u\in H^1(\Omega): u\vert_{\partial\Omega}=0\}$ such that\\
$\displaystyle b(u^*,v) := \int_\Omega \nabla u^*\cdot\nabla v = \int_\Omega fv =: l(v), \; \forall v\in \V=H^1_0(\Omega)$,
\end{tabular}
\right.}
\end{equation}
\end{subequations} with identical trial and test spaces that we equip with the inner product induced by the bilinear form, 
\begin{equation}
(\cdot,\cdot)_{\U}=(\cdot,\cdot)_{\V}=b(\cdot,\cdot)=(\nabla\;\cdot,\nabla\;\cdot)_{L^2(\Omega)}.
\end{equation}
\end{example*} 

\subsection{Robust residual minimization}\label{section2.1}
By invoking the Riesz Representation Theorem \cite{oden2017applied}, for each $u\in\U$, there exists a unique $r(u) \in \V$ satisfying 
\begin{equation}\label{Riesz_theorem}
(r(u),v)_\V= b(u,v)-l(v), \qquad \forall v\in \V.
\end{equation} In other words, $r(u)$ is the Riesz representative of the \emph{residual} $\mathcal{R}(u) := b(u,\cdot)-l(\cdot)$ that lies on the topological dual of $\V$. Then, the hypotheses of the Babǔska-Lax-Milgram theorem \cite{babuvska1971error} imply the following \emph{robustness} equivalence between $r(u)\in\V$ and $e(u):=u-u^*\in\U$:
\begin{equation}\label{robustness}
\frac{1}{\mu}\Vert r(u)\Vert_{\V} \leq \Vert e(u)\Vert_{\U} \leq \frac{1}{\gamma}\Vert r(u)\Vert_{\V}, \qquad u\in\U.
\end{equation} Here, $0<\gamma \leq \mu$ are the so-called  \emph{coercivity} (inf-sup) and \emph{continuity} (sup-sup) constants, respectively; and $\Vert\cdot\Vert_\V$ and $\Vert\cdot\Vert_\U$ stand for the norms induced from the inner products $(\cdot,\cdot)_\V$ and $(\cdot,\cdot)_\U$, respectively. As a result, minimizing (the norm of) $r(u)$ encompasses the minimization of (the norm of) $e(u)$.

Thus, we can reformulate \eqref{variational_problem} as the following residual minimization problem:
\begin{equation}\label{residual_minimization}
    u^* =\arg\min_{u\in\U} \Vert r(u)\Vert_\V^2.
\end{equation}

\subsection{Discretization for the test space}\label{section2.2}

Let $\{v_m\}_{m=1}^M\subset\V$ be a linearly independent set of functions and let
\begin{equation}\label{test_discretization}
   \V_M := \text{span}\{v_m\}_{m=1}^M.
\end{equation} Then, given $u\in\U$, we define the orthogonal projection of $r(u)\in\V$ onto $\V_M$, denoted $r_M(u)$, as the unique element in $\V_M$ satisfying
\begin{equation}\label{phi_M_definition}
(r_M(u),v_m)_\V=b(u,v_m)-l(v_m), \qquad 1\leq m\leq M.
\end{equation} Thus,
\begin{equation}\label{Gram_norm}
    \Vert r_M(u)\Vert_\V^2 = \mathbf{r}(u)^T \mathbf{G}^{-1} \mathbf{r}(u) =: \Vert \mathbf{r}(u)\Vert_{\mathbf{G}^{-1}}^2,
\end{equation} where $\mathbf{G} = [(v_m, v_n)_\V]_{mn}$ is the Gram matrix, $\mathbf{r}(u) := [b(u,v_m)-l(v_m)]_m$ is the residual vector, and $\Vert\cdot\Vert_{\mathbf{G}^{-1}}$ is a vector norm induced from the inverse of the symmetric and positive-definite matrix $\mathbf{G}$.

As a result, we end up with a discretized residual minimization scheme given by
\begin{equation}\label{residual_minimization_test_discrete}
    \min_{u\in\U} \Vert \mathbf{r}(u)\Vert_{\mathbf{G}^{-1}}^2.
\end{equation} By construction, the exact solution $u^*$ is a minimizer of \eqref{residual_minimization_test_discrete}, but it is no longer unique. Any $u\in\U$ satisfying $b(u,v_m)=l(v_m)$ for all $v_m\in\V_m$ is a minimizer regardless there exists $\hat{v}\in\V\setminus\V_M$ such that $b(u,\hat{v})\neq l(\hat{v})$.

\subsection{Discretization for the trial space}\label{section2.3}

We consider a vector-valued neural network $\mathbf{u}^{\boldsymbol{\alpha}}$ defined recursively as follows:
\begin{subequations}\label{neural_network_architecture}
\begin{align}
\mathbf{x}_0(x) &:= x\in\Omega\subset\mathbb{R}^d,\\
\mathbf{x}_l(x) &:= \sigma\left(\mathbf{W}_l\;\mathbf{x}_{l-1}(x) + \mathbf{b}_l\right), \qquad 1\leq l\leq L,\\
\mathbf{u}^{\boldsymbol{\alpha}}(x) &:= \mathbf{x}_{L}(x)\in\mathbb{R}^{N},\label{last_layer}
\end{align}
\end{subequations} where $\mathbf{x}_{l-1}\mapsto \mathbf{W}_l\mathbf{x}_{l-1} + \mathbf{b}_l$ is an affine transformation, $\sigma$ is a non-linear activation function that acts componentwise,  $L$ is the depth of the neural network, $d$ and $N$ are the input and output dimensions, respectively, and
\begin{equation}
\boldsymbol{\alpha} :=\{\mathbf{W}_1, \mathbf{b}_1,\mathbf{W}_2, \mathbf{b}_2,\cdots,\mathbf{W}_L, \mathbf{b}_L\}
\end{equation} is the set of all trainable parameters.

\begin{figure}[htbp]
	\begin{center}
		\begin{tikzpicture}
			
			\draw (0,0.35) node{Input};
	
			\draw (3.75,-0.8) node {$\underbrace{\hspace{6.25cm}}_\text{\hspace{1cm}}$};
			\draw (3.75,-1.1) node {$\color{red} \boldsymbol{\alpha}$};
			\draw (7.85,-0.8) node {$\underbrace{\hspace{1.35cm}}_\text{\hspace{1cm}}$};
   			\draw (7.85,-1.1) node {$\color{blue} \boldsymbol{\omega}$};
            
			\draw(9,0.35) node {Output};
			
			\Vertex[x=10,y=1.5,style={opacity=0.0,color=white!0},size=0.00000001]{inv1}
			\Vertex[x=10,y=0,style={opacity=0.0,color=white!100}]{inv2}
			
			\Vertex[x=0,y=1.5,label=$x$,color=white,size=1,fontsize=\normalsize]{X}
			\Vertex[x=2,y=3,color=white,size=0.75]{A11}
			\Vertex[x=4,y=3,color=white,size=0.75]{A21}
			\Vertex[x=7,y=3,color=white,size=0.75]{A31}
			\Vertex[x=2,y=2,color=white,size=0.75]{A12}
			\Vertex[x=4,y=2,color=white,size=0.75]{A22}
			\Vertex[x=7,y=2,color=white,size=0.75]{A32}
			\Vertex[x=2,y=0,color=white,size=0.75]{A13}
			\Vertex[x=4,y=0,color=white,size=0.75]{A23}
			\Vertex[x=7,y=0,color=white,size=0.75]{A33}
			
			\begin{scope}[shift={(A11.center)}, scale=0.1]
				\draw[color=gray,<->,>={Stealth[length=0.5mm,width=0.5mm]}](-2,0) -- (2,0) ;
				\draw[color=gray,<->,>={Stealth[length=0.5mm,width=0.5mm]}](0,-1.5) -- (0,1.5) ;
				\draw[domain=-2:2,smooth,variable=\x,purple,line width=1pt] plot ({\x},{tanh(\x)});
			\end{scope}
			\begin{scope}[shift={(A21.center)}, scale=0.1]
				\draw[color=gray,<->,>={Stealth[length=0.5mm,width=0.5mm]}](-2,0) -- (2,0) ;
				\draw[color=gray,<->,>={Stealth[length=0.5mm,width=0.5mm]}](0,-1.5) -- (0,1.5) ;
				\draw[domain=-2:2,smooth,variable=\x,purple,line width=1pt] plot ({\x},{tanh(\x)});
			\end{scope}
			\begin{scope}[shift={(A31.center)}, scale=0.1]
				\draw[color=gray,<->,>={Stealth[length=0.5mm,width=0.5mm]}](-2,0) -- (2,0) ;
				\draw[color=gray,<->,>={Stealth[length=0.5mm,width=0.5mm]}](0,-1.5) -- (0,1.5) ;
				\draw[domain=-2:2,smooth,variable=\x,purple,line width=1pt] plot ({\x},{tanh(\x)});
			\end{scope}
			\begin{scope}[shift={(A12.center)}, scale=0.1]
				\draw[color=gray,<->,>={Stealth[length=0.5mm,width=0.5mm]}](-2,0) -- (2,0) ;
				\draw[color=gray,<->,>={Stealth[length=0.5mm,width=0.5mm]}](0,-1.5) -- (0,1.5) ;
				\draw[domain=-2:2,smooth,variable=\x,purple,line width=1pt] plot ({\x},{tanh(\x)});
			\end{scope}
			\begin{scope}[shift={(A13.center)}, scale=0.1]
				\draw[color=gray,<->,>={Stealth[length=0.5mm,width=0.5mm]}](-2,0) -- (2,0) ;
				\draw[color=gray,<->,>={Stealth[length=0.5mm,width=0.5mm]}](0,-1.5) -- (0,1.5) ;
				\draw[domain=-2:2,smooth,variable=\x,purple,line width=1pt] plot ({\x},{tanh(\x)});
			\end{scope}
			\begin{scope}[shift={(A22.center)}, scale=0.1]
				\draw[color=gray,<->,>={Stealth[length=0.5mm,width=0.5mm]}](-2,0) -- (2,0) ;
				\draw[color=gray,<->,>={Stealth[length=0.5mm,width=0.5mm]}](0,-1.5) -- (0,1.5) ;
				\draw[domain=-2:2,smooth,variable=\x,purple,line width=1pt] plot ({\x},{tanh(\x)});
			\end{scope}
			\begin{scope}[shift={(A23.center)}, scale=0.1]
				\draw[color=gray,<->,>={Stealth[length=0.5mm,width=0.5mm]}](-2,0) -- (2,0) ;
				\draw[color=gray,<->,>={Stealth[length=0.5mm,width=0.5mm]}](0,-1.5) -- (0,1.5) ;
				\draw[domain=-2:2,smooth,variable=\x,purple,line width=1pt] plot ({\x},{tanh(\x)});
			\end{scope}
			\begin{scope}[shift={(A32.center)}, scale=0.1]
				\draw[color=gray,<->,>={Stealth[length=0.5mm,width=0.5mm]}](-2,0) -- (2,0) ;
				\draw[color=gray,<->,>={Stealth[length=0.5mm,width=0.5mm]}](0,-1.5) -- (0,1.5) ;
				\draw[domain=-2:2,smooth,variable=\x,purple,line width=1pt] plot ({\x},{tanh(\x)});
			\end{scope}
			\begin{scope}[shift={(A33.center)}, scale=0.1]
				\draw[color=gray,<->,>={Stealth[length=0.5mm,width=0.5mm]}](-2,0) -- (2,0) ;
				\draw[color=gray,<->,>={Stealth[length=0.5mm,width=0.5mm]}](0,-1.5) -- (0,1.5) ;
				\draw[domain=-2:2,smooth,variable=\x,purple,line width=1pt] plot ({\x},{tanh(\x)});
			\end{scope}
			
			\Vertex[x=9,y=1.5,label=$u^{\boldsymbol{\alpha},\boldsymbol{\omega}}(x)$,color=white,size=1.5,fontsize=\normalsize]{u}
			
			\Vertex[x=5.5,y=3,style={color=white},label=$\color{black}\hdots$,size=0.75]{Ah1}
			\Vertex[x=5.5,y=2,style={color=white},label=$\color{black}\hdots$,size=0.75]{Ah2}
			\Vertex[x=5.5,y=0,style={color=white},label=$\color{black}\hdots$,size=0.75]{Ah3}

			\Edge[color=red, Direct,lw=0.5pt](X)(A11)
			\Edge[color=red, Direct,lw=0.5pt](X)(A12)
			\Edge[color=red, Direct,lw=0.5pt](X)(A13)
			
			\Edge[color=white,label={$\color{black}\vdots$}](A12)(A13)
			\Edge[color=white,label={$\color{black}\vdots$}](A22)(A23)
			\Edge[color=white,label={$\color{black}\vdots$}](Ah2)(Ah3)
			\Edge[color=white,label={$\color{black}\vdots$}](A32)(A33)
			
			\Edge[color=red, Direct,lw=0.25pt](A11)(A21)
			\Edge[color=red, Direct,lw=0.25pt](A11)(A22)
			\Edge[color=red, Direct,lw=0.25pt](A11)(A23)
			\Edge[color=red, Direct,lw=0.25pt](A11)(A21)
			\Edge[color=red, Direct,lw=0.25pt](A12)(A21)
			\Edge[color=red, Direct,lw=0.25pt](A12)(A22)
			\Edge[color=red, Direct,lw=0.25pt](A12)(A23)
			\Edge[color=red, Direct,lw=0.25pt](A13)(A21)
			\Edge[color=red, Direct,lw=0.25pt](A13)(A22)
			\Edge[color=red, Direct,lw=0.25pt](A13)(A23)
			
			\Edge[color=red, Direct,lw=0.25pt](A23)(Ah1)
			\Edge[color=red, Direct,lw=0.25pt](A23)(Ah2)
			\Edge[color=red, Direct,lw=0.25pt](A23)(Ah3)
			\Edge[color=red, Direct,lw=0.25pt](A21)(Ah1)
			\Edge[color=red, Direct,lw=0.25pt](A21)(Ah2)
			\Edge[color=red, Direct,lw=0.25pt](A21)(Ah3)
			\Edge[color=red, Direct,lw=0.25pt](A23)(Ah1)
			\Edge[color=red, Direct,lw=0.25pt](A23)(Ah2)
			\Edge[color=red, Direct,lw=0.5pt](A23)(Ah3)
			
			\Edge[color=red, Direct,lw=0.5pt](Ah1)(A31)
			\Edge[color=red, Direct,lw=0.5pt](Ah1)(A32)
			\Edge[color=red, Direct,lw=0.5pt](Ah1)(A33)
			\Edge[color=red, Direct,lw=0.5pt](Ah2)(A31)
			\Edge[color=red, Direct,lw=0.5pt](Ah2)(A32)
			\Edge[color=red, Direct,lw=0.5pt](Ah2)(A33)
			\Edge[color=red, Direct,lw=0.5pt](Ah3)(A31)
			\Edge[color=red, Direct,lw=0.5pt](Ah3)(A32)
			\Edge[color=red, Direct,lw=0.5pt](Ah3)(A33)
			
			\Edge[color=blue, Direct,lw=0.5pt](A31)(u)
			\Edge[color=blue, Direct,lw=0.5pt](A32)(u)
			\Edge[color=blue, Direct,lw=0.5pt](A33)(u)
		\end{tikzpicture}
		%}
	\end{center}
	\caption{Architecture of $u^{\boldsymbol{\alpha},\boldsymbol{\omega}}$. Arrows represent trainable parameters, and graphics within neurons represent activation functions. $\boldsymbol{\alpha}$ stands for the trainable parameters corresponding to all hidden layers, in red, and $\boldsymbol{\omega}$ for the coefficients in the final linear combination (also interpretable as a non-activated and unbiased scalar-valued layer), in blue.}
	\label{figArchitecture}
\end{figure}
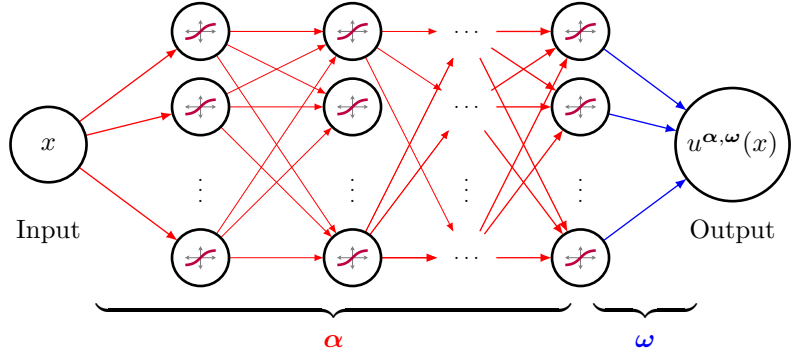

Now, following the ``adaptive basis viewpoint'' of \cite{cyr2020robust}, we interpret each component of $\mathbf{u}^{\boldsymbol{\alpha}}$ as a spanning function for a discretization of $\U$. Then, by linearly combining the components of $\mathbf{u}^{\boldsymbol{\alpha}}$,
\begin{equation}\label{NN_linear_combination}
    u^{\boldsymbol{\alpha},\boldsymbol{\omega}}(x) := \boldsymbol{\omega}\cdot\mathbf{u}^{\boldsymbol{\alpha}}(x) = \sum_{n=1}^N \omega_n \cdot u^{\boldsymbol{\alpha}}_n(x),\qquad \boldsymbol{\omega} := [\omega_n]_{n=1}^N,\quad \mathbf{u}^{\boldsymbol{\alpha}} := [u^{\boldsymbol{\alpha}}_n]_{n=1}^N,
\end{equation} we obtain a scalar-valued neural network $u^{\boldsymbol{\alpha},\boldsymbol{\omega}}:\Omega\longrightarrow\mathbb{R}$ that can be viewed as an element in a finite-dimensional vector space,
\begin{equation}
u^{\boldsymbol{\alpha},\boldsymbol{\omega}}\in \U^{\boldsymbol{\alpha}}_N,\qquad \U^{\boldsymbol{\alpha}}_N := \text{span}\{u_n^{\boldsymbol{\alpha}}\}_{n=1}^N.
\end{equation} We highlight that for each given $\boldsymbol{\alpha}$, $\text{dim}(\U^{\boldsymbol{\alpha}}_N)\leq N$ because $\{u_n^{\boldsymbol{\alpha}}\}_{n=1}^N$ might be linearly dependent. \Cref{figArchitecture} depicts the architecture of $u^{\boldsymbol{\alpha},\boldsymbol{\omega}}$.

\begin{remark}[Imposition of boundary conditions]\label{imposition_Dirichlet}
Depending on the variational formulation of the PDE, space $\U$ might include boundary conditions. To design a conforming discretization of the trial space when having homogeneous Dirichlet boundary conditions, we consider the following (re)definition of the output layer, cf. \eqref{last_layer}:
\begin{equation}\label{neural_network_architecture_dirichlet}
\mathbf{u}^{\boldsymbol{\alpha}}(x):=\xi_D(x) \cdot \mathbf{x}_L(x),\qquad x\in\Omega,
\end{equation} where $\xi_D:\Omega\longrightarrow\mathbb{R}$ is a non-trainable smooth cut-off function that vanishes only on the Dirichlet boundary \cite{berrone2023enforcing}. The multiplication between $\xi_D(x)$ and $\mathbf{x}_L(x)$ is understood componentwise: $\xi_D(x)$ multiplies each component function of $\mathbf{x}_L(x)$.
\end{remark}

\subsection{Optimization}\label{subsection2.4}

Combining \hyperref[section2.2]{Sections~\ref{section2.2}} and \ref{section2.3}, we obtain the RVPINN method described in \cite{rojas2024robust} as the following loss-function minimization:
\begin{equation}\label{loss_function}
    \min_{\boldsymbol{\alpha}, \boldsymbol{\omega}} \mathcal{L}(\boldsymbol{\alpha},\boldsymbol{\omega}), \qquad \mathcal{L}(\boldsymbol{\alpha},\boldsymbol{\omega}) := \Vert \mathbf{r}(u^{\boldsymbol{\alpha},\boldsymbol{\omega}})\Vert_{\mathbf{G}^{-1}}^2.
\end{equation}

\hyperref[section2.4.1]{Sections~\ref{section2.4.1}} and \ref{section2.4.2} consider two training approaches to carry out \eqref{loss_function} via iterative methods.

\subsubsection{Conventional GD-based optimizer}\label{section2.4.1}
Consider $\boldsymbol{\theta}:=(\boldsymbol{\alpha},\boldsymbol{\omega})$ as a single set of trainable parameters and use a GD-based iterative method on $\boldsymbol{\theta}$,
\begin{align}\label{baseline_GD}
\boldsymbol{\theta}(t+1) = \boldsymbol{\theta}(t) - \boldsymbol{\rho}(t) \cdot \nabla_{\boldsymbol{\theta}}\mathcal{L}(\boldsymbol{\theta}(t)), \qquad\qquad t\geq 1.
\end{align} Here, $\boldsymbol{\rho}(t)$ denotes the learning rate at the $t$-th iteration and $\nabla_{\boldsymbol{\theta}}$ denotes the gradient with respect to $\boldsymbol{\theta}$. We highlight that enhanced GD-based optimizers such as Adam \cite{kingma2014adam} can be simply read as in \eqref{baseline_GD} where $\boldsymbol{\rho}(t)$ absorbs the desired effect typically in the form of accumulated first- and second-order momentum estimators in previous iterations.

\subsubsection{The hybrid LS/GD optimizer}\label{section2.4.2}

We split \eqref{loss_function} in an \emph{outer-inner minimization} scheme as follows:
\begin{equation}\label{loss_function2}
    \underbrace{\min_{\boldsymbol{\alpha}} \underbrace{\min_{\boldsymbol{\omega}} \mathcal{L}(\boldsymbol{\alpha},\boldsymbol{\omega})}_{\text{inner}}}_{\text{outer}}.
\end{equation} By linearity, we can rewrite $\mathcal{L}$ in the following matrix-vector form:
\begin{equation}
\mathcal{L}(\boldsymbol{\alpha},\boldsymbol{\omega})= \Vert\underbrace{\mathbf{B}^{\boldsymbol{\alpha}}\boldsymbol{\omega}-\mathbf{l}}_{\mathbf{r}(u^{\boldsymbol{\alpha},\boldsymbol{\omega}})}\Vert_{\mathbf{G}^{-1}}^2,\qquad \mathbf{B}^{\boldsymbol{\alpha}} := [b(u^{\boldsymbol{\alpha}}_n, v_m)]_{mn},\quad \mathbf{l} := [l(v_m)]_m,
\end{equation} showing that $\mathcal{L}$ is quadratic with respect to $\boldsymbol{\omega}$. Then, we can solve the inner minimization in \eqref{loss_function2} using Least Squares (LS),
\begin{equation}\label{loss_function3}
\boldsymbol{\omega}^{\boldsymbol{\alpha}} := \arg\min_{\boldsymbol{\omega}} \Vert \mathbf{B}^{\boldsymbol{\alpha}}\boldsymbol{\omega}-\mathbf{l} \Vert_{\mathbf{G}^{-1}}^2.
\end{equation} The uniqueness of $\boldsymbol{\omega}^{\boldsymbol{\alpha}}$ is subject to the linear independence of $\{u_n^\alpha\}_{n=1}^N$.

\begin{remark}[Petrov-Galerkin viewpoint]
Given $\boldsymbol{\alpha}$, we encounter the following Petrov-Galerkin problem:
\begin{equation}\label{Petrov-Galerkin}
\displaystyle{\left\{\begin{tabular}{l}
Find $u^{\boldsymbol{\alpha},\boldsymbol{\omega}}=\boldsymbol{\omega}^{\boldsymbol{\alpha}}\cdot \mathbf{u}^{\boldsymbol{\alpha}}\in\U^{\boldsymbol{\alpha}}_N$ such that\vspace{0.5em}\\ 
$b(u^{\boldsymbol{\alpha},\boldsymbol{\omega}}, v_m) = l(v_m),\; 1\leq m\leq M$,
\end{tabular}
\right.} 
\end{equation} which, for $M>N$, produces the overdetermined system of linear equations $\mathbf{B}^{\boldsymbol{\alpha}}\boldsymbol{\omega}^{\boldsymbol{\alpha}} =\mathbf{l}$, whose minimum-residual solution results from \eqref{loss_function3}. 
\end{remark}

\begin{remark}[Least-squares regularization]
To avoid computational issues emerging from possible linear dependencies on  $\{u_n^\alpha\}_{n=1}^N$, we add the following regularization:
\begin{equation}
\boldsymbol{\omega}^{\boldsymbol{\alpha}} = \arg\min_{\boldsymbol{\omega}} \mathcal{L}(\boldsymbol{\alpha},\boldsymbol{\omega}) + \lambda \Vert \boldsymbol{\omega} \Vert^2,
\end{equation} where $\Vert\cdot\Vert$ denotes the Euclidean norm and $\lambda>0$ is small. This regularization ensures the uniqueness of the minimizer. \end{remark}

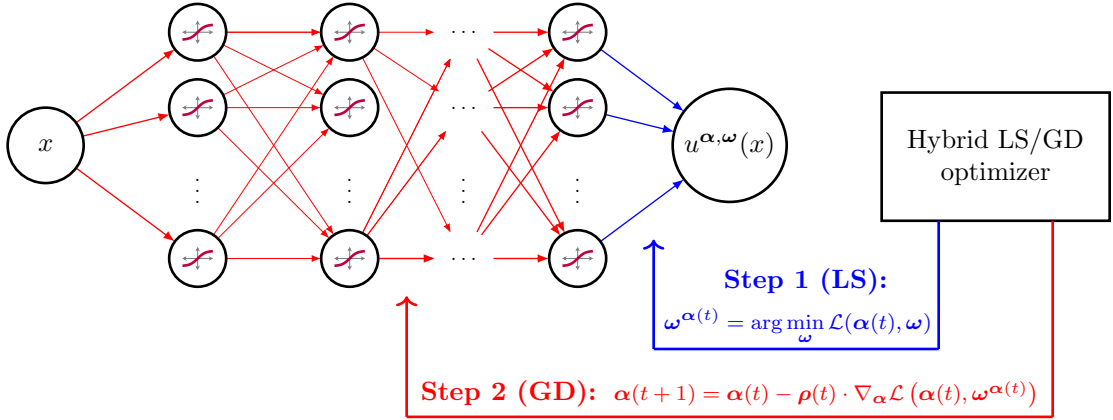
\begin{figure}[htbp]
	\begin{center}
		\begin{tikzpicture}
			\Vertex[x=10,y=1.5,style={opacity=0.0,color=white!0},size=0.00000001]{inv1}
			\Vertex[x=10,y=0,style={opacity=0.0,color=white!100}]{inv2}
            \draw[line width = 1pt, black, text width = 3cm, text centered, inner sep=0.5cm] (11,2.2) rectangle (14,0.5) node[pos=.5] {Hybrid LS/GD optimizer};
			\draw[red, line width=1pt] (13.25,0.5) -- (13.25,-2.1) ;
			\draw[red, line width=1pt] (13.25,-2.1) -- (4.75,-2.1);
			\draw[->,red, line width=1pt] (4.75,-2.1) -- (4.75,-0.5);
			\draw (9,-1.8) node {\color{red}\textbf{Step 2 (GD):\;\;}{\footnotesize$\boldsymbol{\alpha}(t+1) = \boldsymbol{\alpha}(t) - \boldsymbol{\rho}(t) \cdot \nabla_{\boldsymbol{\alpha}} \mathcal{L}\left(\boldsymbol{\alpha}(t),\boldsymbol{\omega}^{\boldsymbol{\alpha}(t)}\right)$}};

			\draw[blue, line width=1pt] (11.75,0.5) -- (11.75,-1.2);
			\draw[blue, line width=1pt] (11.75,-1.2) -- (8,-1.2);
			\draw[->, blue, line width=1pt] (8,-1.2) -- (8,0.25);
			\draw (10,-0.3) node {\textcolor{blue}{\textbf{Step 1 (LS):}}};
			\draw (9.9,-0.9) node {\footnotesize\textcolor{blue}{$\displaystyle\boldsymbol{\omega}^{\boldsymbol{\alpha}(t)} = \arg\min_{\boldsymbol{\omega}} \mathcal{L}(\boldsymbol{\alpha}(t),\boldsymbol{\omega})$}};
			
			\Vertex[x=0,y=1.5,label=$x$,color=white,size=1,fontsize=\normalsize]{X}
			\Vertex[x=2,y=3,color=white,size=0.75]{A11}
			\Vertex[x=4,y=3,color=white,size=0.75]{A21}
			\Vertex[x=7,y=3,color=white,size=0.75]{A31}
			\Vertex[x=2,y=2,color=white,size=0.75]{A12}
			\Vertex[x=4,y=2,color=white,size=0.75]{A22}
			\Vertex[x=7,y=2,color=white,size=0.75]{A32}
			\Vertex[x=2,y=0,color=white,size=0.75]{A13}
			\Vertex[x=4,y=0,color=white,size=0.75]{A23}
			\Vertex[x=7,y=0,color=white,size=0.75]{A33}
			
			\begin{scope}[shift={(A11.center)}, scale=0.1]
				\draw[color=gray,<->,>={Stealth[length=0.5mm,width=0.5mm]}](-2,0) -- (2,0)  ;
				\draw[color=gray,<->,>={Stealth[length=0.5mm,width=0.5mm]}](0,-1.5) -- (0,1.5) ;
				\draw[domain=-2:2,smooth,variable=\x,purple,line width=1pt] plot ({\x},{tanh(\x)});
			\end{scope}
			\begin{scope}[shift={(A21.center)}, scale=0.1]
				\draw[color=gray,<->,>={Stealth[length=0.5mm,width=0.5mm]}](-2,0) -- (2,0)  ;
				\draw[color=gray,<->,>={Stealth[length=0.5mm,width=0.5mm]}](0,-1.5) -- (0,1.5) ;
				\draw[domain=-2:2,smooth,variable=\x,purple,line width=1pt] plot ({\x},{tanh(\x)});
			\end{scope}
			\begin{scope}[shift={(A31.center)}, scale=0.1]
				\draw[color=gray,<->,>={Stealth[length=0.5mm,width=0.5mm]}](-2,0) -- (2,0)  ;
				\draw[color=gray,<->,>={Stealth[length=0.5mm,width=0.5mm]}](0,-1.5) -- (0,1.5) ;
				\draw[domain=-2:2,smooth,variable=\x,purple,line width=1pt] plot ({\x},{tanh(\x)});
			\end{scope}
			\begin{scope}[shift={(A12.center)}, scale=0.1]
				\draw[color=gray,<->,>={Stealth[length=0.5mm,width=0.5mm]}](-2,0) -- (2,0)  ;
				\draw[color=gray,<->,>={Stealth[length=0.5mm,width=0.5mm]}](0,-1.5) -- (0,1.5) ;
				\draw[domain=-2:2,smooth,variable=\x,purple,line width=1pt] plot ({\x},{tanh(\x)});
			\end{scope}
			\begin{scope}[shift={(A13.center)}, scale=0.1]
				\draw[color=gray,<->,>={Stealth[length=0.5mm,width=0.5mm]}](-2,0) -- (2,0)  ;
				\draw[color=gray,<->,>={Stealth[length=0.5mm,width=0.5mm]}](0,-1.5) -- (0,1.5) ;
				\draw[domain=-2:2,smooth,variable=\x,purple,line width=1pt] plot ({\x},{tanh(\x)});
			\end{scope}
			\begin{scope}[shift={(A22.center)}, scale=0.1]
				\draw[color=gray,<->,>={Stealth[length=0.5mm,width=0.5mm]}](-2,0) -- (2,0)  ;
				\draw[color=gray,<->,>={Stealth[length=0.5mm,width=0.5mm]}](0,-1.5) -- (0,1.5) ;
				\draw[domain=-2:2,smooth,variable=\x,purple,line width=1pt] plot ({\x},{tanh(\x)});
			\end{scope}
			\begin{scope}[shift={(A23.center)}, scale=0.1]
				\draw[color=gray,<->,>={Stealth[length=0.5mm,width=0.5mm]}](-2,0) -- (2,0)  ;
				\draw[color=gray,<->,>={Stealth[length=0.5mm,width=0.5mm]}](0,-1.5) -- (0,1.5) ;
				\draw[domain=-2:2,smooth,variable=\x,purple,line width=1pt] plot ({\x},{tanh(\x)});
			\end{scope}
			\begin{scope}[shift={(A32.center)}, scale=0.1]
				\draw[color=gray,<->,>={Stealth[length=0.5mm,width=0.5mm]}](-2,0) -- (2,0)  ;
				\draw[color=gray,<->,>={Stealth[length=0.5mm,width=0.5mm]}](0,-1.5) -- (0,1.5) ;
				\draw[domain=-2:2,smooth,variable=\x,purple,line width=1pt] plot ({\x},{tanh(\x)});
			\end{scope}
			\begin{scope}[shift={(A33.center)}, scale=0.1]
				\draw[color=gray,<->,>={Stealth[length=0.5mm,width=0.5mm]}](-2,0) -- (2,0)  ;
				\draw[color=gray,<->,>={Stealth[length=0.5mm,width=0.5mm]}](0,-1.5) -- (0,1.5) ;
				\draw[domain=-2:2,smooth,variable=\x,purple,line width=1pt] plot ({\x},{tanh(\x)});
			\end{scope}
			
			\Vertex[x=9,y=1.5,label=$u^{\boldsymbol{\alpha},\boldsymbol{\omega}}(x)$,color=white,size=1.5,fontsize=\normalsize]{tildeu}
			
			\Vertex[x=5.5,y=3,style={color=white},label=$\color{black}\hdots$,size=0.75]{Ah1}
			\Vertex[x=5.5,y=2,style={color=white},label=$\color{black}\hdots$,size=0.75]{Ah2}
			\Vertex[x=5.5,y=0,style={color=white},label=$\color{black}\hdots$,size=0.75]{Ah3}

			\Edge[color=red, Direct,lw=0.5pt](X)(A11)
			\Edge[color=red, Direct,lw=0.5pt](X)(A12)
			\Edge[color=red, Direct,lw=0.5pt](X)(A13)
			
			\Edge[color=white,label={$\color{black}\vdots$}](A12)(A13)
			\Edge[color=white,label={$\color{black}\vdots$}](A22)(A23)
			\Edge[color=white,label={$\color{black}\vdots$}](Ah2)(Ah3)
			\Edge[color=white,label={$\color{black}\vdots$}](A32)(A33)
			
			\Edge[color=red, Direct,lw=0.25pt](A11)(A21)
			\Edge[color=red, Direct,lw=0.25pt](A11)(A22)
			\Edge[color=red, Direct,lw=0.25pt](A11)(A23)
			\Edge[color=red, Direct,lw=0.25pt](A11)(A21)
			\Edge[color=red, Direct,lw=0.25pt](A12)(A21)
			\Edge[color=red, Direct,lw=0.25pt](A12)(A22)
			\Edge[color=red, Direct,lw=0.25pt](A12)(A23)
			\Edge[color=red, Direct,lw=0.25pt](A13)(A21)
			\Edge[color=red, Direct,lw=0.25pt](A13)(A22)
			\Edge[color=red, Direct,lw=0.25pt](A13)(A23)
			
			\Edge[color=red, Direct,lw=0.25pt](A23)(Ah1)
			\Edge[color=red, Direct,lw=0.25pt](A23)(Ah2)
			\Edge[color=red, Direct,lw=0.25pt](A23)(Ah3)
			\Edge[color=red, Direct,lw=0.25pt](A21)(Ah1)
			\Edge[color=red, Direct,lw=0.25pt](A21)(Ah2)
			\Edge[color=red, Direct,lw=0.25pt](A21)(Ah3)
			\Edge[color=red, Direct,lw=0.25pt](A23)(Ah1)
			\Edge[color=red, Direct,lw=0.25pt](A23)(Ah2)
			\Edge[color=red, Direct,lw=0.5pt](A23)(Ah3)
			
			\Edge[color=red, Direct,lw=0.5pt](Ah1)(A31)
			\Edge[color=red, Direct,lw=0.5pt](Ah1)(A32)
			\Edge[color=red, Direct,lw=0.5pt](Ah1)(A33)
			\Edge[color=red, Direct,lw=0.5pt](Ah2)(A31)
			\Edge[color=red, Direct,lw=0.5pt](Ah2)(A32)
			\Edge[color=red, Direct,lw=0.5pt](Ah2)(A33)
			\Edge[color=red, Direct,lw=0.5pt](Ah3)(A31)
			\Edge[color=red, Direct,lw=0.5pt](Ah3)(A32)
			\Edge[color=red, Direct,lw=0.5pt](Ah3)(A33)
			
			\Edge[color=blue, Direct,lw=0.5pt](A31)(tildeu)
			\Edge[color=blue, Direct,lw=0.5pt](A32)(tildeu)
			\Edge[color=blue, Direct,lw=0.5pt](A33)(tildeu)
			
		\end{tikzpicture}

	\end{center}
	\caption{Methodology of the hybrid LS/GD optimizer.}
	\label{figArchitecture2}
\end{figure}

Hence, the hybrid LS/GD optimizer is a two-step iterative method that, at each training iteration $t\geq 1$, (Step 1) finds $\boldsymbol{\omega}^{\boldsymbol{\alpha}(t)}$ using an LS solver and then (Step 2) updates $\boldsymbol{\alpha}(t)$ via a usual GD-based scheme,
\begin{align}\label{GD_alpha}
\boldsymbol{\alpha}(t+1) = \boldsymbol{\alpha}(t) - \boldsymbol{\rho}(t) \cdot \nabla_{\boldsymbol{\alpha}}\mathcal{L}\left(\boldsymbol{\alpha}(t),\boldsymbol{\omega}^{\boldsymbol{\alpha}(t)}\right).
\end{align} \Cref{figArchitecture2} illustrates the hybrid scheme.

\begin{remark}[Total derivative viewpoint]\label{total_partial_derivative}
As $\boldsymbol{\omega}^{\boldsymbol{\alpha}}$ is determined by $\boldsymbol{\alpha}$, we may view $\mathcal{L}(\boldsymbol{\alpha},\boldsymbol{\omega}^{\boldsymbol{\alpha}})$ as a function exclusively of $\boldsymbol{\alpha}$. Then, its \emph{total derivative} reads in terms of \emph{partial derivatives}\footnote{Here, by \emph{partial derivative} $\frac{\partial\mathcal{L}}{\partial\boldsymbol{\alpha}}(\boldsymbol{\alpha},\boldsymbol{\omega}^{\boldsymbol{\alpha}})$, we refer to the technique that treats $\boldsymbol{\omega}^{\boldsymbol{\alpha}}$ as a constant during differentiation. Similarly for $\frac{\partial\mathcal{L}}{\partial\boldsymbol{\omega}^{\boldsymbol{\alpha}}}(\boldsymbol{\alpha},\boldsymbol{\omega}^{\boldsymbol{\alpha}})$.} as follows:
\begin{equation}
\frac{d \mathcal{L}}{d \boldsymbol{\alpha}} = \left(\frac{\partial\mathcal{L}}{\partial \boldsymbol{\alpha}} \cdot \frac{d\boldsymbol{\alpha}}{d\boldsymbol{\alpha}}\right) + \left(\frac{\partial\mathcal{L}}{\partial \boldsymbol{\omega}^{\boldsymbol{\alpha}}} \cdot \frac{d\boldsymbol{\omega}^{\boldsymbol{\alpha}}}{d\boldsymbol{\alpha}}\right).
\end{equation}
Due to the definition of $\boldsymbol{\omega}^{\boldsymbol{\alpha}}$ via the optimality condition \eqref{loss_function3} (disregarding the effect of the regularizer), the second addend vanishes, showing that total and partial derivatives coincide. In \Cref{section4.4}, we will take advantage of this to enhance computational savings.
\end{remark}

\section{Automatic differentiation: forward and backward modes}\label{section3}

The above optimization algorithms require not only derivatives with respect to the trainable parameters but also with respect to the spatial variable corresponding to the input of the neural network. We use \emph{automatic differentiation} (AD) to compute these derivatives during training. AD mainly operates in two modes: \emph{forward} and \emph{backward} (or \emph{reverse}). The most popular choice is the backward mode, also known as the backpropagation algorithm, because it is the most efficient for computing the gradient of the loss function with respect to the trainable parameters. However, the forward mode might be orders of magnitude more efficient to compute the derivatives of the neural network with respect to its input.

In the \emph{forward accumulation}, numerical seeds at the input of the function traverse the chain-rule graph from the input to the output; in the \emph{backward accumulation} mode, once the forward-pass evaluation is made, numerical seeds deposited at the output follow the opposite path (from output to input) \cite{griewank2008evaluating, baydin2018automatic, margossian2019review}. In each case, the cost of evaluating the derivative is bounded by the multiplication of the following three terms: (i) the cost of evaluating the function itself,  (ii) the dimension of the seed (i.e., the dimension of the input or the output in forward and backward accumulations, respectively), and (iii) a constant. This constant is smaller than or equal to three under specific theoretical conditions (Baur-Strassen Theorem \cite{baur1983complexity}). We summarize this in \Cref{table_AD}.

\begin{table}[htbp]
\centering
\begin{tabular}{ c | c c }
\toprule
\textbf{Cost}								& Forward  AD & Backward AD\\ \midrule
$\phantom{\nabla^2} \mathbf{u}^{\boldsymbol{\alpha}}(x)$ & $C_{net}$ & $C_{net}$  \\ \midrule
$\phantom{^2}\nabla \mathbf{u}^{\boldsymbol{\alpha}}(x)$ & $\mathcal{O}(dC_{net})$ & $\mathcal{O}(N C_{net})$ \\ \bottomrule
\end{tabular}
\caption{Costs of computing $\mathbf{u}^{\boldsymbol{\alpha}}(x) = [u_n^{\boldsymbol{\alpha}}(x)]_{n=1}^N$ and $\nabla \mathbf{u}^{\boldsymbol{\alpha}}(x) = [\nabla u_n^{\boldsymbol{\alpha}}(x)]_{n=1}^N$ using forward and backward AD. $C_{net}$ is the reference cost of evaluating the neural network at a single input sample $x\in\Omega\subset\mathbb{R}^d$. This cost depends on the specifications of the chosen architecture (number of trainable parameters, chosen activation functions, layer distribution, etc.).}
\label[table]{table_AD}
\end{table}

\begin{figure}[htbp]
\centering
     \begin{subfigure}{0.48\textwidth}
     \begin{tikzpicture}
     \begin{axis}[
        xlabel = {input dimension ($d$)},
        ylabel = {output dimension ($N$)},
       height=.75*\textwidth,
        width=.75*\textwidth,
        grid=major,
        xmode=log,
        log base x=2,
        ymode=log,
        log base y=2,
        zmode=log,
        log base z=2,
        view={0}{90},
        colorbar,
        colorbar style={yticklabel=$2^{\pgfmathprintnumber{\tick}}$},
        title={$\dfrac{\# \text{ FLOPs of }\nabla\mathbf{u}^{\boldsymbol{\alpha}}(x) }{ \# \text{ FLOPs of }\mathbf{u}^{\boldsymbol{\alpha}}(x)}$},
	]
	 	
        \addplot3[surf, scatter, only marks, opacity=1, scatter/@pre marker code/.append style={/tikz/mark size=3}] table[x=input_dim, y=output_dim, z=fwd, col sep=comma]{FIGURES/section3/ratios_AD.csv};
     \end{axis}
     \end{tikzpicture}
     \caption{Forward AD.}
     \end{subfigure}\hfill
     \begin{subfigure}{0.48\textwidth}
     \begin{tikzpicture}
     \begin{axis}[
        xlabel = {input dimension ($d$)},
        ylabel = {output dimension ($N$)},
       height=.75*\textwidth,
        width=.75*\textwidth,
        grid=major,
        xmode=log,
        log base x=2,
        ymode=log,
        log base y=2,
        zmode=log,
        log base z=2,
        view={0}{90},
        colorbar,
        colorbar style={yticklabel=$2^{\pgfmathprintnumber{\tick}}$},
        title={$\dfrac{\# \text{ FLOPs of }\nabla\mathbf{u}^{\boldsymbol{\alpha}}(x) }{ \# \text{ FLOPs of }\mathbf{u}^{\boldsymbol{\alpha}}(x)}$},
	]
	 	
        \addplot3[surf, scatter, only marks, opacity=1, scatter/@pre marker code/.append style={/tikz/mark size=3}] table[x=input_dim, y=output_dim, z=bwd, col sep=comma]{FIGURES/section3/ratios_AD.csv};
     \end{axis}
     \end{tikzpicture}
     \caption{Backward AD.}
     \end{subfigure}\caption{Costs of computing the gradient of a vector-valued neural network $\mathbb{R}^d\ni x \mapsto \mathbf{u}^{\boldsymbol{\alpha}}(x)\in \mathbb{R}^N$ with five hidden layers of width $2^{10}$. Experimentation is performed for $d, N\in\{2^0,2^2,\ldots,2^9\}$.}
\label{fig:flopsAD_experiment}
\end{figure}
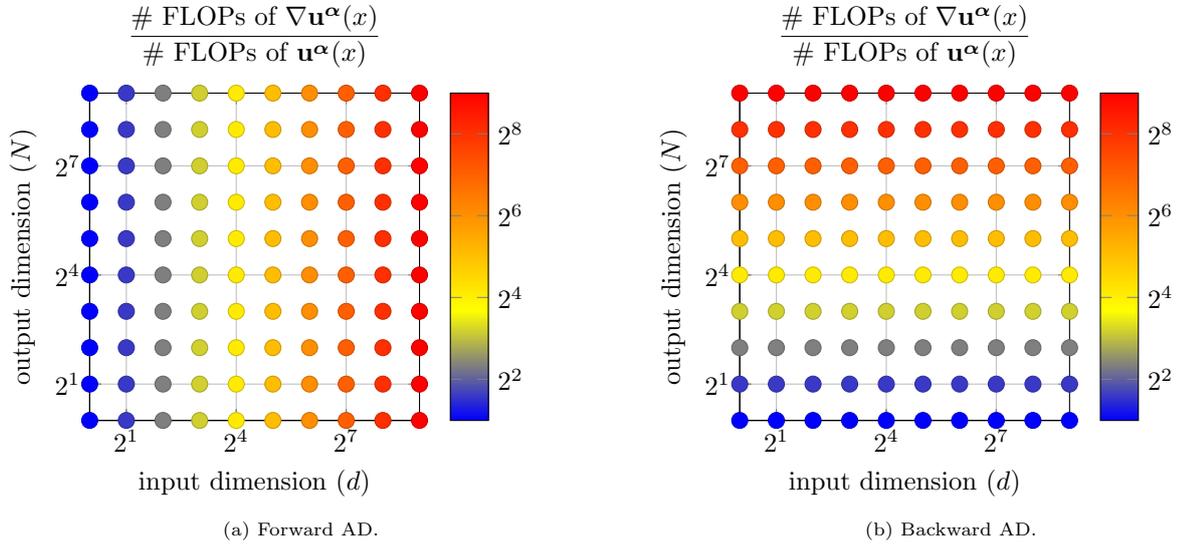

To illustrate the above, we consider a collection of (vector-valued) neural networks with different input and output dimensions and measure the amount of floating-point operations (FLOPs)\footnote{We carried out the experiments in an Ubuntu 22.04.4 LTS system running Linux kernel version 5.15.0-102-generic on an x86\_64 architecture, equipped with an AMD EPYC 9474F 48-Core Processor (2 threads per core) operating at a frequency of 4.11 GHz with 377 GB of RAM. For measuring FLOPs, we utilized the \texttt{perf} profiler tool for Linux (\url{https://perf.wiki.kernel.org}, Last accessed: April 2024) via the CPU-dependent command \texttt{perf stat -e fp\_ret\_sse\_avx\_ops.all \textit{[python-path] [python-script]}}.} when using AD with respect to the input variable. \Cref{fig:flopsAD_experiment} shows the ratios between the gradient computation and the network evaluation using either forward or backward AD. When using forward AD, we observe that the ratio is effectively independent of the output dimension but scales linearly with the input dimension. Conversely, in backward AD, the scaling increases linearly with the output dimension and remains independent of the input dimension. 

For reproducibility, we refer to our GitHub Repository \cite{Uriarte_Optimizing-VPINNs-using-LS_2024}.

\section{Computational costs and implementations}\label{section4}

For simplicity, from now on we restrict to the Deep Fourier Residual (DFR) method \cite{taylor2023deep}, which is a particular case of (R)VPINNs using orthonormal, sinusoidal basis functions so $\mathbf{G}=\mathbf{I}=\mathbf{G}^{-1}$.

\begin{remark}[Spectral viewpoint]\label{Fourier_spectrum}
The Riesz representative of the residual, $r(u^{\boldsymbol{\alpha},\boldsymbol{\omega}})$, and its discrete-level counterpart, $r_M(u^{\boldsymbol{\alpha},\boldsymbol{\omega}})$, relate as follows (recall \Cref{section2.2}): 
\begin{align}
    r(u^{\boldsymbol{\alpha},\boldsymbol{\omega}}) &= \sum_{m = 1}^\infty \{b(u^{\boldsymbol{\alpha},\boldsymbol{\omega}},v_m)-l(v_m)\} v_m,\\
    r_M(u^{\boldsymbol{\alpha},\boldsymbol{\omega}}) &= \sum_{m=1}^M \{b(u^{\boldsymbol{\alpha},\boldsymbol{\omega}},v_m)-l(v_m)\} v_m,
\end{align} where $\{v_m\}_{m\geq 1}$ is a set of orthonormal sinusoidal functions for the test space and $\{b(u^{\boldsymbol{\alpha},\boldsymbol{\omega}},v_m)-l(v_m)\}_{m\geq 1}$ is its set of associated spectral coefficients. We refer to \cite{taylor2023deep} for a detailed discussion concerning the eigenvalues of the underlying residual operator $\mathcal{R}:u\mapsto b(u,\cdot)-l(\cdot)$.

By orthonormality, we obtain
\begin{equation}\label{Fourier_norm_relation}
    \underbrace{\sum_{m = 1}^\infty \{b(u^{\boldsymbol{\alpha},\boldsymbol{\omega}},v_m)-l(v_m)\}^2}_{\Vert r(u^{\boldsymbol{\alpha},\boldsymbol{\omega}}) \Vert_{\mathbb{V}}^2} = \underbrace{\sum_{m=1}^M \{b(u^{\boldsymbol{\alpha},\boldsymbol{\omega}},v_m)-l(v_m)\}^2}_{\Vert r_M(u^{\boldsymbol{\alpha},\boldsymbol{\omega}}) \Vert_{\mathbb{V}}^2  \; = \; \mathcal{L}(\boldsymbol{\alpha},\boldsymbol{\omega})} \; + \; \underbrace{\sum_{m = M+1}^\infty \{b(u^{\boldsymbol{\alpha},\boldsymbol{\omega}},v_m)-l(v_m)\}^2}_{\Vert r(u^{\boldsymbol{\alpha},\boldsymbol{\omega}}) \; - \; r_M(u^{\boldsymbol{\alpha},\boldsymbol{\omega}}) \Vert_{\mathbb{V}}^2},
\end{equation} which shows that the loss function represents a limited portion of $\Vert r(u^{\boldsymbol{\alpha},\boldsymbol{\omega}}) \Vert_{\mathbb{V}}^2$. We will exploit this relation in the numerical experiments.
\end{remark}

\subsection{Conventional GD-based optimizer}\label{section4.1}
The loss function is given by
\begin{align}\label{exact_loss_formula1}
\mathcal{L}(\boldsymbol{\alpha},\boldsymbol{\omega}) = \Vert \mathbf{r}(u^{\boldsymbol{\alpha},\boldsymbol{\omega}})\Vert^2 = \sum_{m=1}^M \left( \int_{\Omega} \nabla u^{\boldsymbol{\alpha},\boldsymbol{\omega}}\cdot \nabla v_m - fv_m \right)^2.
\end{align} 

To approximate the above collection of integrals, we consider a quadrature rule as follows:
\begin{align}\label{loss_formula1}
\mathcal{L}(\boldsymbol{\alpha},\boldsymbol{\omega}) \approx \sum_{m=1}^M \Bigg( \underbrace{\sum_{k=1}^K \gamma_k \Big\lbrace \nabla u^{\boldsymbol{\alpha},\boldsymbol{\omega}}(x_k) \cdot \nabla v_m(x_k) - f(x_k) \cdot v_m(x_k)\Big\rbrace}_{\text{quadrature rule}} \Bigg)^2,
\end{align} where $\gamma_k$ is the quadrature weight associated to the integration point $x_k\in\Omega$ for $1\leq k\leq K$. For simplicity,  we use the same quadrature points and weights in all integrals. Thus, the computational cost for evaluating the loss function is dominated by the following two terms:
\begin{equation}\label{cost1}
\underbrace{\mathcal{O}(KC_{net})}_{\left\lbrace\nabla u^{\boldsymbol{\alpha},\boldsymbol{\omega}}(x_k)\right\rbrace_{k=1}^K} + \underbrace{\mathcal{O}(KM),}_{\left\lbrace\nabla v_m(x_k)\right\rbrace_{k,m=1}^{K,M}}
\end{equation} where $C_{net}>>1$ denotes the cost of evaluating $u^{\boldsymbol{\alpha},\boldsymbol{\omega}}$ over a single input sample $x_k$. Here:
\begin{itemize}
\item $\mathcal{O}(KC_{net})$ refers to the cost of computing $\nabla u^{\boldsymbol{\alpha},\boldsymbol{\omega}}$ over the batch of integration points $\{x_k\}_{k=1}^K$. The cost of computing the spatial gradient of a scalar-valued function, such as $u^{\boldsymbol{\alpha},\boldsymbol{\omega}}$, via AD is comparable to the cost of simply evaluating it (recall \Cref{section3}).

\item $\mathcal{O}(KM)$ indicates the cost of evaluating the collection of $M$ test functions over the batch of integration points $\{x_k\}_{k=1}^K$. We assume that the cost of evaluating the gradient of a single test function at a single integration point is negligible with respect to the cost of evaluating $u^{\boldsymbol{\alpha},\boldsymbol{\omega}}$. Note that $\mathcal{O}(KM)$ also refers to the summation cost in \eqref{loss_formula1} for $1\leq k\leq K$ and $1\leq m\leq M$.
\end{itemize}

Because the computation of the gradient of the loss is comparable to that of the evaluation itself, the final cost of applying a GD-based optimizer at a single iteration results as \eqref{cost1}.

\subsection{The hybrid LS/GD optimizer}\label{section4.2}

Now, the loss has the following form when making explicit the summation over $\boldsymbol{\omega}=[\omega_n]_{n=1}^N$:
\begin{align}\label{loss_formula2}
\mathcal{L}(\boldsymbol{\alpha},\boldsymbol{\omega}) = \Vert \mathbf{B}^{\boldsymbol{\alpha}}\boldsymbol{\omega}-\mathbf{l}\Vert^2 \approx \sum_{m=1}^M \Bigg( \sum_{n=1}^N \omega_n \cdot \underbrace{\sum_{k=1}^K \gamma_k \left\{ \nabla u^{\boldsymbol{\alpha}}_n(x_k)\cdot \nabla v_m(x_k) \right\}}_{\text{entry-wise approx. of } \mathbf{B}^{\boldsymbol{\alpha}}} - \underbrace{\sum_{k=1}^K \gamma_k \left\{ f(x_k) \cdot v_m(x_k)\right\}}_{\text{entry-wise approx. of } \mathbf{l}} \Bigg)^2.
\end{align} 

\begin{remark}[Matrix construction vs. action]\label{loss_evaluations_remark}
In \eqref{loss_formula2}, we first construct the full bilinear-form matrix
\begin{equation}
\mathbf{B}^{\boldsymbol{\alpha}} \approx \left[ \sum_{k=1}^K \gamma_k \left\lbrace \nabla u^{\boldsymbol{\alpha}}_n(x_k)\cdot \nabla v_m(x_k)\right\rbrace \right]_{mn}
\end{equation} and later calculate the action $\mathbf{B}^{\boldsymbol{\alpha}}\boldsymbol{\omega}$. In \eqref{loss_formula1}, instead, we avoid the full matrix construction and directly compute its action.
\end{remark}

The computational cost of the loss evaluation and subsequent optimizer iteration is now dominated by the following three terms:

\begin{equation}\label{cost2}
\underbrace{\mathcal{O}(KC_{net})}_{\substack{\text{computation of} \\[0.5ex] \left\lbrace \nabla u^{\boldsymbol{\alpha}}_n(x_k) \right\rbrace_{k,n=1}^{K,N}}} + \underbrace{\mathcal{O}(KMN)}_{\substack{\text{full matrix} \\[0.5ex] \mathbf{B}^{\boldsymbol{\alpha}} \text{ construction}}} + \;\; \underbrace{\mathcal{O}(MN^2)}_{\text{LS solver}}.
\end{equation} Here:
\begin{itemize}
\item $\mathcal{O}(KC_{net})$ refers to the efficient cost of computing $\nabla\mathbf{u}^{\boldsymbol{\alpha}} = [\nabla u^{\boldsymbol{\alpha}}_n]_{n=1}^N$ over the batch of integration points, which occurs when using \emph{forward-mode AD}---when using backward-mode AD, this cost becomes $\mathcal{O}(KNC_{net})$, recall \Cref{section3}. Note that we have dropped the $d$ multiplicative coefficient with respect to \Cref{table_AD} because its impact is negligible in low-dimensional domains, which is where the scope of this work is focused.

\item Given the trial and test evaluations, $\mathcal{O}(KMN)$ represents the cost of constructing the full bilinear-form matrix $\mathbf{B}^{\boldsymbol{\alpha}}$.

\item Given $\mathbf{B}^{\boldsymbol{\alpha}}$ and $\mathbf{l}$, which have been computed during loss evaluation, $\mathcal{O}(MN^2)$ is the cost of employing a direct LS solver to find $\boldsymbol{\omega}^{\boldsymbol{\alpha}}$.
\end{itemize}

\begin{remark}\label{assumptions_fair_comparison}
For the hybrid LS/GD optimizer to add negligible cost over a conventional GD-based optimizer, we consider the following conditions:
\begin{itemize}
\item $MN < C_{net}$: The multiplication of the number of test and trial spanning functions, $MN$, is smaller than the reference cost of the neural network evaluation, $C_{net}$.  In practical terms, this is achieved by considering a limited set of basis functions $\{v_m\}_{m=1}^M$ for the test space discretization.
\item $N < M \lesssim K$: The number of integration points, $K$, is of the same or greater order as the number of basis functions for the discretization of the test space, $M$, which is greater than the number of spanning functions for the discretization of the trial space, $N$.
\end{itemize} 

When the above conditions are satisfied, the dominant term in both the conventional GD-based and hybrid LS/GD optimizers becomes $\mathcal{O}(KC_{net})$. Thus, the cost of using LS is negligible.

\end{remark}

\subsection{VPINNs and UltraPINNs}\label{section4.3}
 
In our weak-form model problem \eqref{example02}, we have identical trial and test spaces $\U=H^1_0(\Omega)=\V$ equipped with the inner product induced by the bilinear form,
\begin{subequations}
\begin{flalign}\label{weak_implementation}
b(u,v)=\int_\Omega \nabla u\cdot\nabla v, \qquad u\in H^1_0(\Omega) = \U, \quad v\in H^1_0(\Omega) = \V.&&
\end{flalign} In (R)VPINNs, it is only necessary to evaluate the bilinear form over a finite set of test functions. In particular, if the test functions are smooth, as in the case of the DFR method, integration by parts of \eqref{weak_implementation} yields
\begin{flalign}\label{ultraweak_implementation}
b(u,v_m)=\int_{\Omega} -u\cdot\Delta v_m,\qquad u\in H^1_0(\Omega) = \U, \quad v_m\in \V_M.&&
\end{flalign}
\end{subequations} Effectively, \eqref{weak_implementation} may be read as an alternative integration rule of \eqref{ultraweak_implementation} that avoids the need to differentiate the trial function, without modifying the underlying (weak) variational formulation of the PDE itself. We call the latter an \emph{ultraweak implementation} and refer to the resulting (R)VPINN approach as an \emph{UltraPINN}. As a result, we obtain a faster construction of matrix $\mathbf{B}^{\boldsymbol{\alpha}}$ when using the hybrid LS/GD optimizer in UltraPINNs. Thus, the computational costs are those of \eqref{cost2} with a smaller constant in the second term.

For simplicity, VPINNs will henceforth refer only to the typical weak implementation approach of RVPINNs. 

\subsection{Implementation aspects}\label{section4.4}

We use \emph{Keras 3.3} \cite{chollet2015keras} with the backend of \emph{Tensorflow (TF) 2.16} \cite{tensorflow2015-whitepaper} to implement the workflow of our experiments in \emph{Python 3.12}. We prioritize the use of Keras to construct our models while leaving TF for tensor-specific tasks such as AD. For backward AD, we employ the \texttt{tf.GradientTape} environment and deactivate the \texttt{watch\_accessed\_variables} option to purposely control the traced variables during the forward pass. For forward AD, we use \texttt{tf.autodiff.ForwardAccumulator}.

For optimization, we redefine the base \texttt{train\_step} function of the master class \texttt{keras.Model}. In the conventional GD-based optimization, we evaluate the loss function within the \texttt{tf.GradientTape} environment tracing all trainable parameters, $(\boldsymbol{\alpha},\boldsymbol{\omega}) \mapsto \mathcal{L}(\boldsymbol{\alpha},\boldsymbol{\omega})$. Then, we extract the gradient to subsequently apply an iteration of our chosen GD-based optimizer. In the hybrid LS/GD optimizer, we only trace $\boldsymbol{\alpha}$. However, the way in which we implement this tracing (i.e., use the \texttt{tf.GradientTape} environment) is critical in terms of performance. Indeed, we explore the following two implementation approaches after constructing $\mathbf{B}^{\boldsymbol{\alpha}}$ and $\mathbf{l}$ (Step 1) and finding $\boldsymbol{\omega}^{\boldsymbol{\alpha}}$ via the \texttt{tf.linalg.lstsq} solver (Step 2):
\begin{itemize}
\item[(a)] We compute $\Vert \mathbf{B}^{\boldsymbol{\alpha}} \boldsymbol{\omega}^{\boldsymbol{\alpha}}-\mathbf{l}\Vert^2$ to obtain the final loss evaluation (Step 3a). Here, the tracing of $\boldsymbol{\alpha}$ is made throughout the entire process (from Step 1 to Step 3a), constructing a computable graph $\boldsymbol{\alpha}\mapsto\boldsymbol{\omega}^{\boldsymbol{\alpha}}\mapsto\Vert \mathbf{B}^{\boldsymbol{\alpha}} \boldsymbol{\omega}^{\boldsymbol{\alpha}}-\mathbf{l}\Vert^2$ that includes the LS solver.

\item[(b)] We (re-)evaluate $\mathcal{L}(\boldsymbol{\alpha}, \boldsymbol{\omega}^{\boldsymbol{\alpha}})$ using the implementation that avoids the construction of the bilinear-form matrix (see \Cref{loss_evaluations_remark}), establishing the tracing environment at this step only (Step 3b). As a result, the LS solver is untraced, and the resulting computational graph $\boldsymbol{\alpha}\mapsto \mathcal{L}(\boldsymbol{\alpha},\boldsymbol{\omega}^{\boldsymbol{\alpha}})$ never reads $\boldsymbol{\omega}^{\boldsymbol{\alpha}}$ as $\boldsymbol{\alpha}$-dependent but as a given constant.
\end{itemize}

\begin{remark}[Revisiting total vs. partial derivative viewpoint]
Approaches (a) and (b) might be interpreted as computations of the total and partial derivatives of $\mathcal{L}$ with respect to $\boldsymbol{\alpha}$, respectively (recall \Cref{total_partial_derivative}). Indeed, while approach (a) records the full dependence $\boldsymbol{\alpha}\mapsto\boldsymbol{\omega}^{\boldsymbol{\alpha}}$ that entails a significant number of operations, approach (b) simply ignores this. Because both produce the same effect (total and partial derivatives coincide), the lack of tracing in the latter results in a significant memory-saving advantage.
\end{remark}

We employ approach (b) during implementation (see \hyperref[train_step_guidelines]{Algorithm~\ref{train_step_guidelines}}).  For further details, see \cite{Uriarte_Optimizing-VPINNs-using-LS_2024}.

\begin{algorithm}[H]
\SetAlgoLined
Compute $\mathbf{B}^{\boldsymbol{\alpha}}$ and $\mathbf{l}$ using forward-mode AD (\texttt{tf.autodiff.ForwardAccumulator}).\\
Compute $\boldsymbol{\omega}^{\boldsymbol{\alpha}} = \texttt{tf.linalg.lstsq}(\mathbf{B}^{\boldsymbol{\alpha}}, \mathbf{l})$.\\
Assign $\boldsymbol{\omega}^{\boldsymbol{\alpha}}$ to the coefficient values of the output layer in the neural network.\\
Within the backward-mode AD environment (\texttt{tf.GradientTape}) and watching/tracing only $\boldsymbol{\alpha}$, \hspace{1cm} evaluate $\mathcal{L}(\boldsymbol{\alpha}, \boldsymbol{\omega}^{\boldsymbol{\alpha}})$ without constructing $\mathbf{B}^{\boldsymbol{\alpha}}$ and $\mathbf{l}$.\\
Compute $\partial\mathcal{L}/\partial\boldsymbol{\alpha}$ and apply the GD optimizer to update $\boldsymbol{\alpha}$.
\caption{Implementation guideline for \texttt{train\_step} using TensorFlow}\label{train_step_guidelines}
\end{algorithm}

\subsection{Illustrative experiments}\label{section4.5}

To illustrate the scaling behavior of the computational time for the previously outlined implementations, we consider a comparison of the various methods. We highlight that these results are independent of the particular PDE considered and hyperparameters of the optimizer, although, of course, the quality of the obtained solutions will be highly dependent upon them. Our focus here is on the computational cost.

We consider the DFR discretization setup given by $M = 5N$, $K =10N$, and $N\in\{2,2^2,\ldots,2^9\}$ with the same neural network architecture as in \Cref{fig:flopsAD_experiment} so that it meets the conditions of \Cref{assumptions_fair_comparison}. \Cref{LSGDvsGD_ratios} displays the quotients when measuring FLOPs of a single-iteration optimization under the following conditions:
\begin{itemize}
    \item \textbf{Numerator.} Our hybrid LS/GD optimizer in UltraPINNs and VPINNs with forward or backward AD for the computation of $\{\nabla\mathbf{u}^{\boldsymbol{\alpha}}(x_k)\}_{k=1}^K$.
    \item \textbf{Denominator.} A conventional GD-based optimization in VPINNs with backward AD for the computation of $\{\nabla u^{\alpha,\omega}(x_k)\}_{k=1}^K$.
\end{itemize}

\begin{table}[H]
\centering
\begin{tabular}{|c|ccccccccc|}
\toprule
$N$                                               & 2    & 4    & 8    & 16   & 32   & 64    & 128   & 256   & 512  \\ \midrule \midrule
UltraPINNs                                                     & 0.6702                 & 0.6699                 & 0.6699                 & 0.6699                  & 0.6707                  & 0.6727                  & 0.6751                   & 0.6884                   & 0.7244                   \\ \midrule
\begin{tabular}[c]{@{}c@{}}VPINNs \\ (forward AD)\end{tabular} & 1.9904                 & 1.9903                 & 1.9976                 & 1.9927                  & 2.0000                  & 2.0000                  & 1.9970                   & 2.0217                   & 2.0476                   \\ \midrule
\begin{tabular}[c]{@{}c@{}}VPINNs\\ (backward AD)\end{tabular} & 1.4904                 & 1.8252                 & 2.4939                 & 3.8264                  & 6.4634                  & 11.8182                 & 22.3547                  & 43.7681                  & 86.0544                  \\ \bottomrule
\end{tabular}
\caption{FLOP cost quotients when applying a single iteration of the hybrid LS/GD optimizer in UltraPINNs and VPINNs versus a conventional GD optimizer in VPINNs. The ``forward AD'' and ``backward AD'' labels in VPINNs refer only to the computation of the derivatives of the network with respect to the input variable. During training, the GD update with respect to the trainable parameters is always performed using backward AD.}
\label{LSGDvsGD_ratios}
\end{table}

\hyperref[GD_measurements]{Tables~\ref{GD_measurements}} and \ref{LSGD_measurements} show absolute FLOP measurements of the quotients shown above. For further details, see our GitHub Repository \cite{Uriarte_Optimizing-VPINNs-using-LS_2024}.

As we observe, using the hybrid LS/GD optimizer in UltraPINNs is cheaper than using a conventional GD-based optimizer in VPINNs with forward AD due to the absence of spatial differentiation. On the other hand, when using the hybrid LS/GD optimizer in VPINNs, the scaling ratio is moderate and constant when using forward AD for spatial differentiation (by a factor of two, approximately) while growing linearly with the value of $N$ when using backward AD. Consequently, we avoid using VPINNs with backward AD from now on, and VPINNs---without further specifications---will refer to the case with spatial forward AD.

\begin{table}[H]
\centering
\begin{tabular}{|cl|cc|ccc|}
\toprule
\multicolumn{2}{|c|}{\textbf{Specifications}} & \multicolumn{2}{c|}{\textbf{FLOPs}} & \multicolumn{3}{c|}{\textbf{Ratios}} \\  \midrule
\textit{Implementation} & $N$ & \textit{(vect-)net } & \textit{total} & $C_{net} = \frac{\textit{net}}{K}$ & $\frac{\textit{(vect-)net}}{KC_{net}}$ & $\frac{total}{KC_{net}}$   \\ \midrule \midrule
\multicolumn{1}{|c}{\multirow{9}{*}{\begin{tabular}[c]{@{}c@{}}VPINNs\\ (backward AD)\\ + \\ GD optimizer \\ (backward AD)\end{tabular}}} & 2 & 1.72e+08 & 1.04e+09 & 8.58e+06 & 1.00e+00 & 6.09e+00\\
\multicolumn{1}{|l}{}                                                                                   & 4 & 3.43e+08 & 2.06e+09 & 8.58e+06 & 1.00e+00 & 6.00e+00\\
\multicolumn{1}{|l}{}                                                                                   & 8 & 6.87e+08 & 4.09e+09 & 8.59e+06 & 1.00e+00 & 5.96e+00\\
\multicolumn{1}{|l}{}                                                                                   & 16 & 1.38e+09 &  8.18e+09 & 8.61e+06 & 1.00e+00 & 5.94e+00\\
\multicolumn{1}{|l}{}                                                                                   & 32 & 2.77e+09 & 1.64e+10 & 8.64e+06 & 1.00e+00 & 5.93e+00\\
\multicolumn{1}{|l}{}                                                                                   & 64 & 5.57e+09 & 3.30e+10 & 8.70e+06 & 1.00e+00 & 5.93e+00\\
\multicolumn{1}{|l}{}                                                                                   & 128 & 1.13e+10 & 6.71e+10 & 8.84e+06 & 1.00e+00 & 5.93e+00\\
\multicolumn{1}{|l}{}                                                                                   & 256 & 2.33e+10 & 1.38e+11 & 9.11e+06 & 1.00e+00 & 5.94e+00\\
\multicolumn{1}{|l}{}                                                                                   & 512 & 4.94e+10 &  2.94e+11 & 9.64e+06 & 1.00e+00 & 5.95e+00\\ \midrule %\bottomrule
\end{tabular}%
%}
\caption{Reference FLOP measurements when performing a single iteration of a conventional GD-based optimizer in VPINNs. \textit{net} is the cost of evaluating $\{u^{\boldsymbol{\alpha},\boldsymbol{\omega}}(x_k)\}_{k=1}^K$. $\textit{total}$ is the final optimization cost that includes the spatial gradient computation, the evaluation of the loss evaluation, and the application of the GD optimizer. The last three columns show the empirical constants regarding the cost order discussion in \Cref{section4.2}.}
\label{GD_measurements}
\vskip 1em
\begin{tabular}{|cl|cc|ccc|}
\midrule
%\toprule
%\multicolumn{6}{|c|}{\textbf{Reference cost when using a GD optimizer in VPINNs}} \\ \toprule
%\multicolumn{2}{|c|}{\textbf{Specifications}} & \multicolumn{2}{c}{\textbf{FLOPs}} & \multicolumn{2}{c|}{\textbf{Ratios}} \\  \midrule
%\textit{Implementation} & $N$ & \textit{net } & \textit{total} & $C_{net} = \frac{\textit{net}}{K}$ & $\frac{total}{KC_{net}}$   \\ \midrule \midrule
\multicolumn{1}{|c}{\multirow{9}{*}{\begin{tabular}[c]{@{}c@{}}UltraPINNs\\ (no AD)\\ + \\ LS/GD optimizer \\ (backward AD)\end{tabular}}} & 2  &1.72e+08 & 6.97e+08 & 8.58e+06 & 1.00e+00 & 4.06e+00 \\
                                                                                   & 4  & 3.43e+08 & 1.38e+09 & 8.58e+06 & 1.00e+00 & 4.01e+00 \\
                                                                                   & 8  & 6.87e+08 & 2.74e+09 & 8.59e+06 & 1.00e+00 & 3.99e+00 \\
                                                                                   & 16 & 1.38e+09 & 5.48e+09 & 8.61e+06 & 1.00e+00 & 3.98e+00 \\
                                                                                   & 32 & 2.77e+09 & 1.10e+10 & 8.64e+06 & 1.00e+00 & 3.97e+00 \\
                                                                                   & 64 & 5.57e+09 & 2.22e+10 & 8.70e+06 & 1.00e+00 & 3.98e+00 \\
                                                                                   & 128 & 1.13e+10 & 4.53e+10 & 8.84e+06 & 1.00e+00 & 4.00e+00 \\
                                                                                   & 256 & 2.33e+10 & 9.50e+10 & 9.11e+06 & 1.00e+00 & 4.07e+00 \\
                                                                                   & 512 & 4.94e+10 & 2.13e+11 & 9.64e+06 & 1.00e+00 & 4.32e+00 \\ \midrule
\multicolumn{1}{|c}{\multirow{9}{*}{\begin{tabular}[c]{@{}c@{}}VPINNs\\ (forward AD)\\ + \\ LS/GD optimizer \\ (backward AD)\end{tabular}}}  & 2 & 5.09e+08 & 2.07e+09 & 8.58e+06 & 2.97e+00 & 1.21e+01 \\
                                                                                   & 4 & 1.02e+09 & 4.10e+09 & 8.58e+06 & 2.97e+00 & 1.19e+01 \\
                                                                                   & 8 & 2.04e+09 & 8.17e+09 & 8.59e+06 & 2.97e+00 & 1.19e+01 \\
                                                                                   & 16 & 4.09e+09 & 1.63e+10 & 8.61e+06 & 2.97e+00 & 1.19e+01 \\
                                                                                   & 32 & 8.21e+09 & 3.28e+10 & 8.64e+06 & 2.97e+00 & 1.18e+01 \\
                                                                                   & 64 & 1.65e+10 & 6.60e+10 & 8.70e+06 & 2.97e+00 & 1.19e+01 \\
                                                                                   & 128 & 3.36e+10 & 1.34e+11 & 8.84e+06 & 2.97e+00 & 1.19e+01 \\
                                                                                   & 256 & 6.92e+10 & 2.79e+11 & 9.11e+06 & 2.97e+00 & 1.19e+01 \\
                                                                                   & 512 & 1.47e+11 & 6.02e+11 & 9.64e+06 & 2.97e+00 & 1.22e+01 \\  \midrule
\multicolumn{1}{|c}{\multirow{9}{*}{\begin{tabular}[c]{@{}c@{}}VPINNs\\ (backward AD)\\ + \\ LS/GD optimizer \\ (backward AD)\end{tabular}}} & 2 & 5.09e+08 & 1.55e+09 & 8.58e+06 & 2.97e+00 & 9.06e+00 \\
                                                                                   & 4 & 1.70e+09 & 3.76e+09 & 8.58e+06 & 4.94e+00 & 1.09e+01 \\
                                                                                   & 8 & 6.10e+09 & 1.02e+10 & 8.59e+06 & 8.88e+00 & 1.48e+01 \\
                                                                                   & 16 & 2.31e+10 & 3.13e+10 & 8.61e+06 & 1.67e+01 & 2.27e+01 \\
                                                                                   & 32 & 8.99e+10 & 1.06e+11 & 8.64e+06 & 3.25e+01 & 3.84e+01 \\
                                                                                   & 64 & 3.57e+11 & 3.90e+11 & 8.70e+06 & 6.40e+01 & 7.00e+01 \\
                                                                                   & 128 & 1.44e+12 & 1.50e+12 & 8.84e+06 & 1.27e+02 & 1.33e+02 \\ 
                                                                                   & 256 & 5.90e+12 & 6.04e+12 & 9.11e+06 & 2.53e+02 & 2.59e+02 \\    
                                                                                   & 512 & 2.49e+13 & 2.53e+13 & 9.64e+06 & 5.05e+02 & 5.12e+02 \\ \bottomrule
\end{tabular}%
%}
\caption{FLOP measurements when performing a single iteration of our LS/GD optimizer in UltraPINNs and VPINNs when using different AD modes for spatial differentiation. \textit{vect-net} is the cost of evaluating $\{\mathbf{u}^{\boldsymbol{\alpha},\boldsymbol{\omega}}(x_k)\}_{k=1}^K$ in UltraPINNs or $\{\nabla\mathbf{u}^{\boldsymbol{\alpha},\boldsymbol{\omega}}(x_k)\}_{k=1}^K$ in VPINNs. \textit{total} is the final cost of applying the LS/GD optimizer, which includes the LS matrix construction and solving, (re-)evaluation of the loss function, and GD update. The last three columns show the empirical constants regarding the cost order discussion in \Cref{section4.2}.}
\label{LSGD_measurements}
\end{table}

\subsection{Limitations}

The proposed method has several limitations that go beyond the scope of this work, but we briefly mention them here for completeness. First, the inversion of the Gram matrix may significantly increase computational costs, and while sparse cases like FEM discretizations could use local factorizations and forward-backward substitution to avoid explicit inversion, this remains challenging with current neural network libraries. Second, stochastic numerical integration, employed during training to prevent overfitting, often demands a large number of points to reduce integration errors, resulting in significant computational overhead. Additionally, these errors can severely impact the accuracy of the LS solver, emphasizing the importance of effectively managing stochastic integration errors.

\section{Numerical experiments}\label{section5}

We conduct numerical experiments in both one- and two-dimensional domains to demonstrate the superiority of the hybrid LS/GD optimizer compared to a conventional GD-based optimizer in UltraPINNs and VPINNs. \Cref{section5.1} outlines the experimental setup, while \hyperref[section5.2]{Sections~\ref{section5.2}}, \ref{section5.3}, \ref{section5.4}, and \ref{section5.5} shows results on various problems with distinct properties of interest, such as high-frequency and singular behaviors.

\subsection{Experimentation setup}\label{section5.1}

We consider Poisson's weak-form equation \eqref{example02} for $\Omega=(0,\pi)^d$, $d\in\{1,2\}$. Because the inner product in $H^1_0(\Omega)=\U=\V$ is chosen as the bilinear form, the error function $e(u)=u-u^*\in\U$ coincides with the Riesz representative of the residual $r(u)\in \V$ for all $u\in \U$.

We construct a set of orthonormal functions in $H^1_0(\Omega)$ following the DFR scheme \cite{taylor2023deep}. Specifically, for $d=1$, we consider
\begin{equation}\label{eq:basis1D}
	v_m := \frac{\tilde{v}_m}{\Vert\tilde{v}_m\Vert_{\mathbb{V}}},\qquad \tilde{v}_m(x):=\sin\left(m x \right), \qquad m\geq 1.
\end{equation} For $d=2$, we construct a tensor-product structure identifying $m=(m_1,m_2)$ and choosing
\begin{equation}\label{eq:basis2D}
\tilde{v}_m(x,y) := \sin(m_1 x)\sin(m_2 y),\qquad  m_1,m_2\geq 1.
\end{equation} In this way, for a given cut-off frequency $M(=M_1 \times M_2\text { in } 2D)$, we obtain our orthonormal discretization for the test space. To implement the derivatives of the test functions ($\nabla v_m$ in VPINNs and $\Delta v_m$ in UltraPINNs), we consider their usual trigonometric closed formulas to avoid applying AD on them during training.

The architecture of the neural network for our experiments consists of three hidden layers of width $N$ with the hyperbolic tangent as the activation function. We consider $\xi_D(x)=x(x-\pi)$ (in 1D) and $\xi_D(x,y)=xy(x-\pi)(y-\pi)$ (in 2D) as the cut-off functions to strongly impose the homogeneous boundary conditions in the network (recall \Cref{imposition_Dirichlet}). Moreover, we use Adam \cite{kingma2014adam} as the conventional GD-based optimizer with an initial learning rate of $10^{-3}$. The goal of this section is to show the effectiveness of using LS followed by Adam (the LS/Adam optimizer) versus Adam alone. To ensure a fair comparison, we fix the total number of training iterations for each experiment and use two independent neural networks with identical architectures and initializations at the start of training. While hyperparameter tuning is essential for achieving superior accuracy, it is not the focus of this study.

We employ a stochastic numerical integration scheme that randomly samples a batch of integration points and weights at each training iteration. It is founded on a stochastic and unbiased composite quadrature rule tailored to be exact for piecewise-linear functions (see details in \ref{appendixIntegration}). In addition, we consider a validation dataset with a larger number of integration points and weights to measure error progress during training. 

All our experiments consist of modifications of the given source $f$ in order to explore the training behavior of the optimizers when approximating manufactured solutions $u^*$ that present distinct properties of interest.

\subsection{A smooth 1D problem}\label{section5.2}

Let
\begin{equation}\label{exact1}
	u^*(x)= \sin(4x)\sin\left( \frac{x}{2} \right),\qquad 0\leq x\leq \pi,
\end{equation} and let us consider the following discretization setup: $N=16$ for the trial space, $M=2N = 32$ for the test space, $K=32M = 1\mathord{,}024$ for numerical integration over a uniform integration grid, and $8K = 8\mathord{,}192$ for the size of the validation set.

\begin{figure}[htbp]
		\centering
        \begin{subfigure}[b]{0.49\textwidth}
		\includegraphics[scale=0.49]{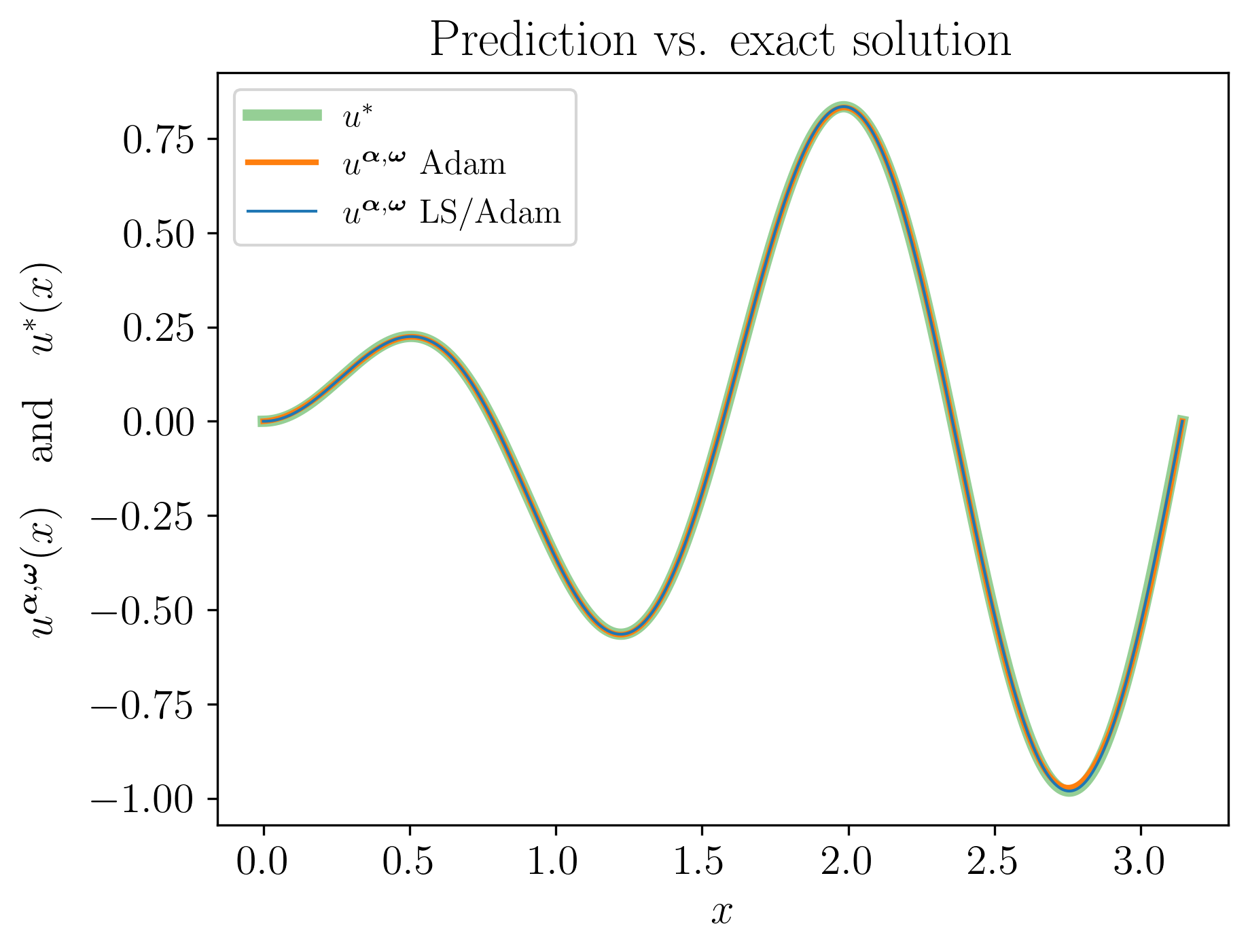}
        \end{subfigure}
		\begin{subfigure}[b]{0.49\textwidth}
		\includegraphics[scale=0.49]{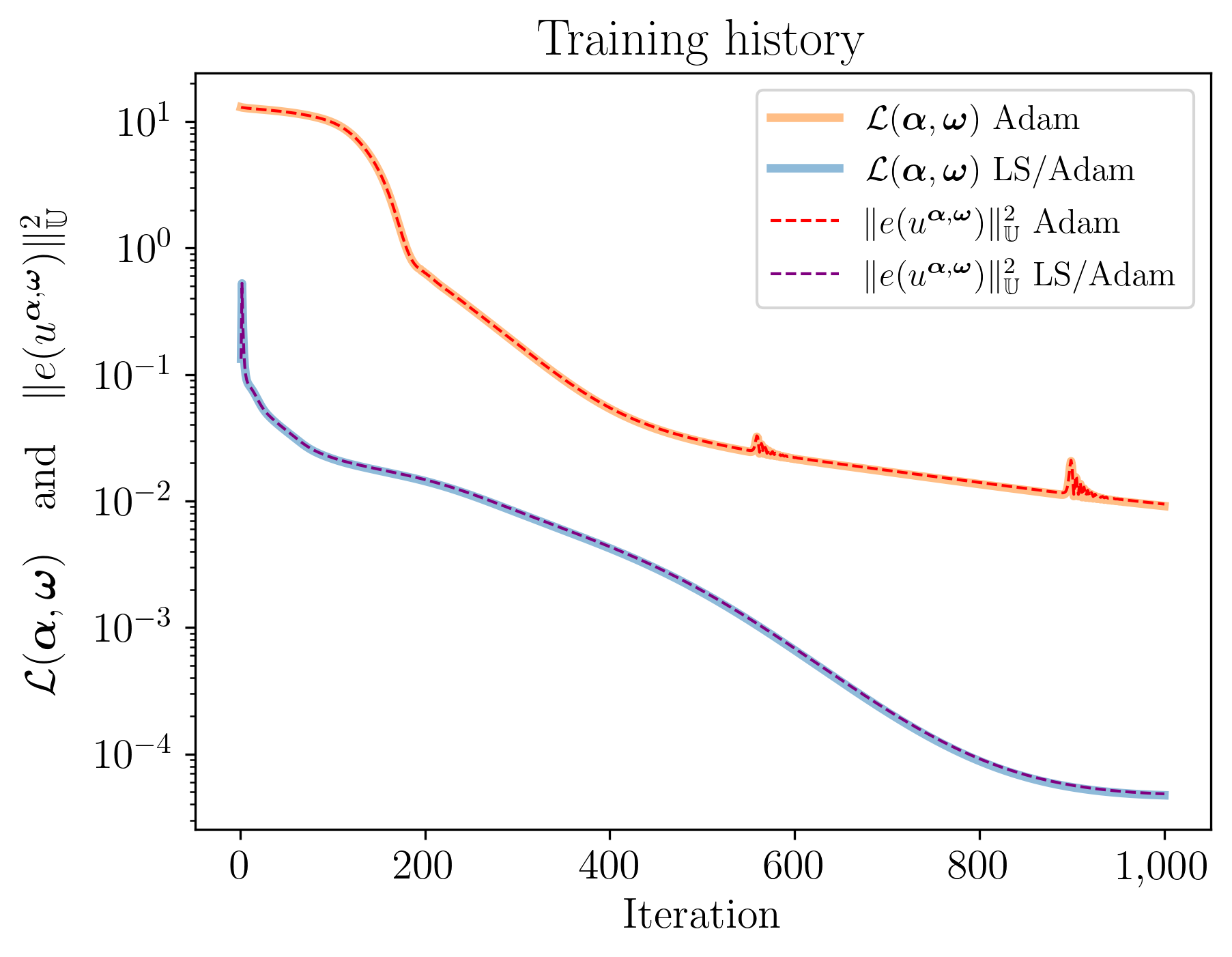}
        \end{subfigure}
        \begin{subfigure}[b]{0.49\textwidth}
		\includegraphics[scale=0.49]{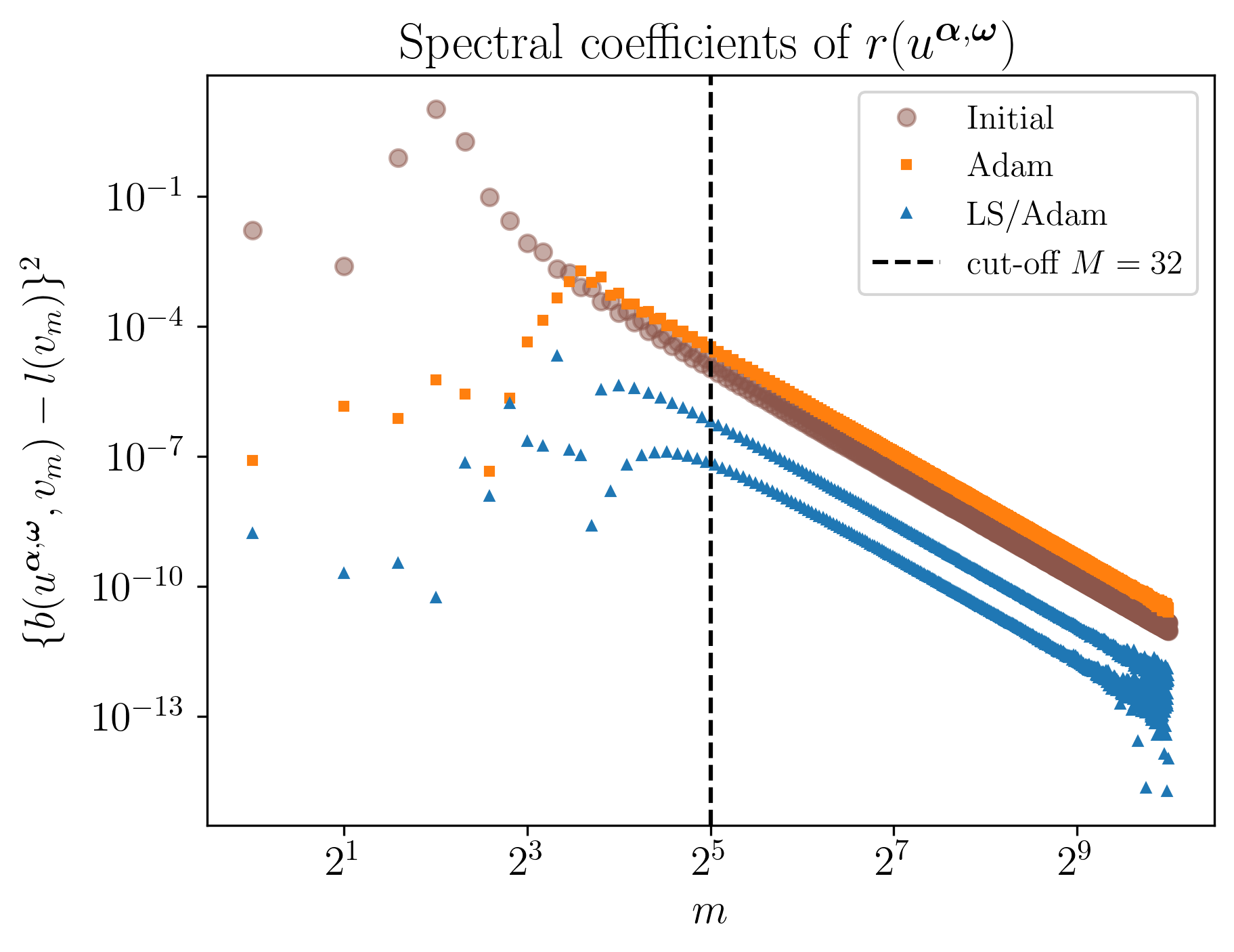}
		\end{subfigure}
  		\begin{subfigure}[b]{0.49\textwidth}
		\includegraphics[scale=0.49]{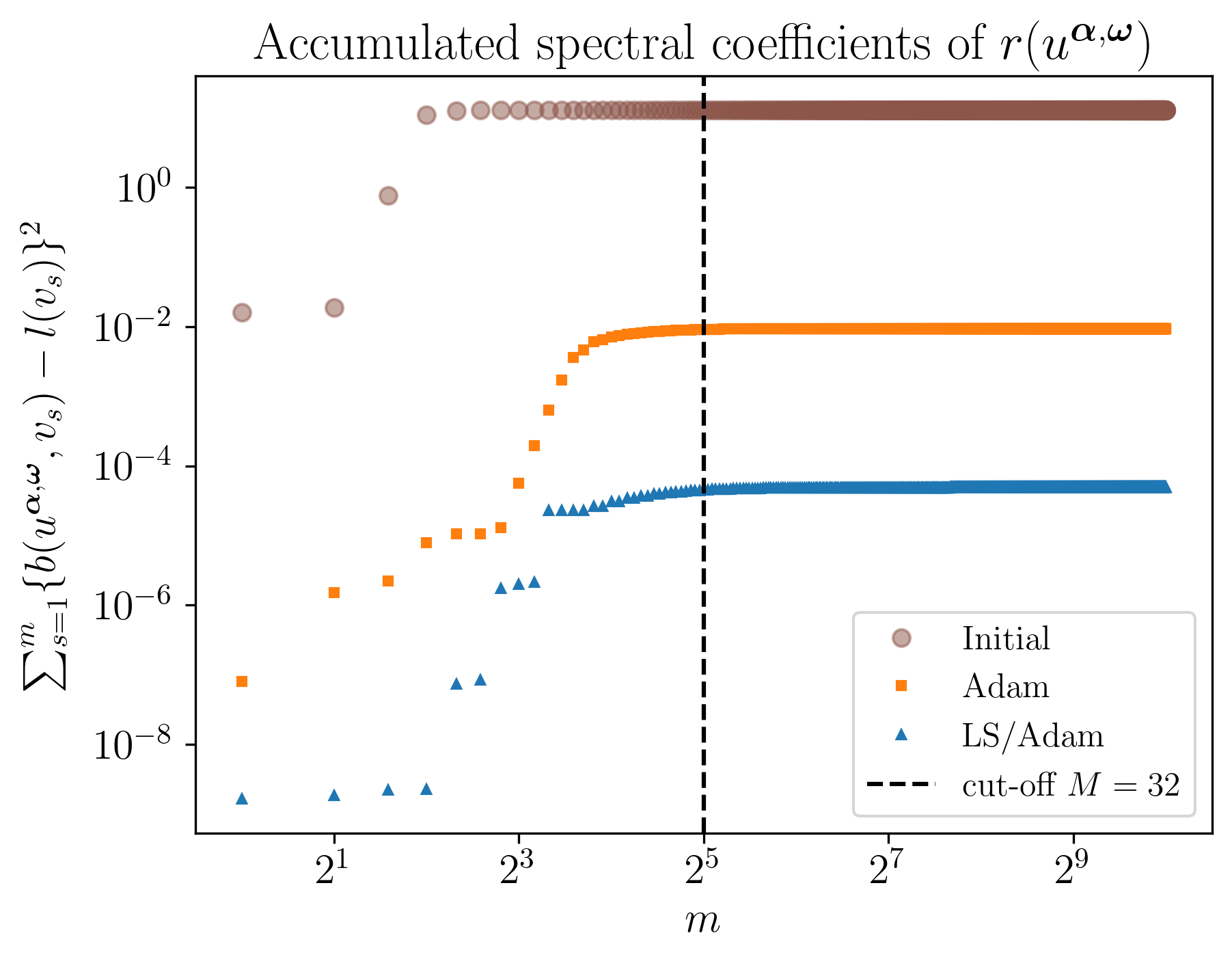}
		\end{subfigure}
		\caption{Results for VPINNs with manufactured solution $u^*(x)=\sin(4x)\sin(x/2)$ using either Adam or LS/Adam. The final relative errors\protect\footnotemark \; are $2.72\%$ and $0.19\%$ for Adam and LS/Adam, respectively.
  }
	\label{fig:exact1weak}
\end{figure}

\begin{figure}[htbp]
		\centering
        \begin{subfigure}[b]{0.49\textwidth}
		\includegraphics[scale=0.49]{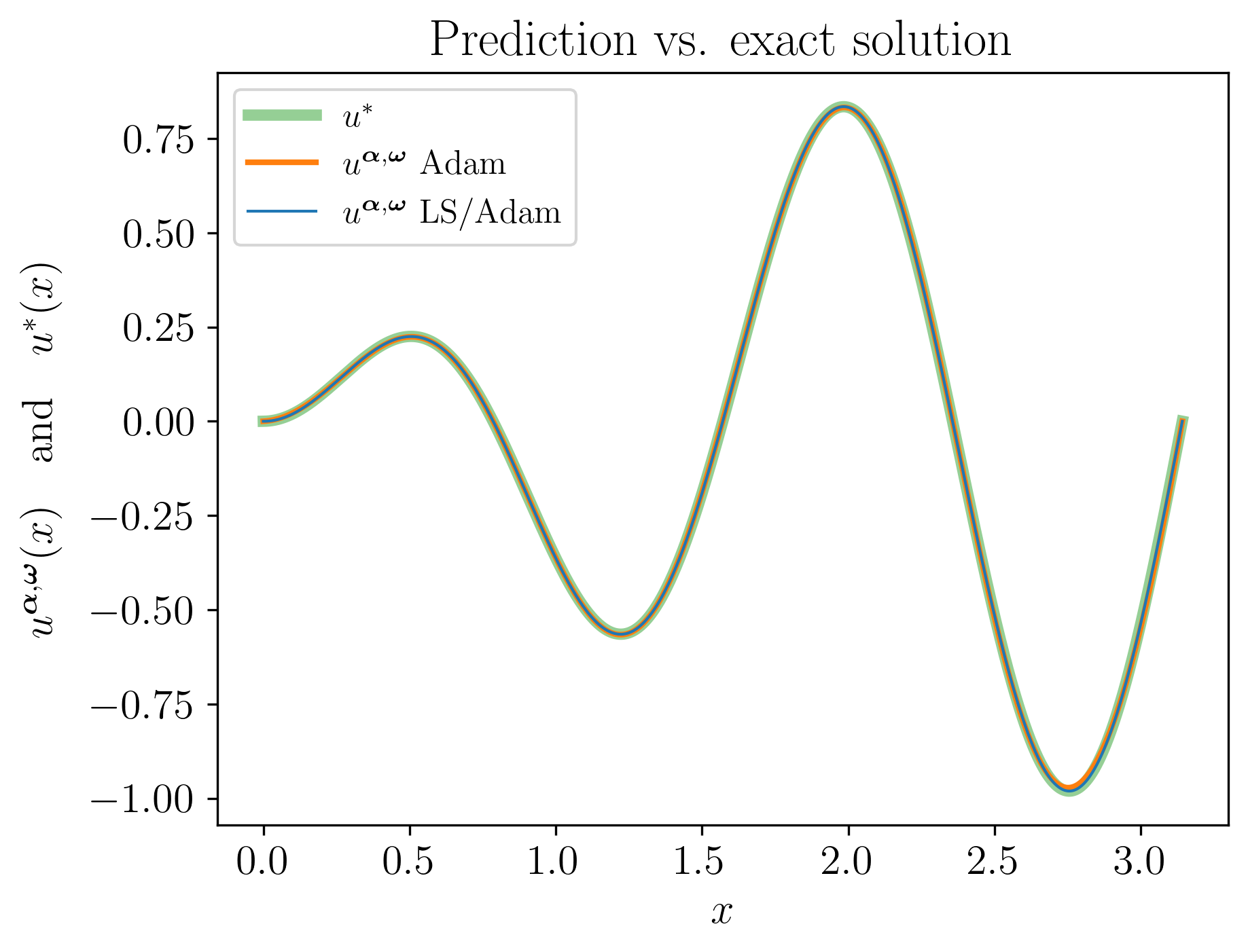}
        \end{subfigure}
		\begin{subfigure}[b]{0.49\textwidth}
		\includegraphics[scale=0.49]{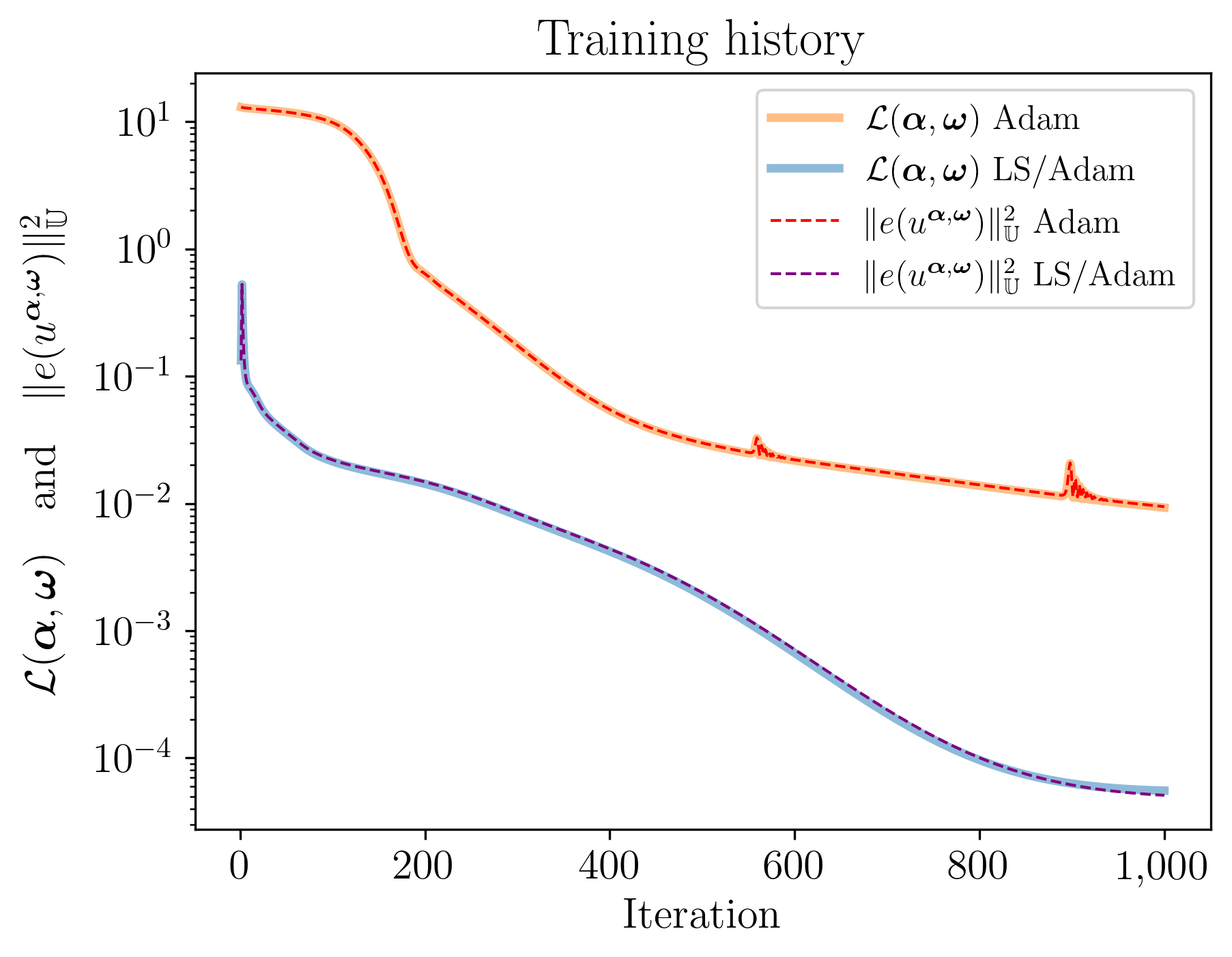}
        \end{subfigure}
        \begin{subfigure}[b]{0.49\textwidth}
		\includegraphics[scale=0.49]{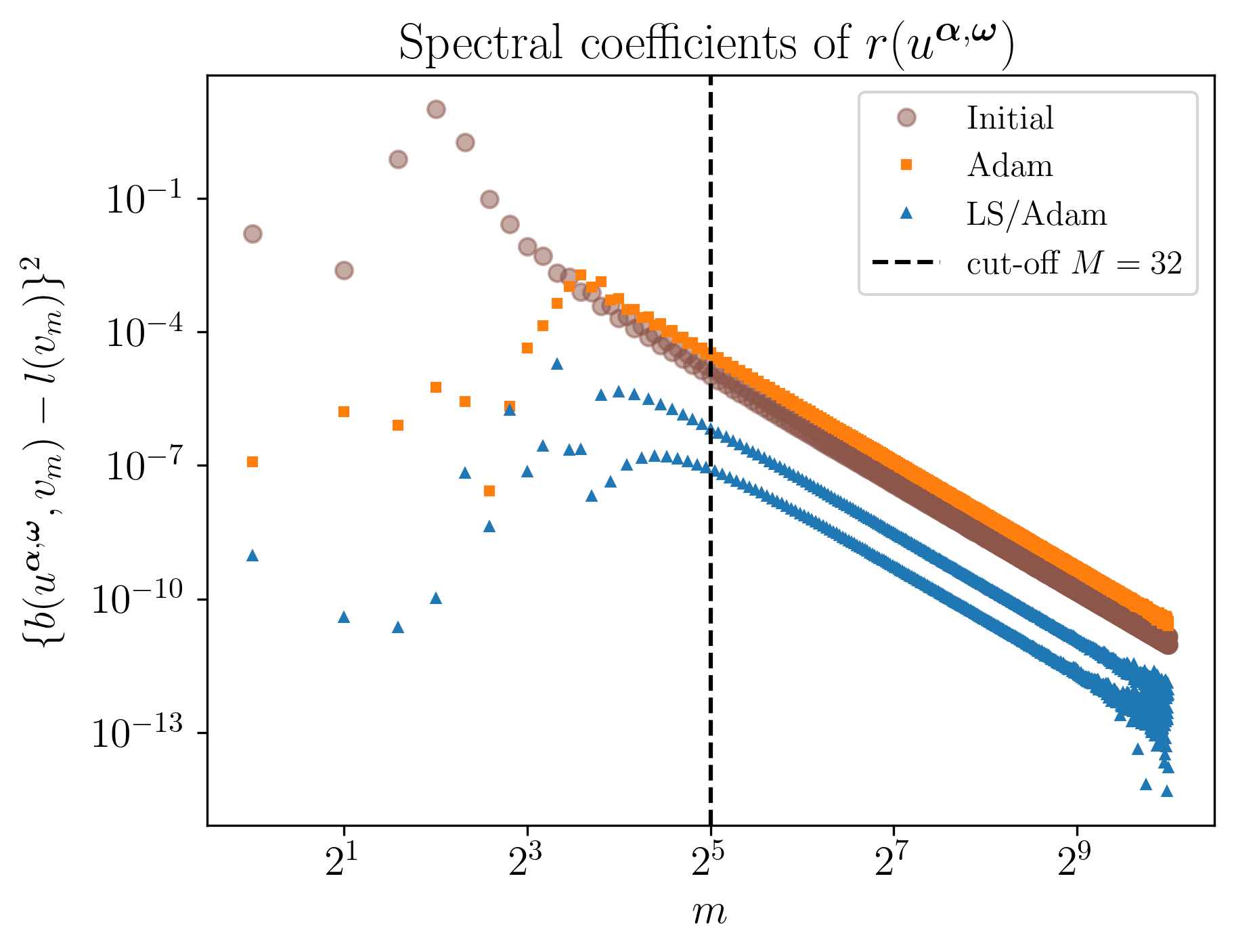}
		\end{subfigure}
  		\begin{subfigure}[b]{0.49\textwidth}
		\includegraphics[scale=0.49]{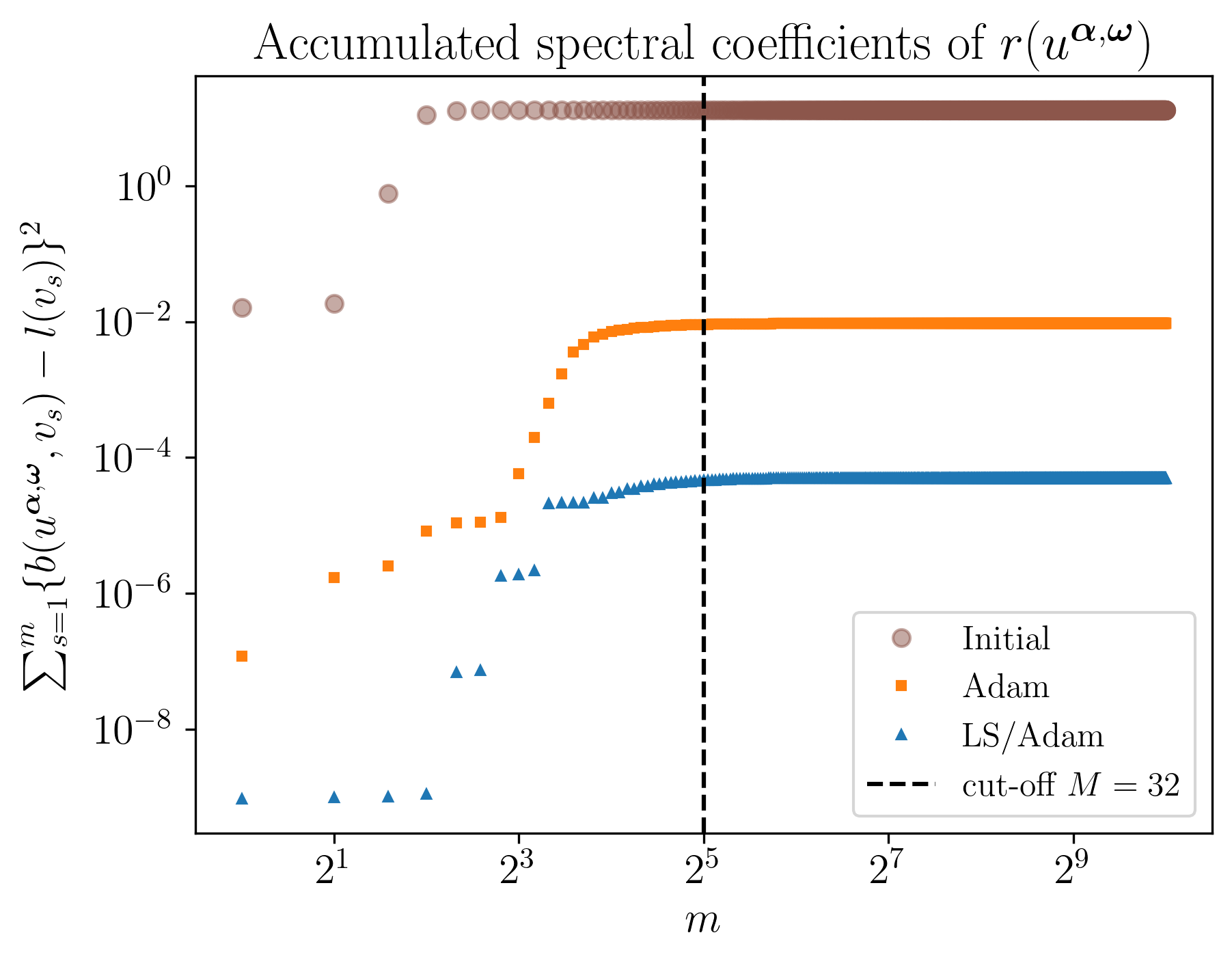}
		\end{subfigure}
\caption{Results for UltraPINNs with manufactured solution $u^*(x)=\sin(4x)\sin(x/2)$ using either Adam or LS/Adam. The resulting final relative errors are $2.72\%$ and $0.20\%$ for Adam and LS/Adam.
}
	\label{fig:exact1ultraweak}
\end{figure}

\hyperref[fig:exact1weak]{Figures~\ref{fig:exact1weak}} and \ref{fig:exact1ultraweak} show the results of VPINNs and UltraPINNs, respectively, after $1\mathord{,}000$ training iterations when using either Adam or LS/Adam.

\footnotetext{By \emph{relative error} we refer to the quotient $\Vert e(u^{\boldsymbol{\alpha},\boldsymbol{\omega}})\Vert_{\mathbb{U}}/\Vert u^*\Vert_{\mathbb{U}}$. The numerator is computed using the validation set of integration points while the denominator is calculated exactly. We highlight the difference between the squared norm of the error displayed in the training history panel and this square-rooted version.}

The results in terms of accuracy during training are virtually identical in VPINNs and UltraPINNs. Moreover, the LS/Adam optimizer significantly boosts the decrease of the loss versus Adam alone, driving an entire order of magnitude improvement from the very first iteration and an almost three-order improvement at the end of training. The correlation between $\mathcal{L}(\boldsymbol{\alpha},\boldsymbol{\omega})$ and $\Vert e(u^{\boldsymbol{\alpha},\boldsymbol{\omega}})\Vert_{\mathbb{U}}^2$ is practically one-to-one during training. This occurs because the contributions of the spectral coefficients (recall \Cref{Fourier_spectrum}) after the test space discretization cut-off are negligible (equivalently, the curves of the accumulated spectral coefficients plateau to the right of the cut-off). In other words, our discretization of the test space adequately captures the significant spectrum of the underlying residual operator.

\subsection{A smooth 1D problem with high frequencies}\label{section5.3}

\begin{figure}[htbp]
	\centering
        \begin{subfigure}[b]{0.49\textwidth}
	    \includegraphics[scale=0.49]{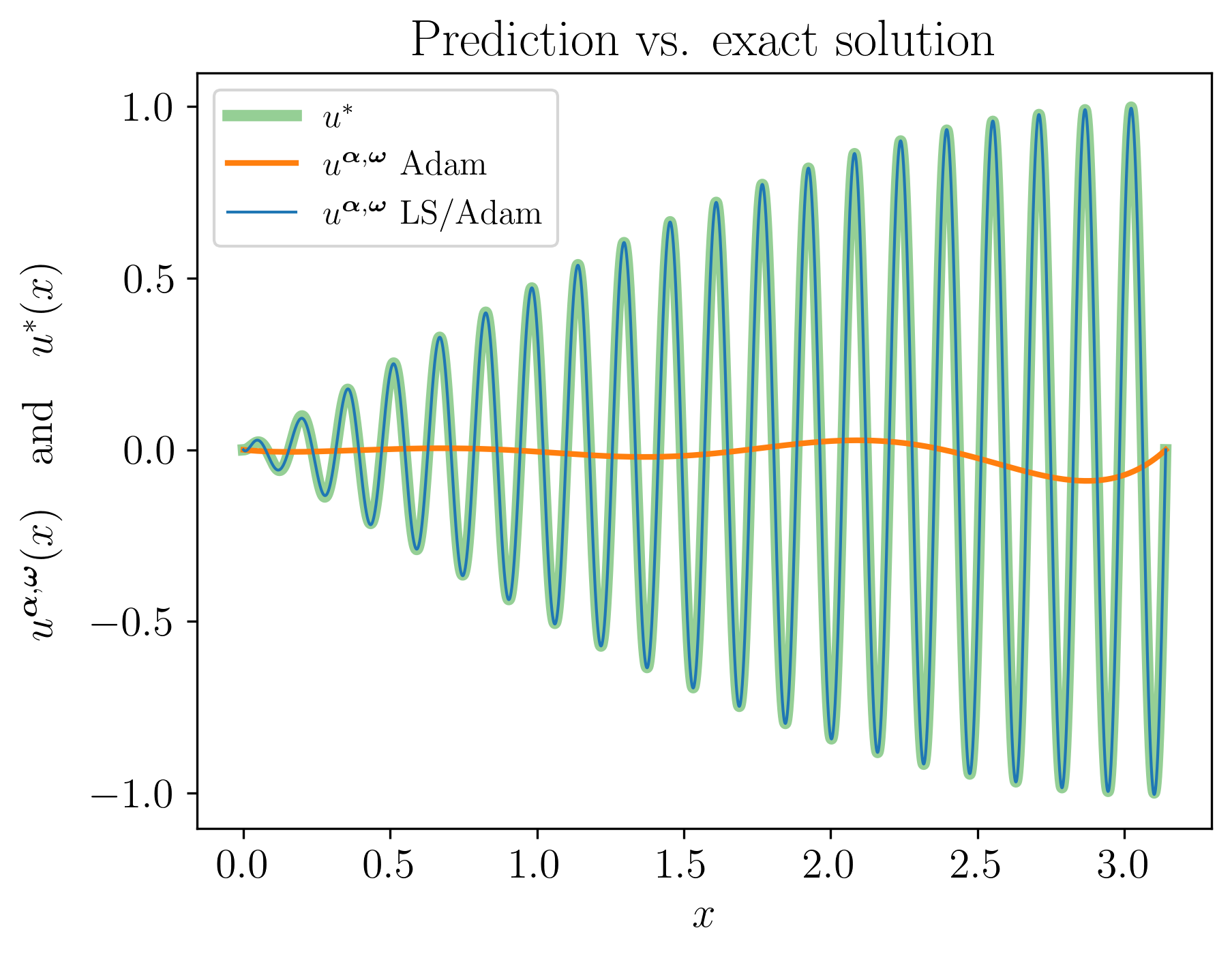}
        \end{subfigure}
		\begin{subfigure}[b]{0.49\textwidth}
		\includegraphics[scale=0.49]{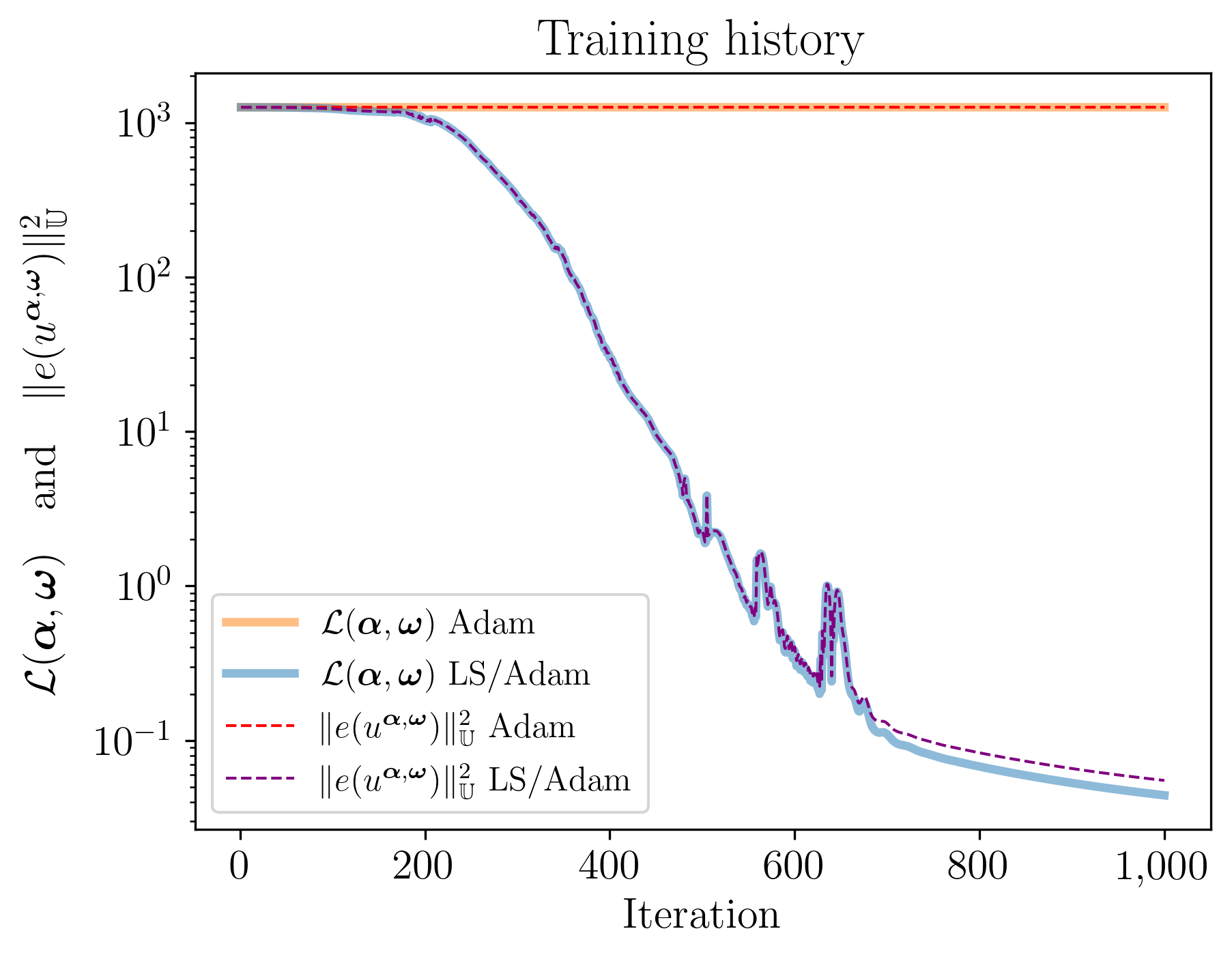}
        \end{subfigure}
        \begin{subfigure}[b]{0.49\textwidth}
		\includegraphics[scale=0.49]{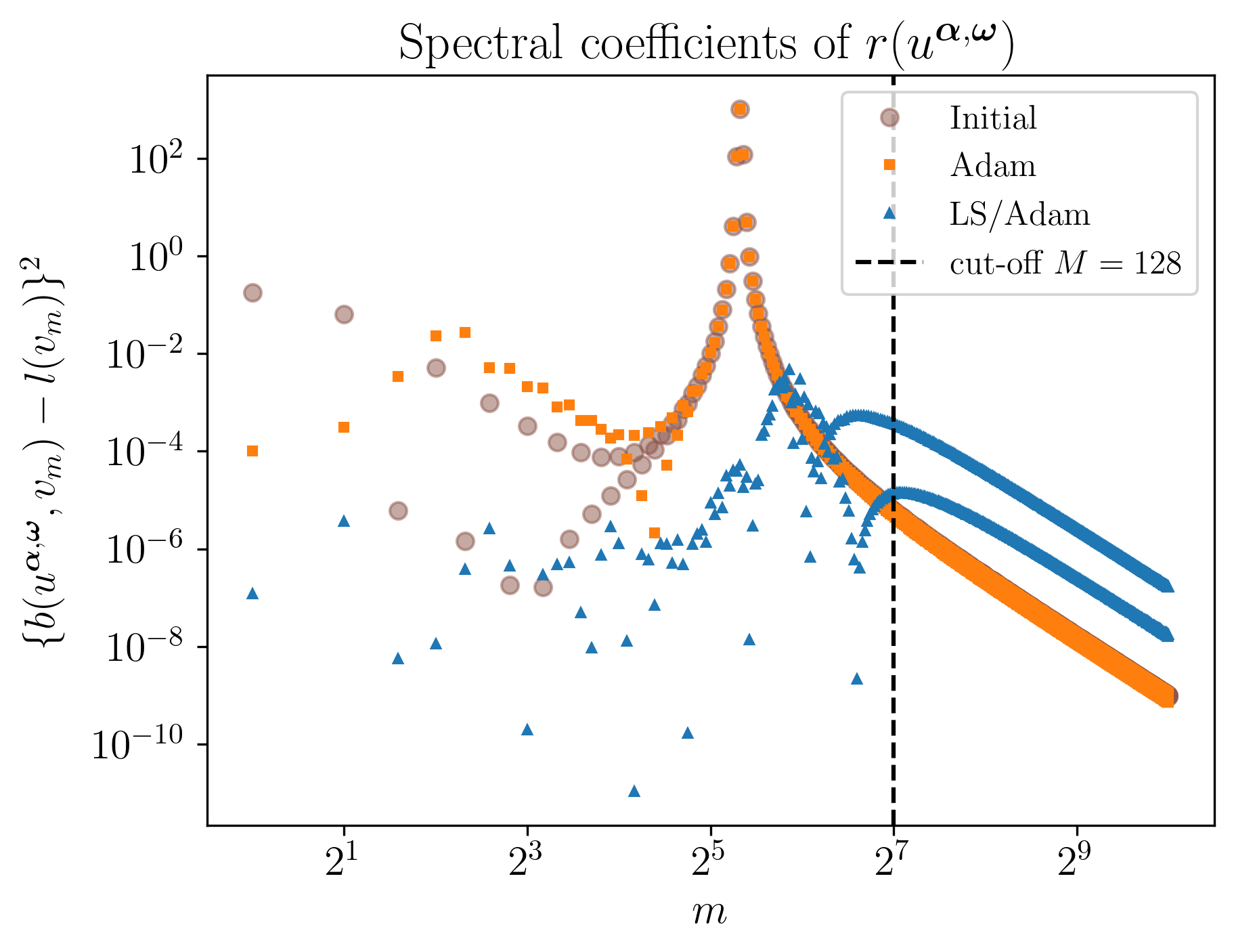}
		\end{subfigure}
  		\begin{subfigure}[b]{0.49\textwidth}
		\includegraphics[scale=0.49]{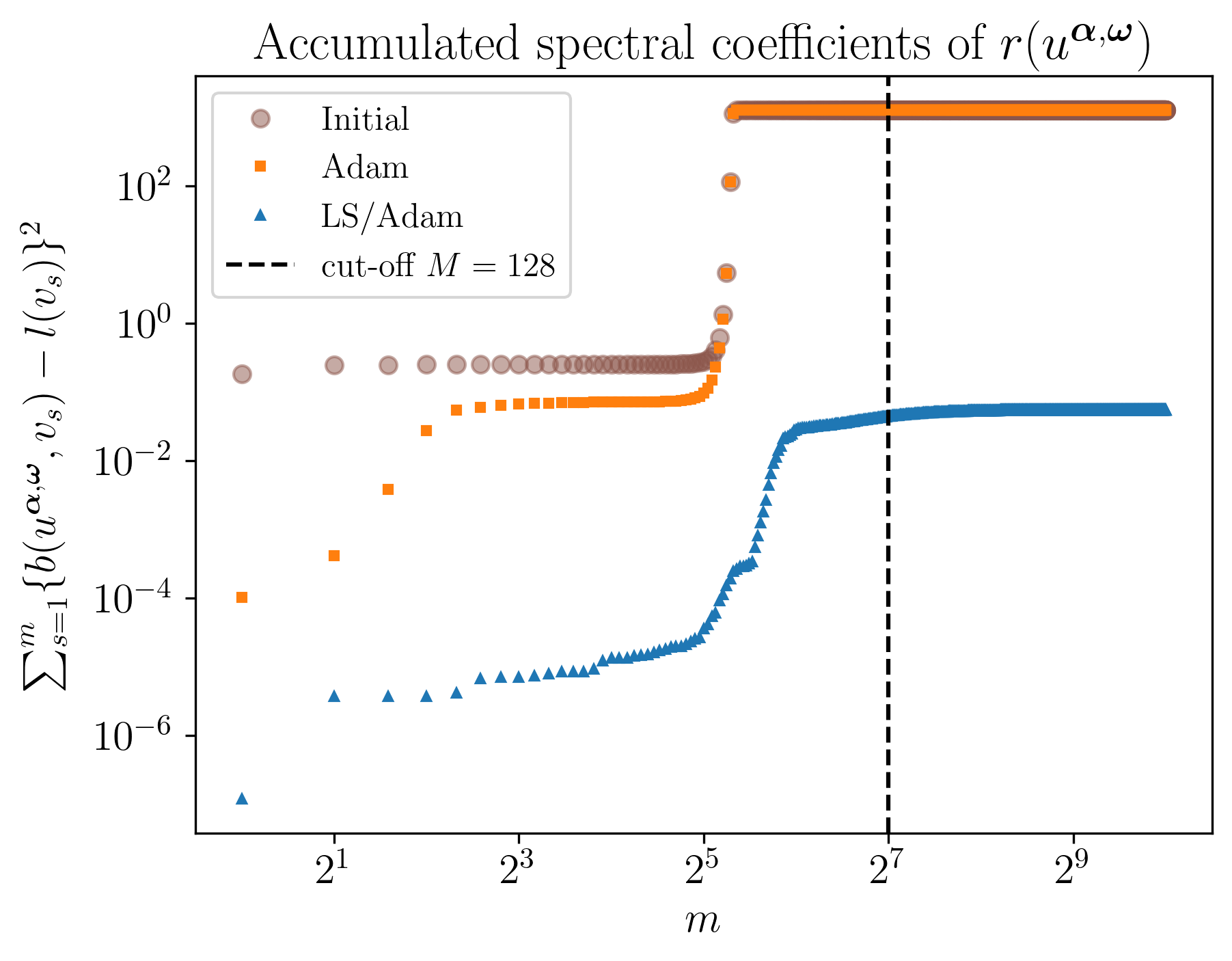}
		\end{subfigure}
		\caption{Results for VPINNs with manufactured solution $u^*(x)=\sin(40x)\sin(x/2)$ using either Adam or LS/Adam. The resulting final relative errors are $99.99\%$ and $0.66\%$ for Adam and LS/Adam, respectively.
  }
	\label{fig:exact2weak}
\end{figure}

Now, we increase the frequency of the manufactured solution by ten,
\begin{equation}\label{exact2}
	u^*(x)= \sin(40x)\sin\left( \frac{x}{2} \right),\qquad 0\leq x\leq \pi.
\end{equation} We increase the neural network discretization by four, $N=64$, and maintain the remaining discretization proportions and training conditions as before. \Cref{fig:exact2weak} shows the results for VPINNs (we obtain virtually identical results with UltraPINNs).

As in the previous example, we observe that the results in terms of accuracy during training are almost identical between VPINNs and UltraPINNs. Nevertheless, while the effect of the LS/Adam optimizer enables decreasing the loss to an acceptable level of convergence, the Adam optimizer alone is unable to significantly improve the network prediction from the initial state. This occurs due to the slow spectral bias of Adam \cite{rahaman2019spectral, cao2019towards,aldirany2024multi}: along the $1\mathord{,}000$-iteration training, Adam has only managed to barely decrease the very first two spectral modes of the error, leaving practically immobile all the remaining ones. Informally speaking, Adam had no time to ``discover'' the modes that most significantly represent the essence of the error (around mode number $40$). In contrast, LS/Adam did manage to decrease this significant frequency by about four orders of magnitude during the $1\mathord{,}000$-iteration training. The correlation between $\mathcal{L}(\boldsymbol{\alpha},\boldsymbol{\omega})$ and $\Vert e(u^{\boldsymbol{\alpha},\boldsymbol{\omega}})\Vert_{\mathbb{U}}^2$ is practically one-to-one during training.

\begin{figure}[htbp]
	\centering
        \begin{subfigure}[b]{0.49\textwidth}
	    \includegraphics[scale=0.49]{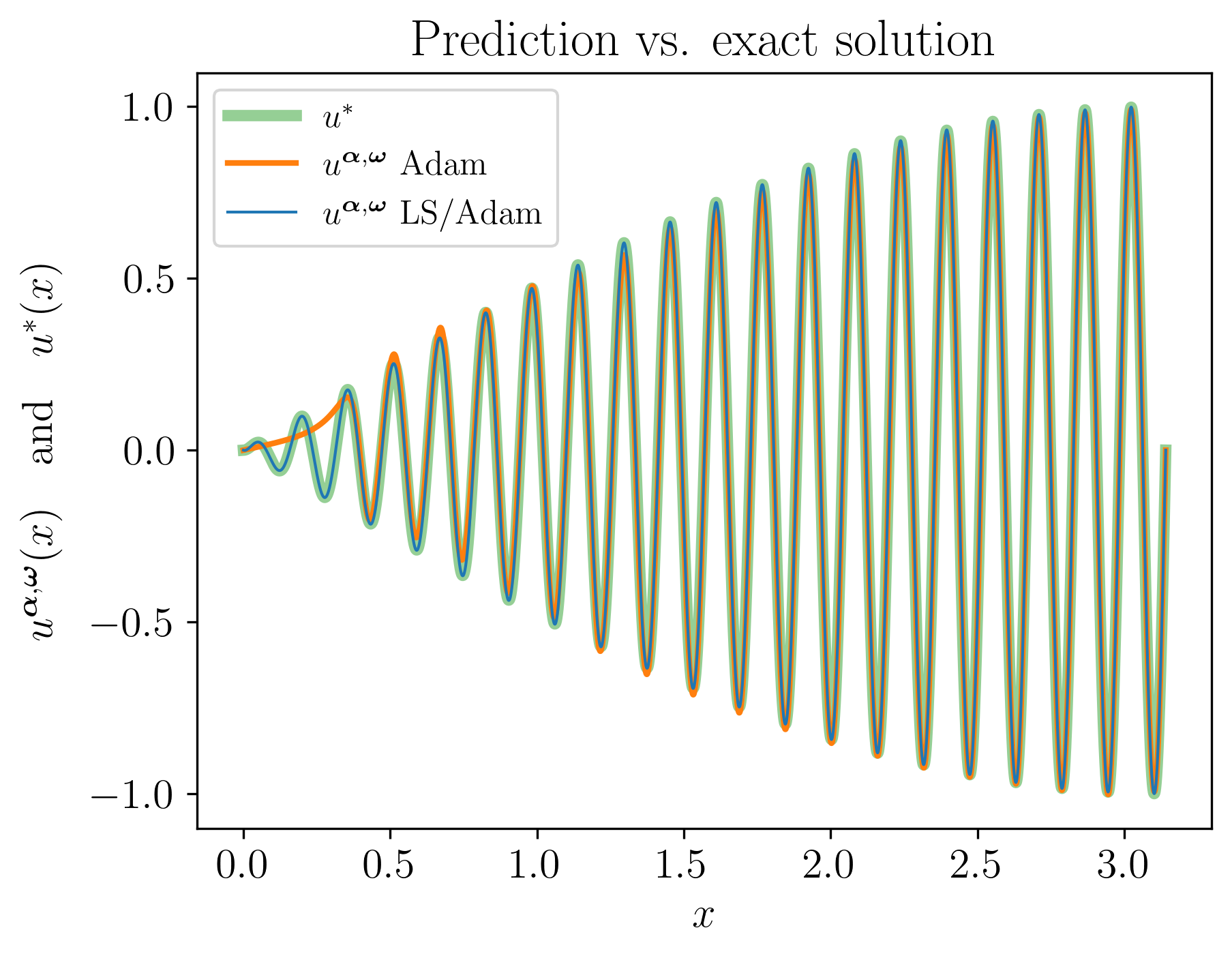}
        \end{subfigure}
		\begin{subfigure}[b]{0.49\textwidth}
		\includegraphics[scale=0.49]{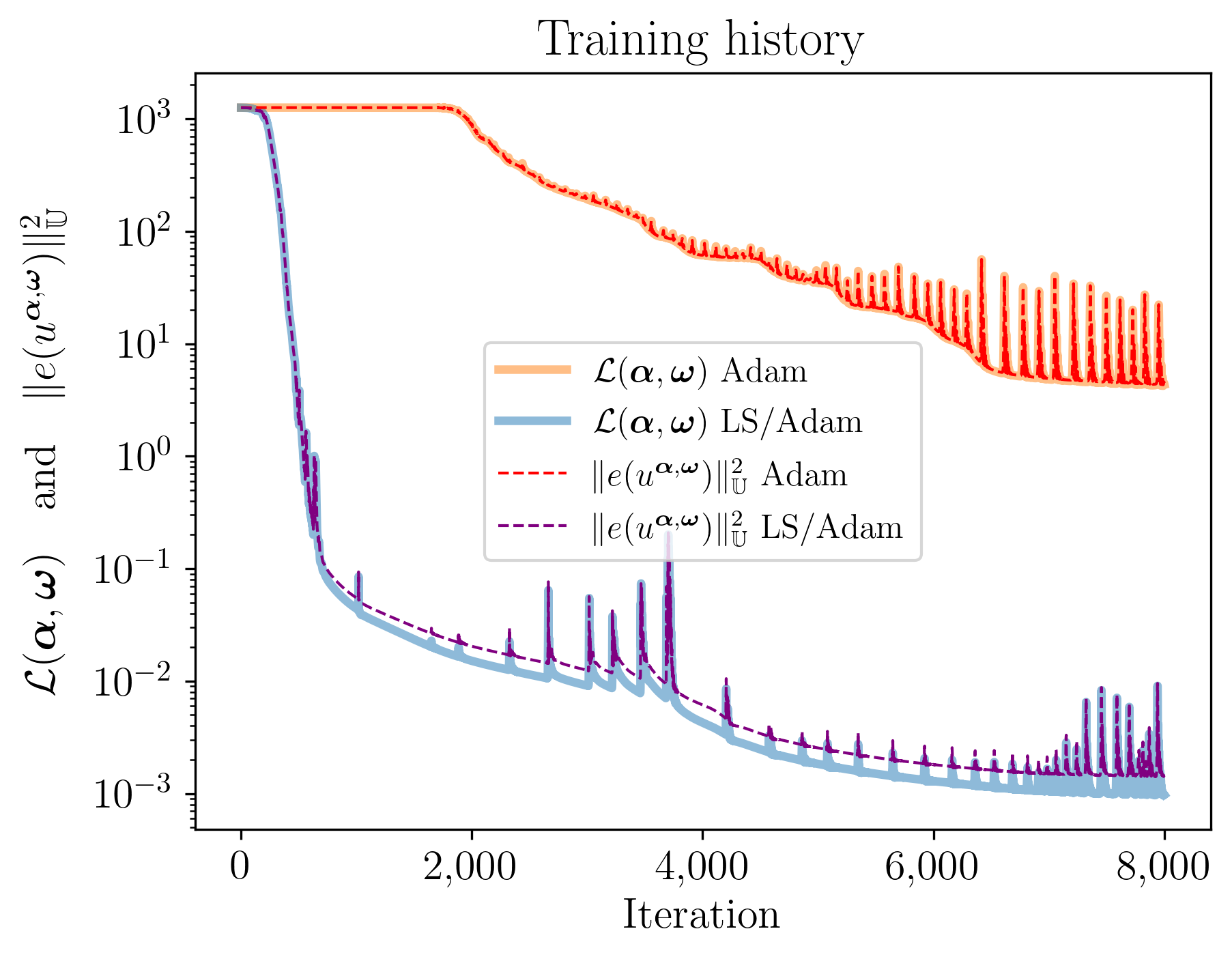}
        \end{subfigure}
        \begin{subfigure}[b]{0.49\textwidth}
		\includegraphics[scale=0.49]{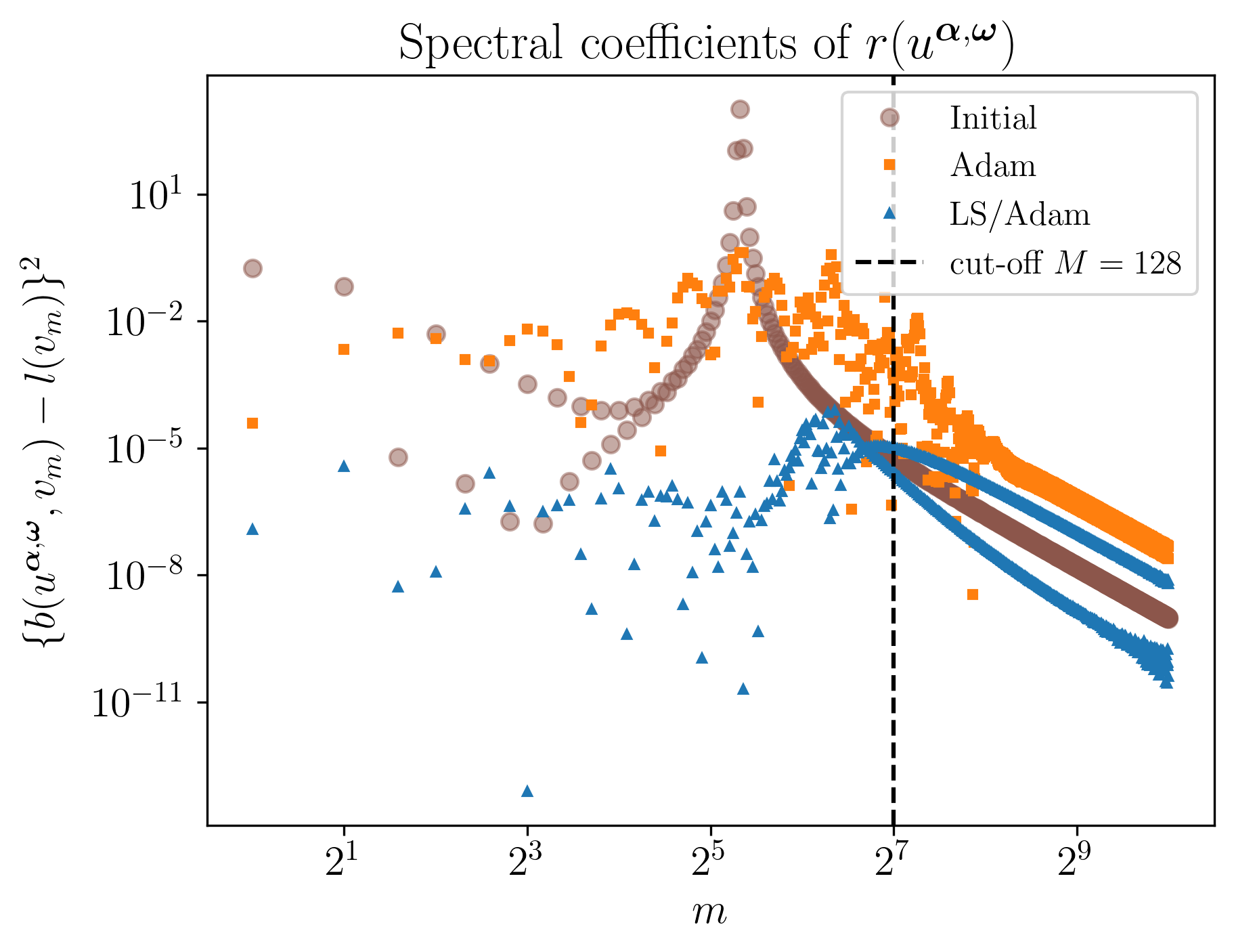}
		\end{subfigure}
  		\begin{subfigure}[b]{0.49\textwidth}
		\includegraphics[scale=0.49]{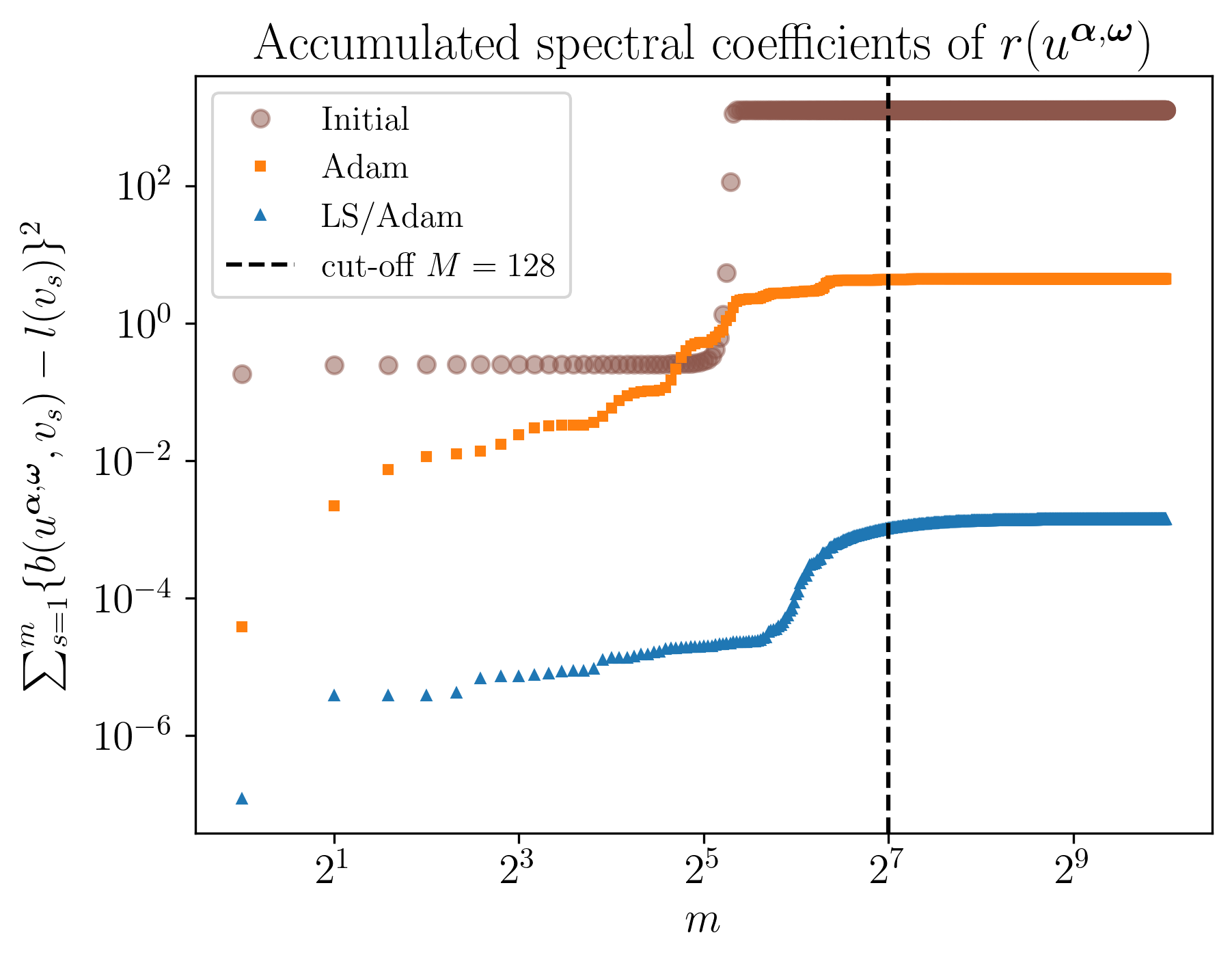}
		\end{subfigure}
		\caption{Results for VPINNs with manufactured solution $u^*(x)=\sin(40x)\sin(x/2)$ using either Adam or LS/Adam. The resulting final relative errors after $8\mathord{,}000$ iterations are $5.94\%$ and $0.11\%$ for Adam and LS/Adam, respectively.
  }
	\label{fig:experiment3weakmore}
\end{figure}

This example illustrates that although the discretization of the test space might be rich enough to adequately represent the error, Adam's optimization alone may be ineffective. By applying LS, we have significantly accelerated the training time required to reach an acceptable accuracy threshold. \Cref{fig:experiment3weakmore} shows the results of a $8\mathord{,}000$-iteration training using VPINNs under the aforementioned conditions. As we observe, Adam alone has needed about $3\mathord{,}000$ iterations to ``discover'' the critical mode of the error and consequently be able to escape from the initial plateau state and minimize the loss slowly.

\subsection{A singular 1D problem with exploding gradient}\label{section5.4}

Now, we consider a case where the exact solution presents a singularity in one of the endpoints of the boundary,
\begin{equation}
	u^*(x)=x^{\beta}(\pi-x),\qquad \beta=0.7,\qquad 0\leq x\leq \pi.
\end{equation}
This problem has as source $f\in H^{-1}(\Omega)\setminus L^2(\Omega)$, which makes the underlying variational problem to be well defined in weak form but not in strong form. 

For the discretization, we consider $N=16$ and $M\in\{32, 64, 128\}$, letting the remaining proportions and conditions as before. We consider an integration grid with half of the integration points following an exponential distribution in $(0,1)$ towards $x=0$ while the remaining half is uniformly distributed in $(1,\pi)$. 

\begin{figure}[htbp]
		\centering
  \begin{subfigure}[b]{0.32\textwidth}
		\includegraphics[scale=0.33]{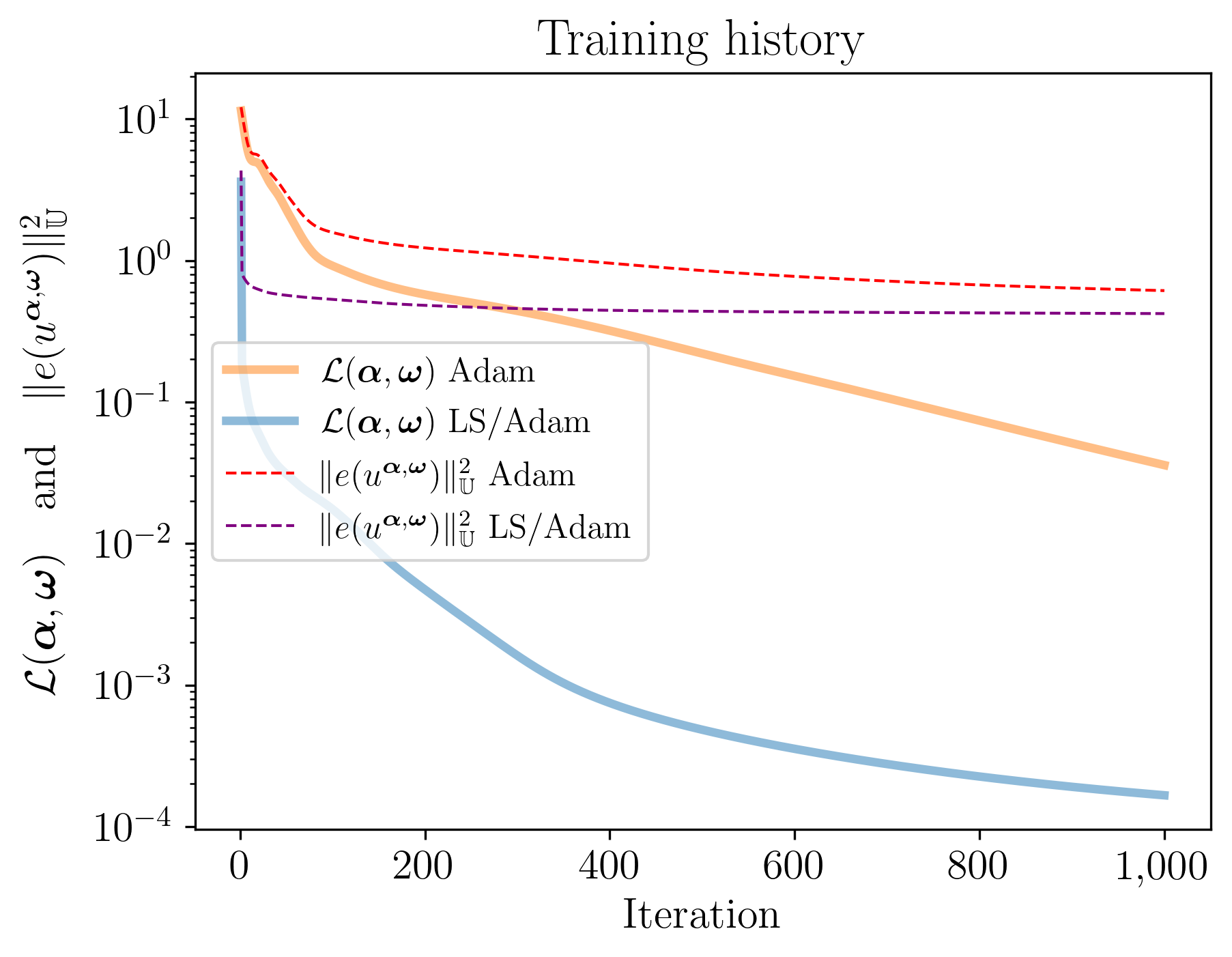}
  \end{subfigure}
  \begin{subfigure}[b]{0.32\textwidth}
		\includegraphics[scale=0.33]{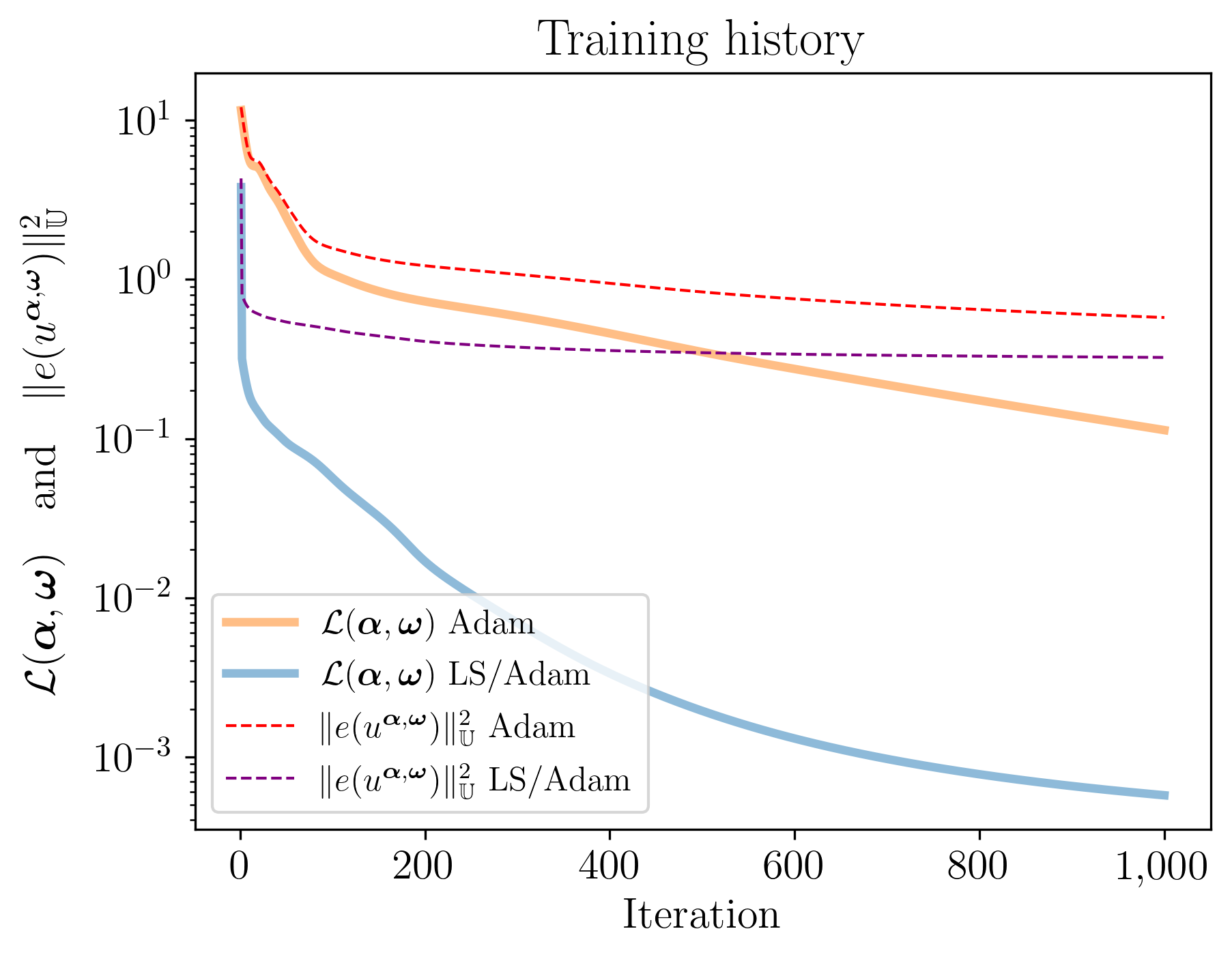}
  \end{subfigure}
  \begin{subfigure}[b]{0.32\textwidth}
  		\includegraphics[scale=0.33]{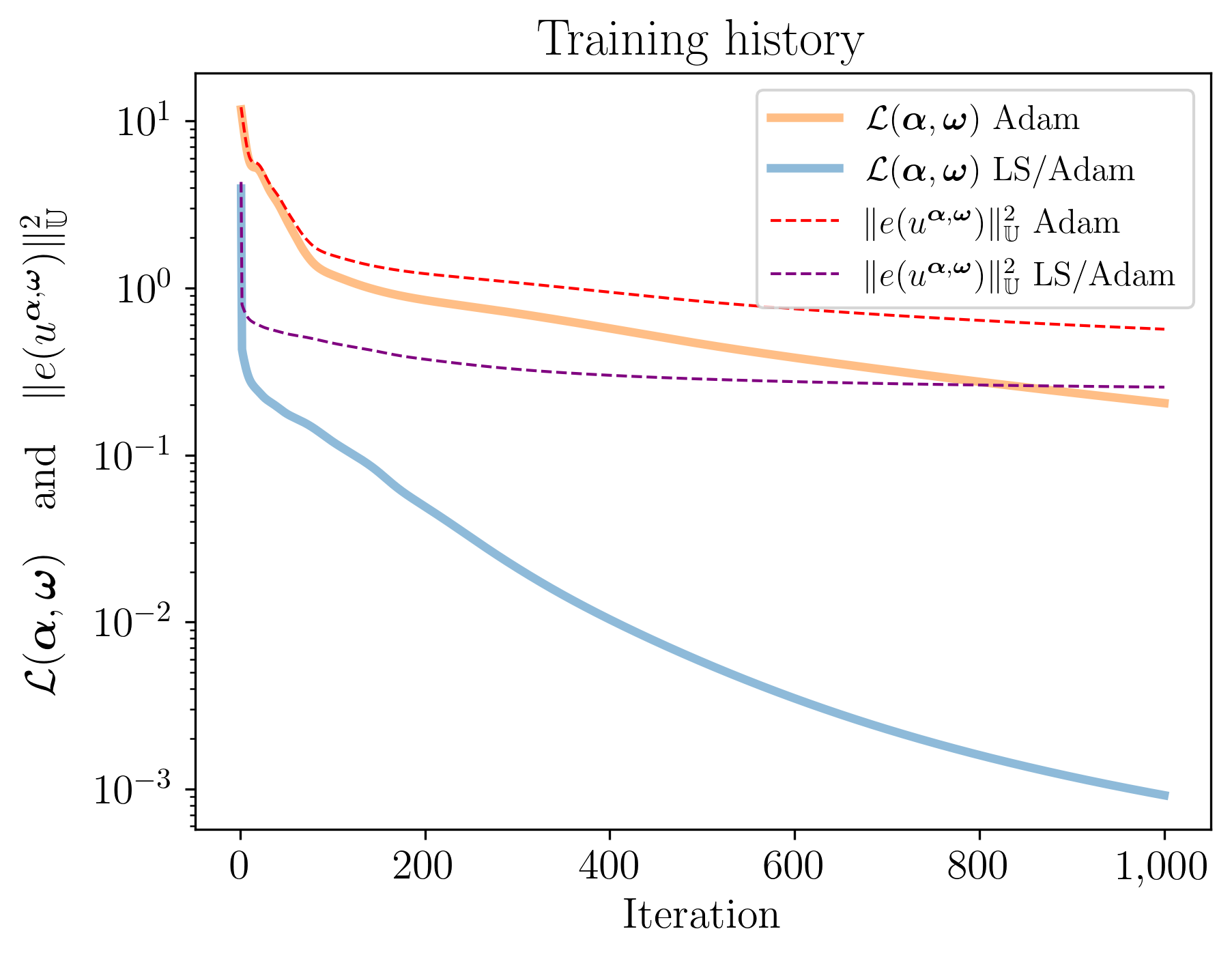}
    \end{subfigure}
      \begin{subfigure}[b]{0.32\textwidth}
		\includegraphics[scale=0.33]{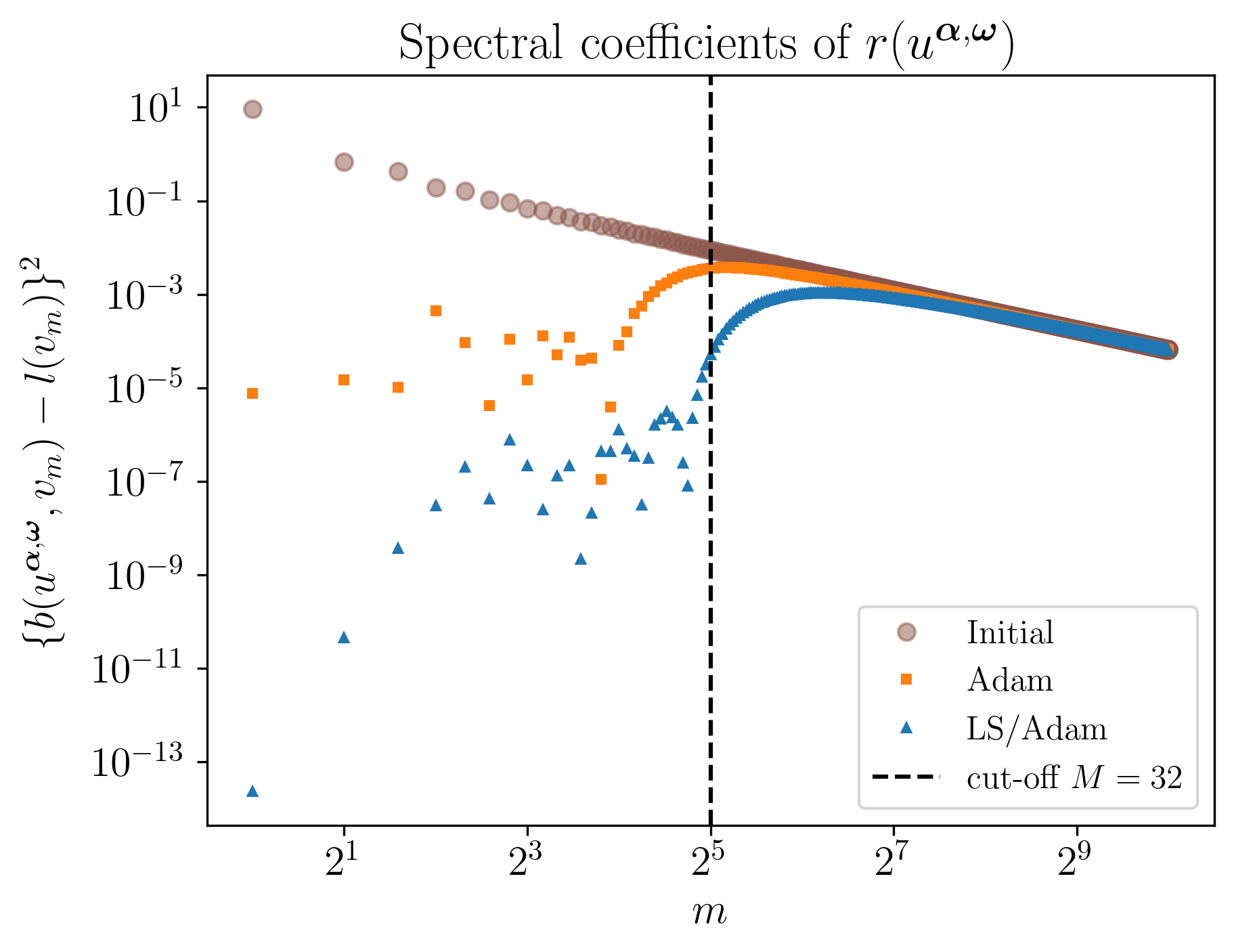}
  \end{subfigure}
  \begin{subfigure}[b]{0.32\textwidth}
		\includegraphics[scale=0.33]{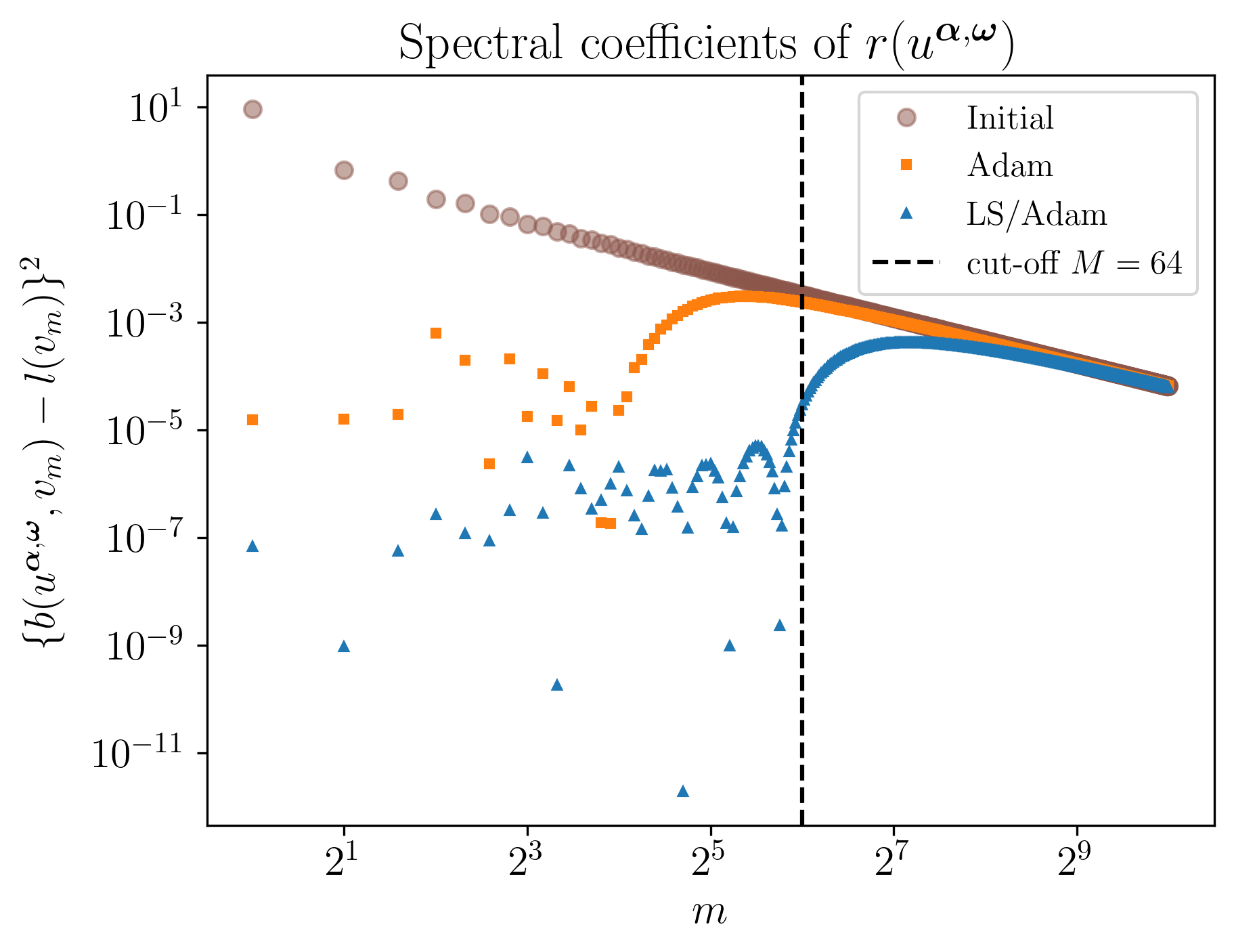}
  \end{subfigure}
  \begin{subfigure}[b]{0.32\textwidth}
  		\includegraphics[scale=0.33]{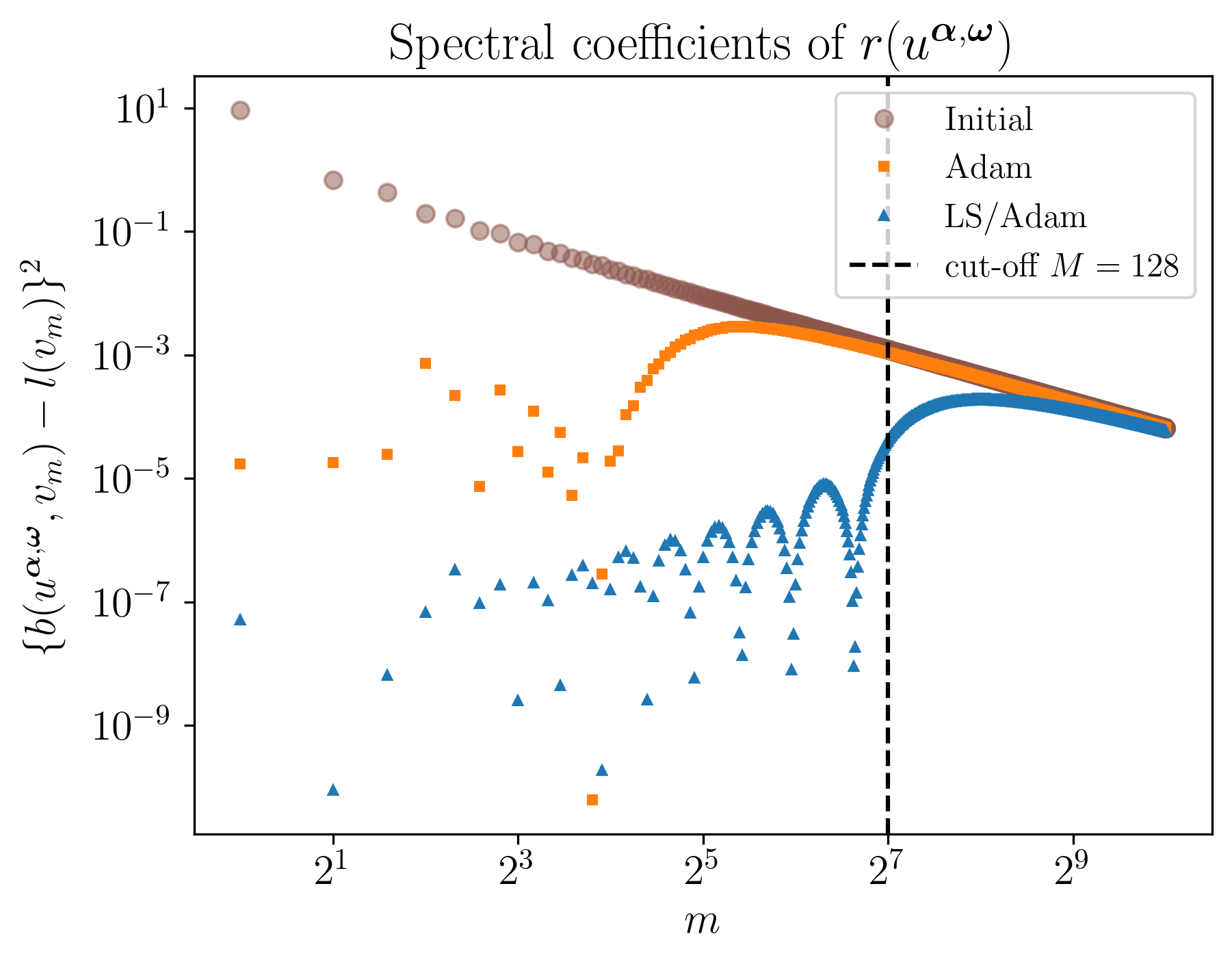}
    \end{subfigure}
      \begin{subfigure}[b]{0.32\textwidth}
		\includegraphics[scale=0.33]{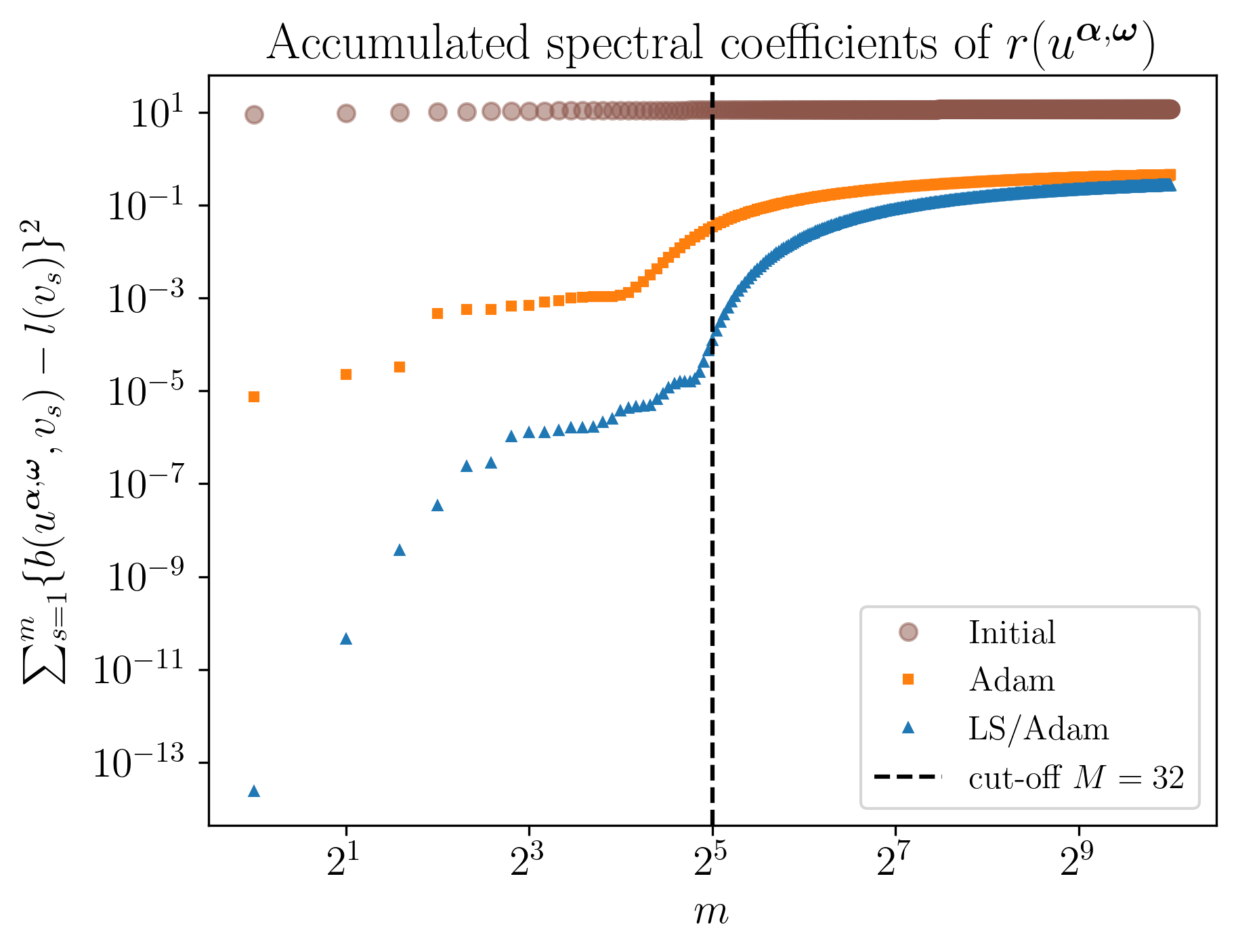}
        \caption{$M=32$}
  \end{subfigure}
  \begin{subfigure}[b]{0.32\textwidth}
		\includegraphics[scale=0.33]{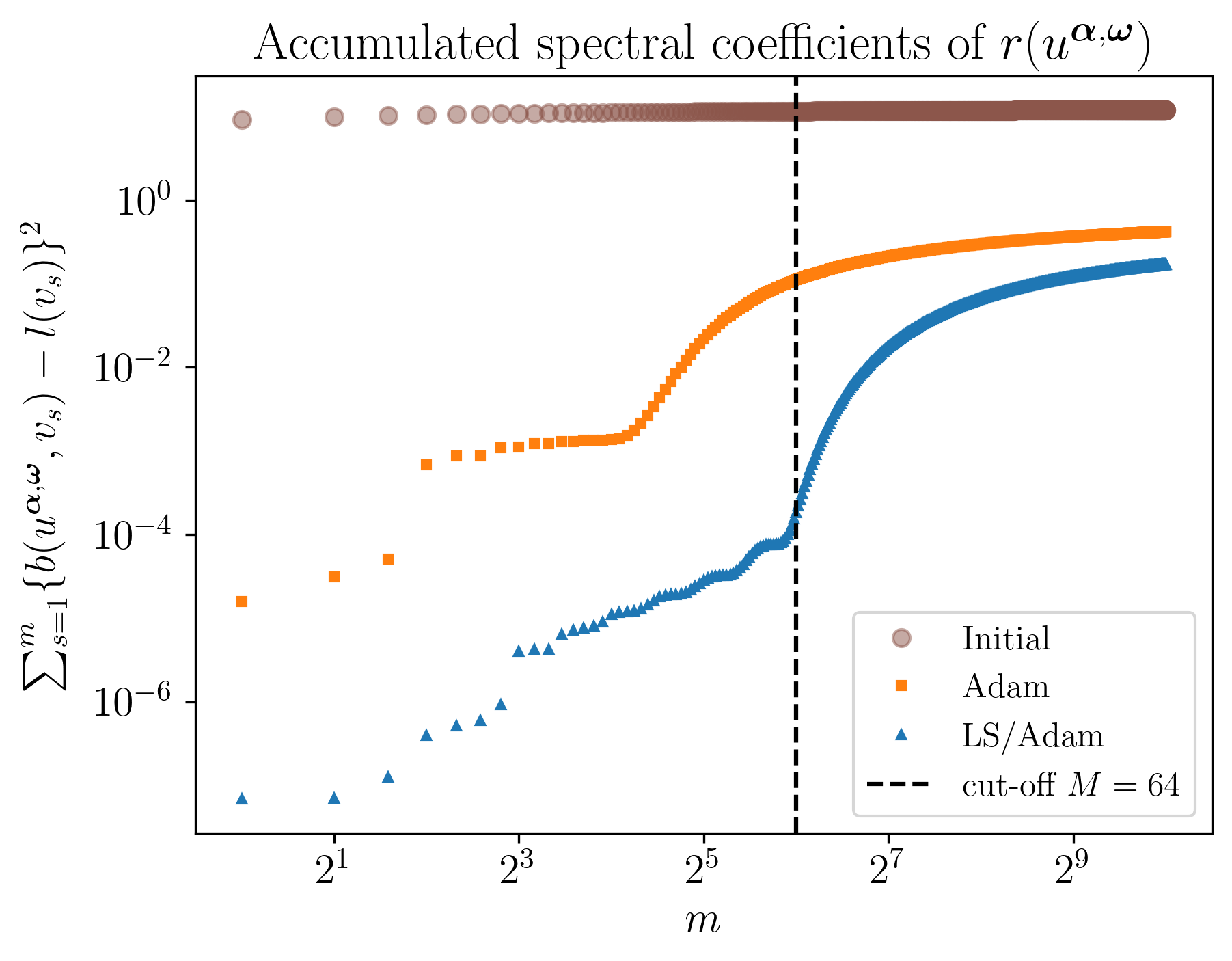}
        \caption{$M=64$}
  \end{subfigure}
  \begin{subfigure}[b]{0.32\textwidth}
  		\includegraphics[scale=0.33]{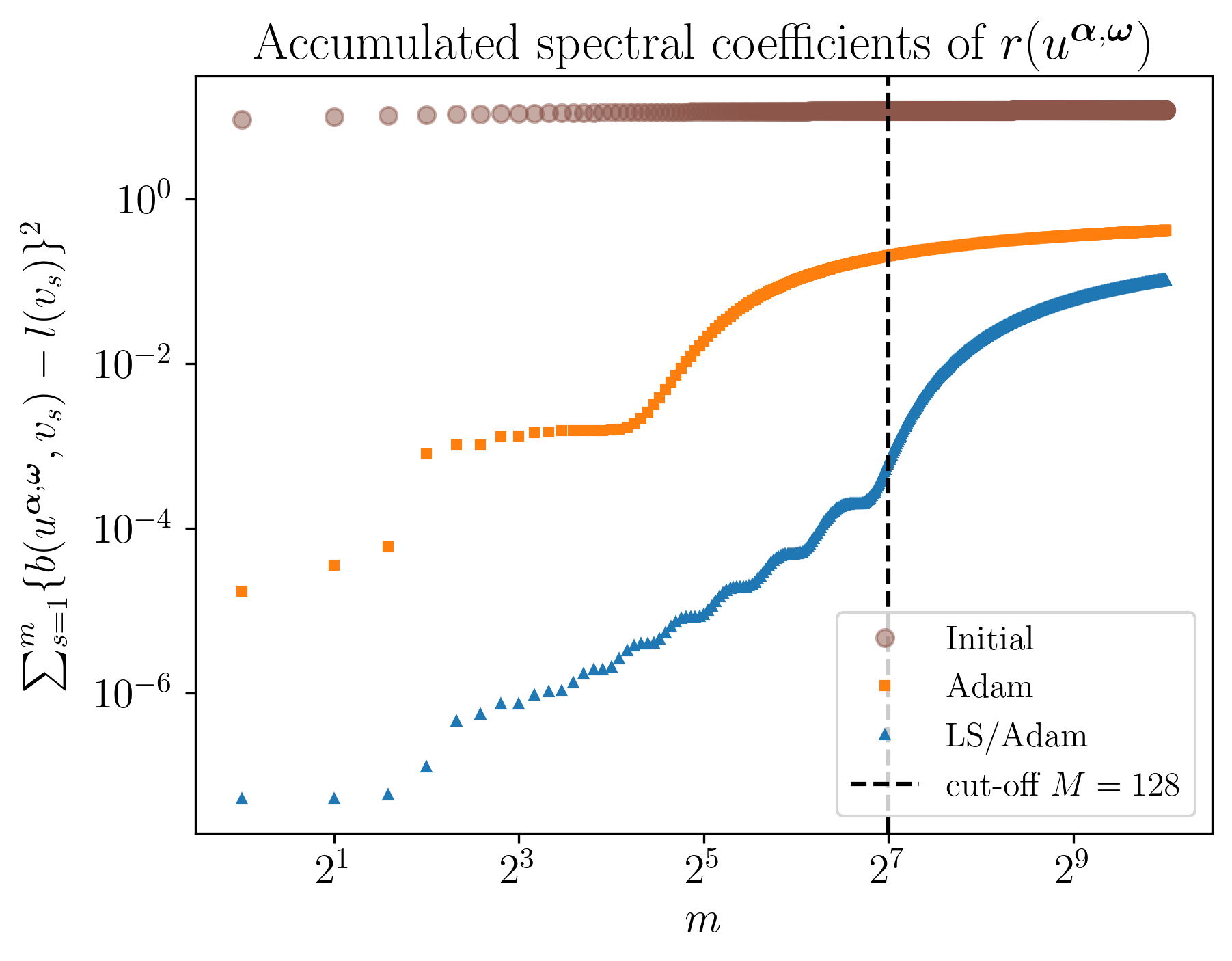}
        \caption{$M=128$}
    \end{subfigure}
		\caption{Training of VPINNs with manufactured solution $u^*(x)=x^\beta(\pi-x)$ using either Adam or LS/Adam for different test space discretizations. The resulting final relative errors for Adam and LS/Adam) are, respectively: $23.20\%$ and $19.25\%$ for $M=32$, $22.52\%$ and $16.87\%$ for $M=64$, and $22.33\%$ and $14.98\%$ for $M=128$.
  }
	\label{fig:experiment4}
\end{figure}

\Cref{fig:experiment4} shows the $1\mathord{,}000$-iteration training history and (accumulated) spectral coefficients when using VPINNs for the different test space discretization cut-offs. Experimentation with UltraPINNs produced worse results due to numerical integration issues around the singularity whose discussion falls beyond the scope of this work.

We observe a much better convergence of the loss with LS/Adam compared to Adam alone, as in the previous problems. Unfortunately, most of the error is above the training cut-off frequency employed in the loss, and thus, the error is dominated by this part of the spectrum of the residual, i.e., $\Vert r(u^{\boldsymbol{\alpha},\boldsymbol{\omega}}) - r_M(u^{\boldsymbol{\alpha},\boldsymbol{\omega}}) \Vert_{\mathbb{V}}^2$ is more significant than $\mathcal{L}(\boldsymbol{\alpha},\boldsymbol{\omega})=\Vert r_M(u^{\boldsymbol{\alpha},\boldsymbol{\omega}}) \Vert_{\mathbb{V}}^2$. To reduce the error below 15\%, we would need to increase the cut-off frequency beyond $M=128$, which would involve large computational costs. Alternatively, we could employ singularity-oriented discretizations (see, e.g., \cite{taylor2024adaptive,taylor2024reconns}), but this falls outside the scope of this work.

\subsection{Two smooth 2D problems}\label{section5.5}

\begin{figure}[htbp]
		\centering
  \begin{subfigure}[b]{0.49\textwidth}
		\includegraphics[scale=0.49]{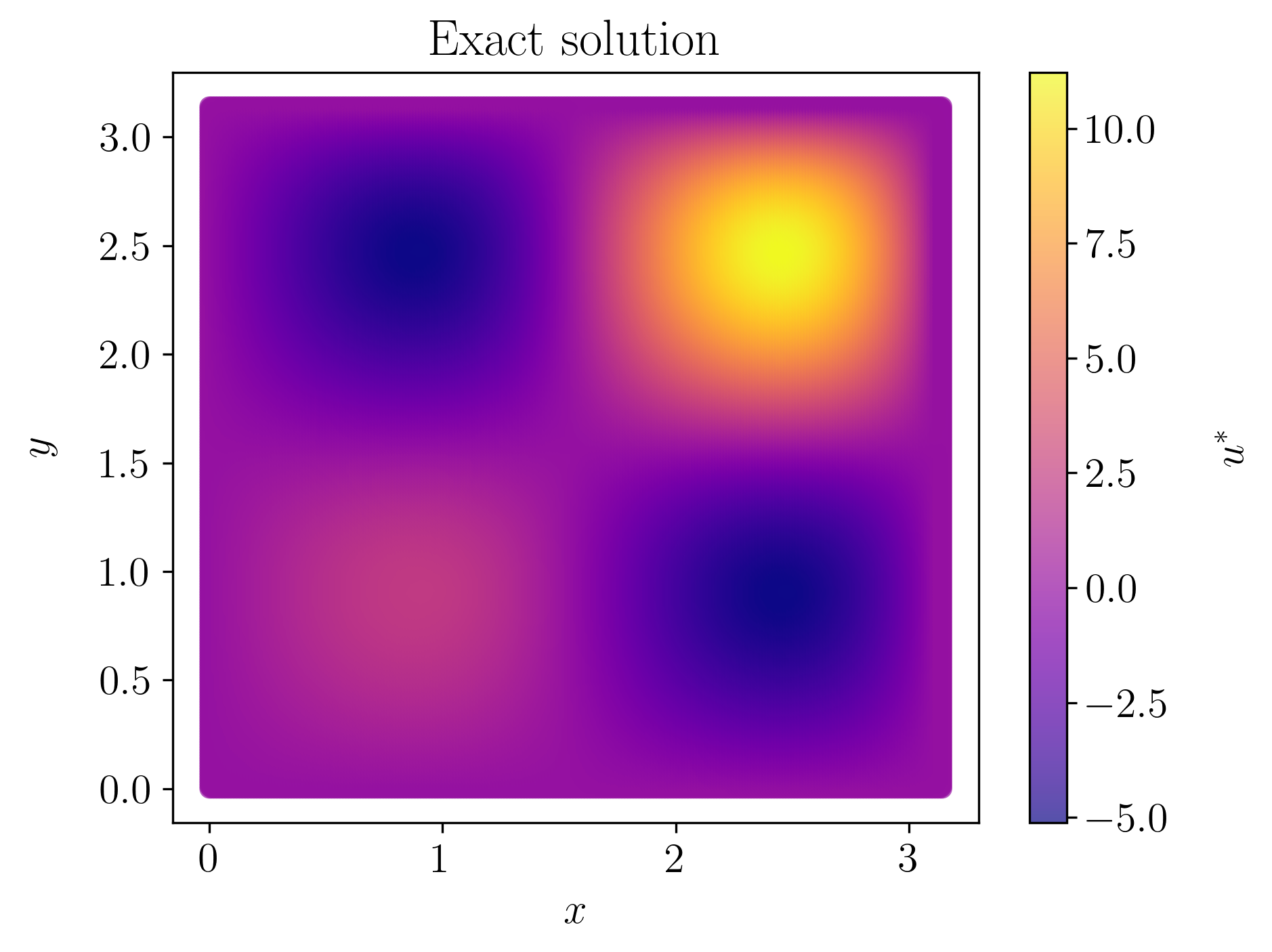}
  \end{subfigure}\hfill
  \begin{subfigure}[b]{0.49\textwidth}
		\includegraphics[scale=0.49]{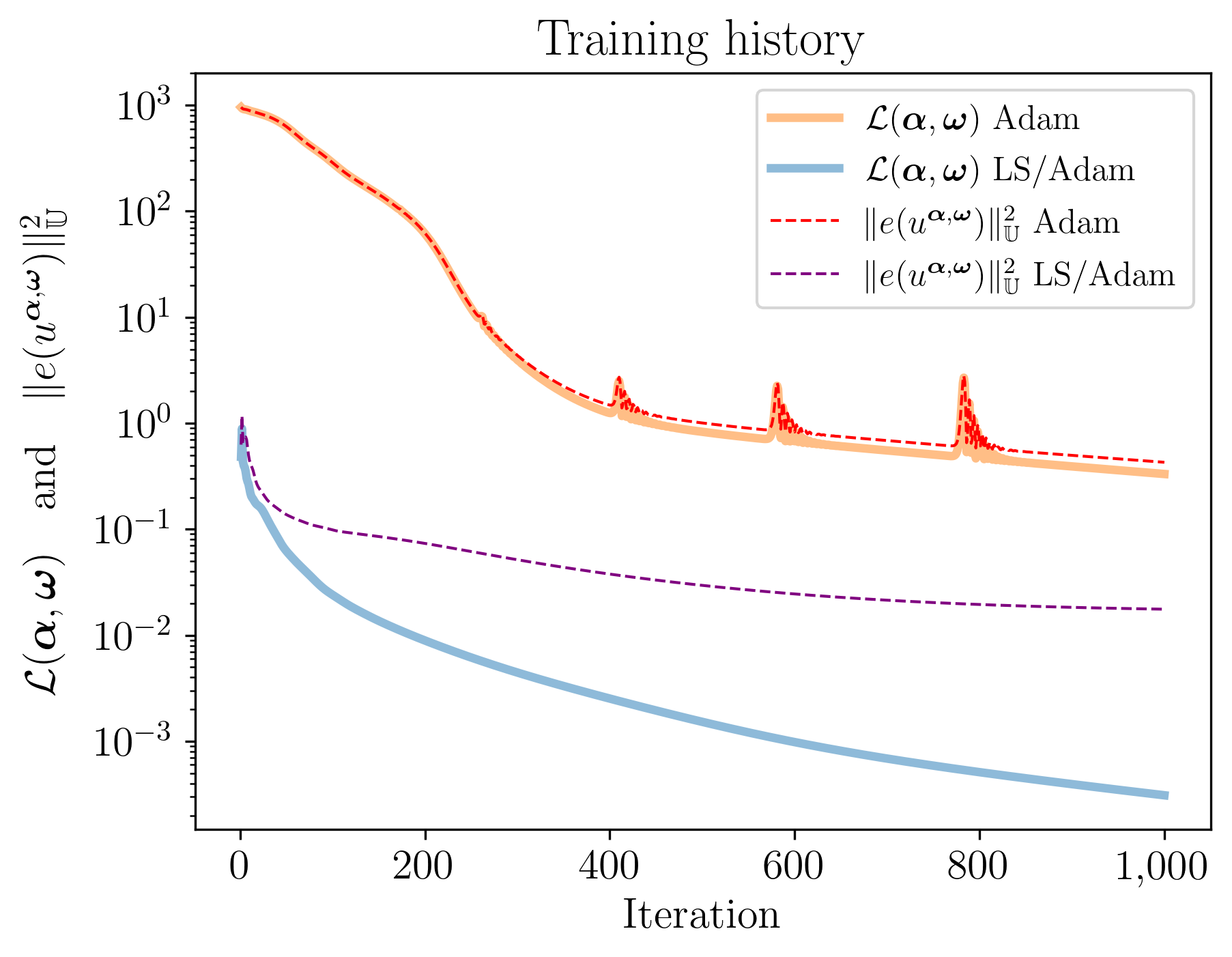}
  \end{subfigure}
  \begin{subfigure}[b]{0.49\textwidth}
  		\includegraphics[scale=0.49]{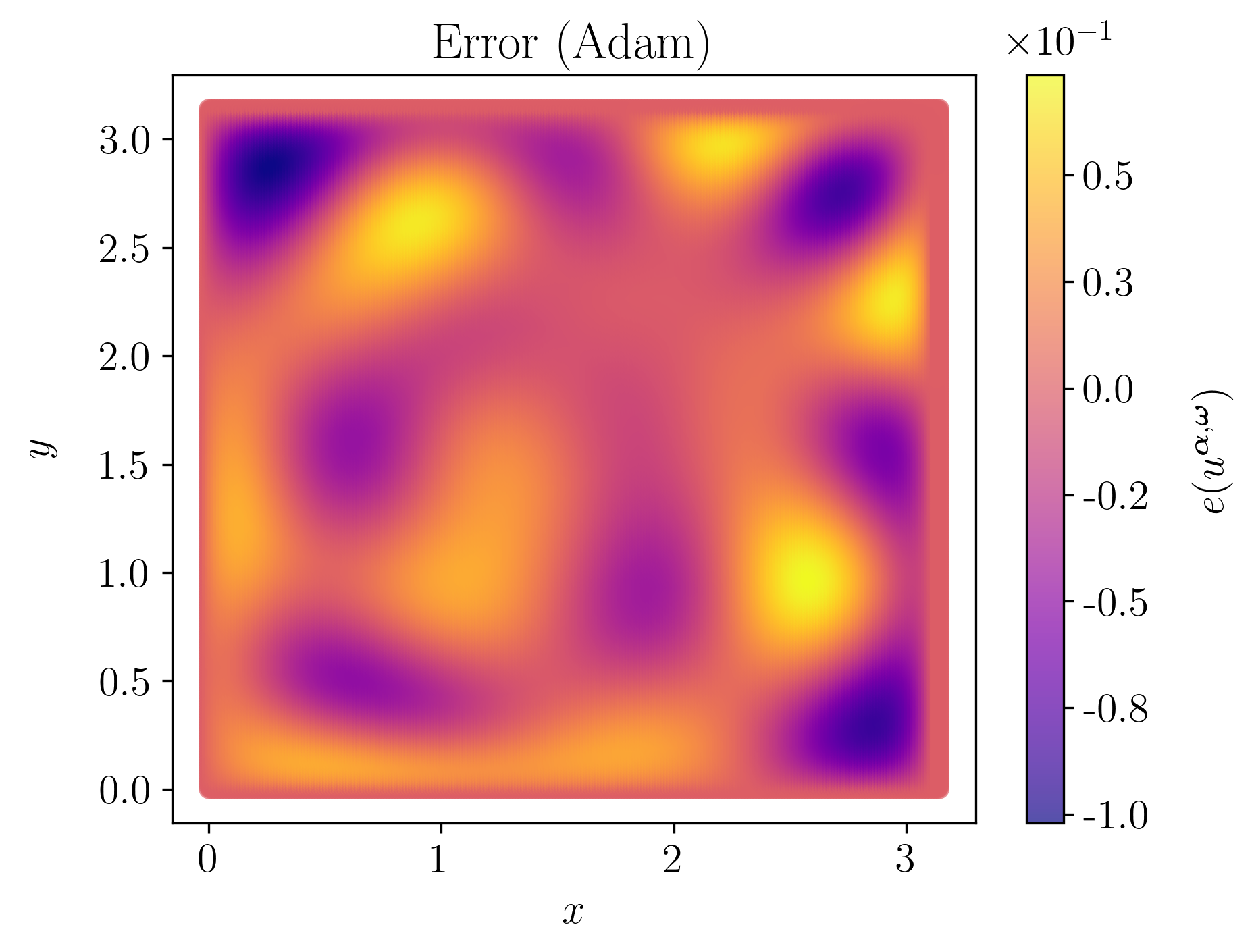}
    \end{subfigure}\hfill
      \begin{subfigure}[b]{0.49\textwidth}
		\includegraphics[scale=0.49]{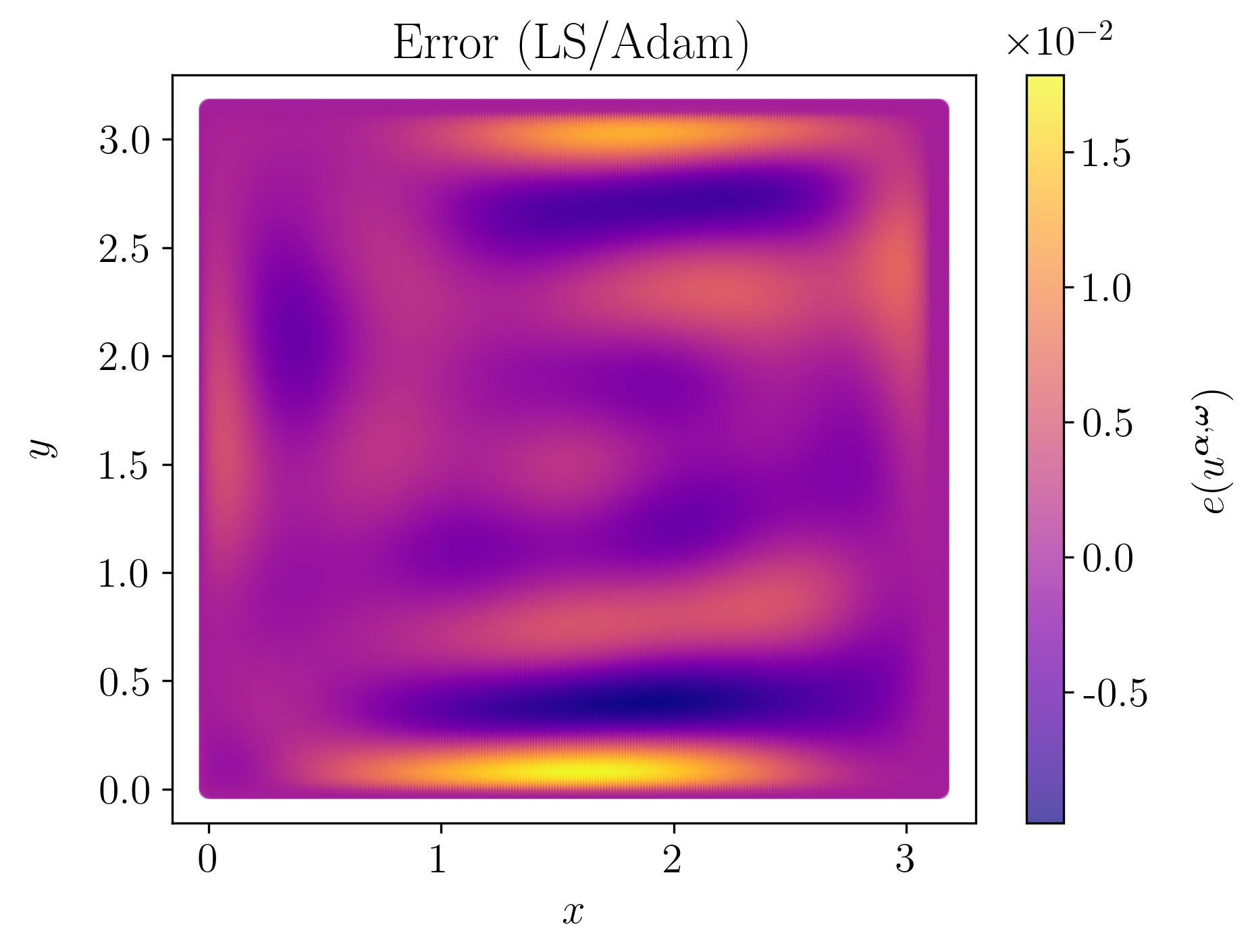}
  \end{subfigure}
\caption{Training of VPINNs with manufactured solution $u^*(x,y)=\sin(2x)\sin(2y)e^{\frac{x+y}{2}}$ using either Adam or LS/Adam. The resulting final relative errors are $2.15\%$ and $0.44\%$ for Adam and LS/Adam, respectively.}
\label{fig:experiment5one}
\end{figure}

Let
\begin{equation}
	u^*(x)=\sin(2\eta x)\sin(2y)e^{\frac{x+y}{2}},\qquad 0\leq x,y\leq \pi,
\end{equation} and consider the cases $\eta\in\{1,10\}$ over $1\mathord{,}000$ training iterations. For $\eta=1$, we consider the discretization setup given by $N=32$ for the trial space, $M = 8\times 8$ for the test space, $K=128\times 128$ for integration during training, and $256\times 256$ for validation. \Cref{fig:experiment5one} shows the obtained results.

As we observe, an increasing discrepancy between the error and the loss arises from very early iterations when using LS/Adam. This occurs due to the limited richness of the considered test space discretization that covers low frequencies only. 

\begin{figure}[htbp]
		\centering
  \begin{subfigure}[b]{0.49\textwidth}
		\includegraphics[scale=0.49]{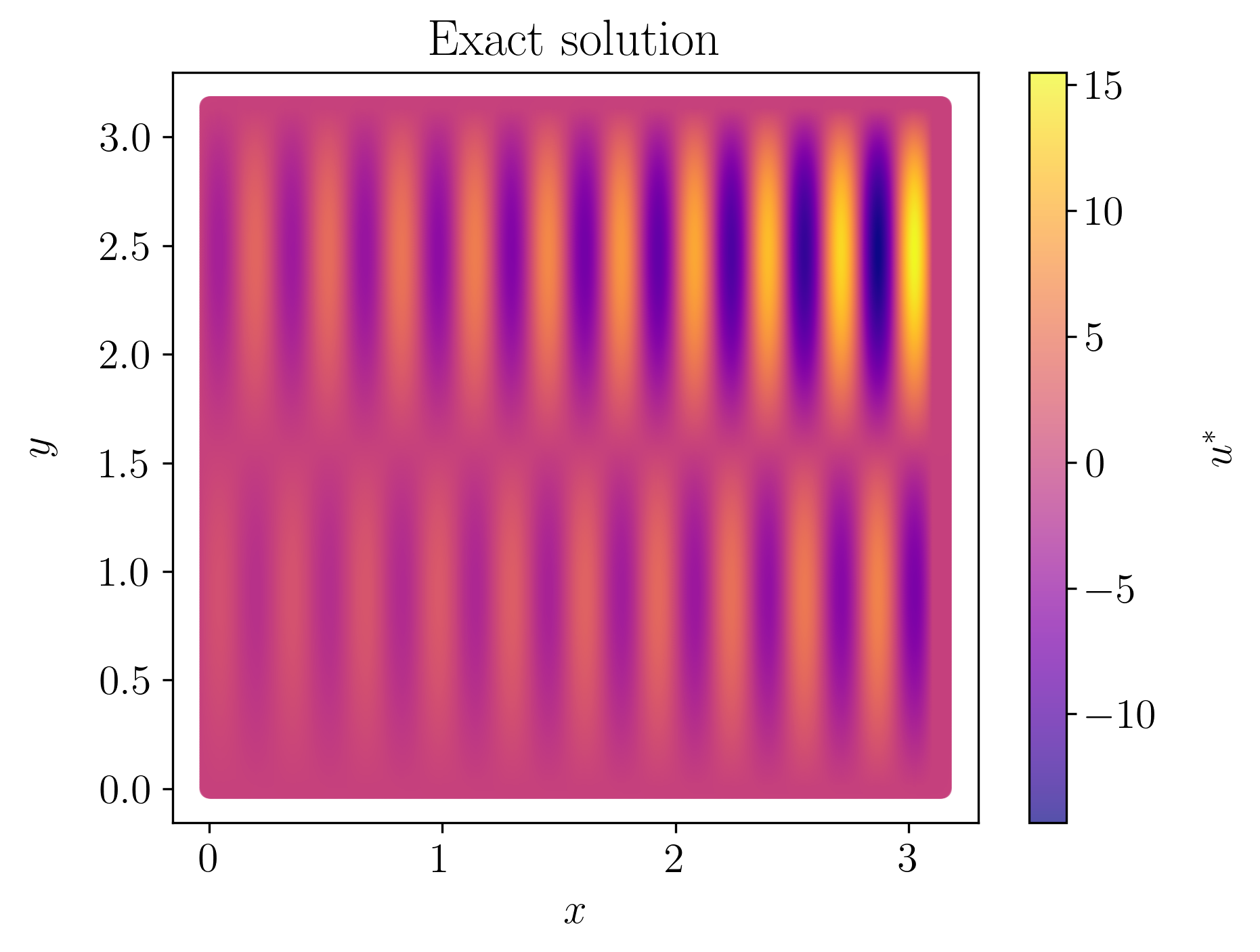}
  \end{subfigure}\hfill
  \begin{subfigure}[b]{0.49\textwidth}
		\includegraphics[scale=0.49]{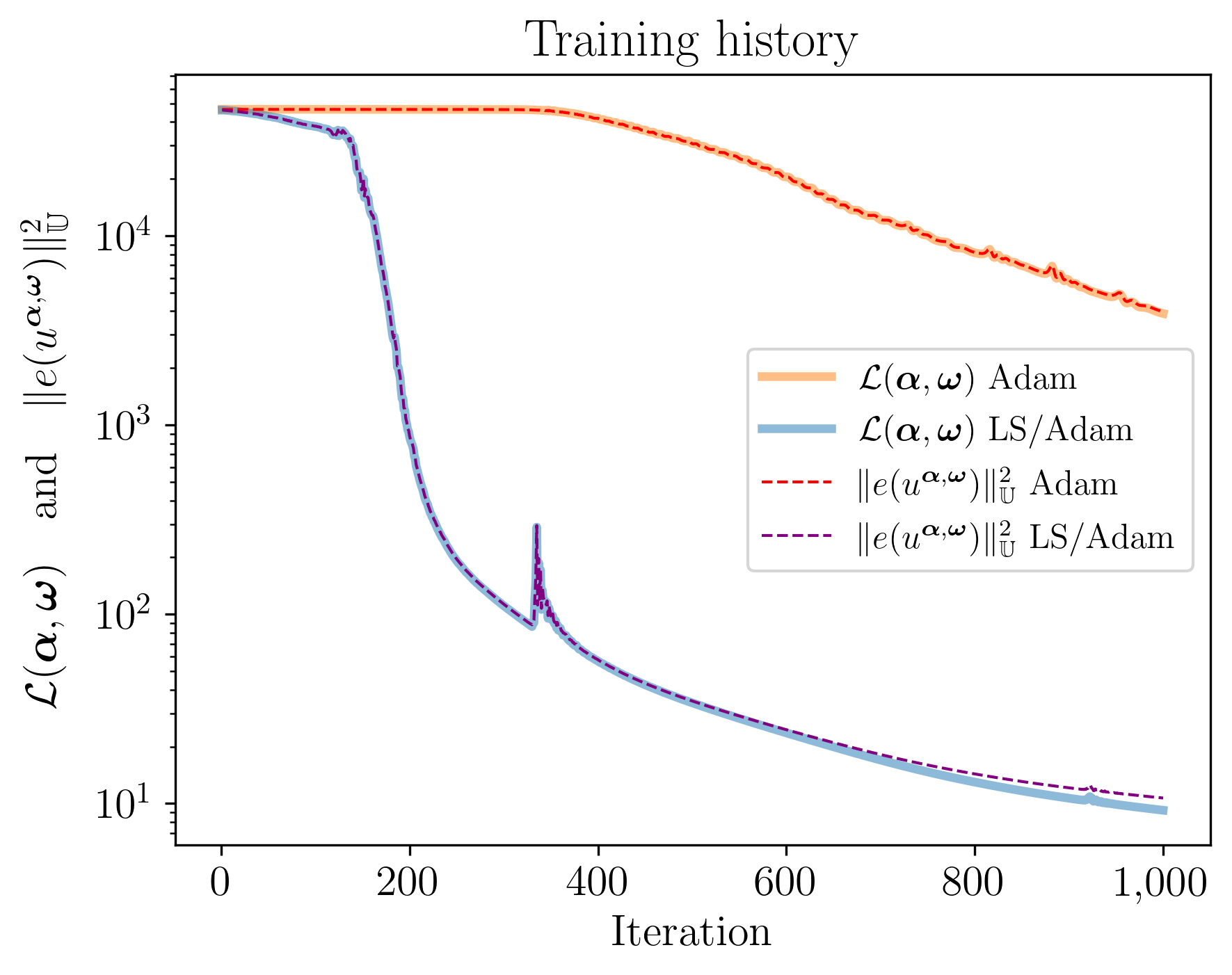}
  \end{subfigure}
  \begin{subfigure}[b]{0.49\textwidth}
  		\includegraphics[scale=0.49]{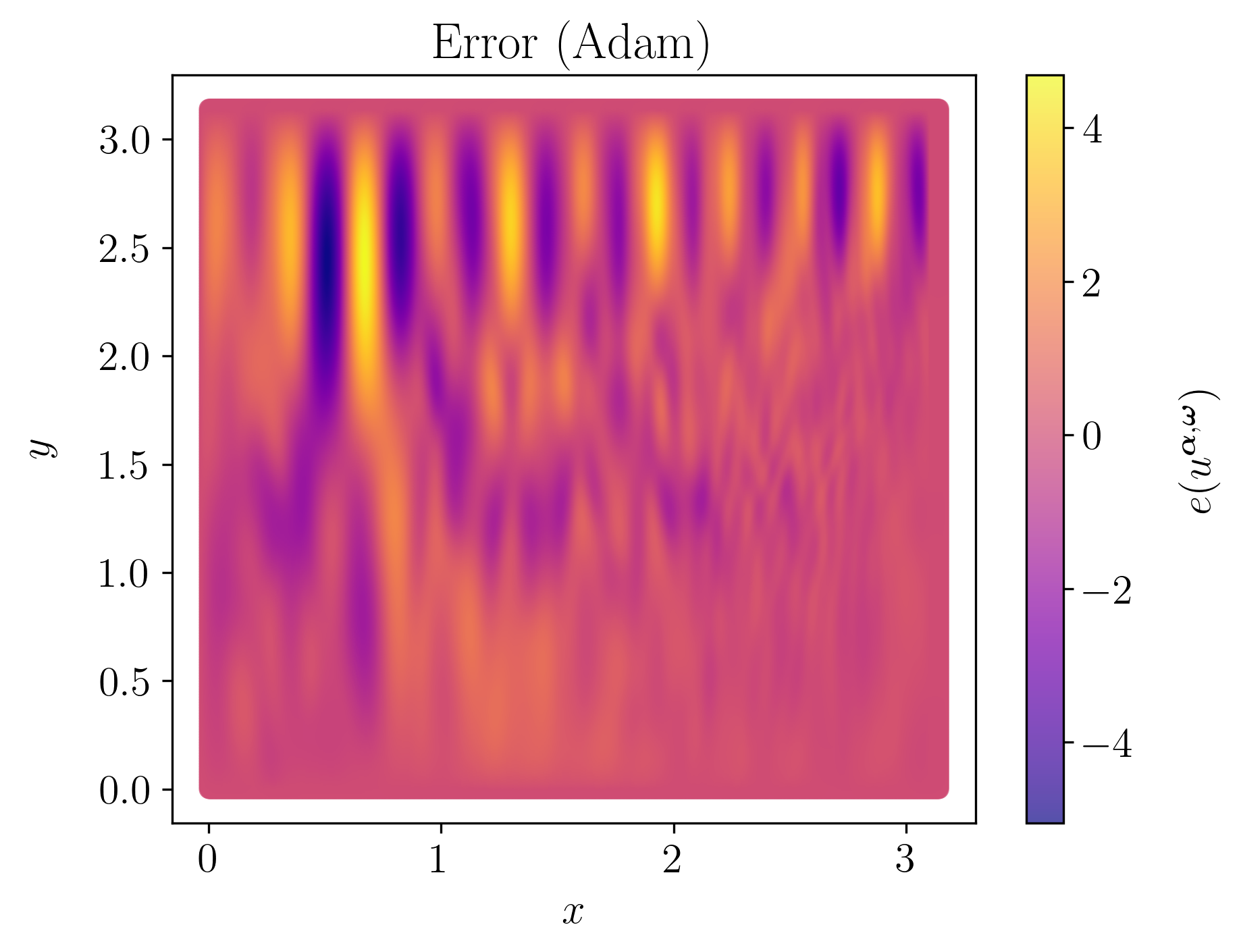}
    \end{subfigure}\hfill
      \begin{subfigure}[b]{0.49\textwidth}
		\includegraphics[scale=0.49]{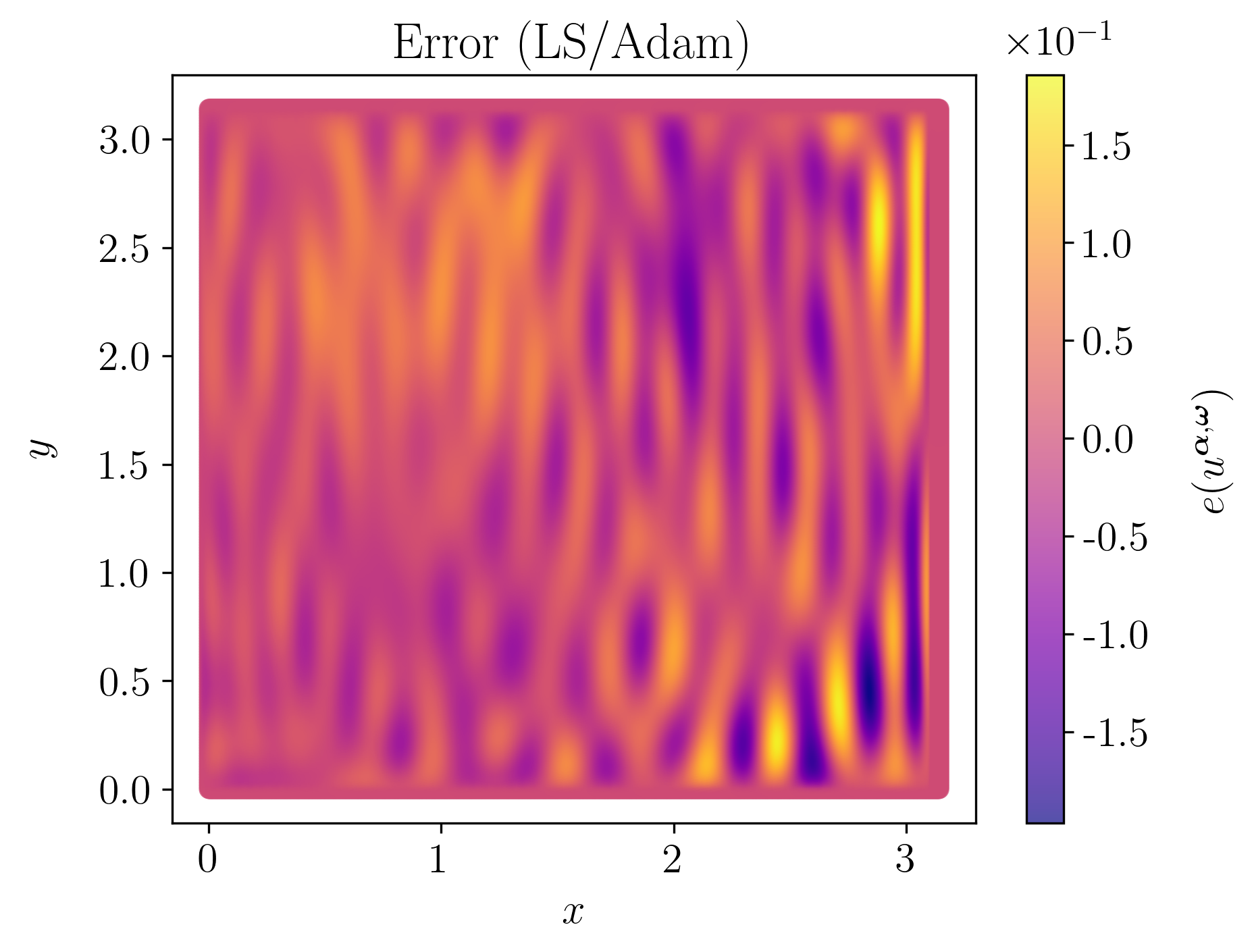}
  \end{subfigure}
\caption{Training of VPINNs with manufactured solution $u^*(x,y)=\sin(20x)\sin(2y)e^{\frac{x+y}{2}}$ using either Adam or LS/Adam. The resulting final relative errors are $29.24\%$ and $1.51\%$ for Adam and LS/Adam, respectively.}
\label{fig:experiment5ten}
\end{figure}

For $\eta=10$, we consider $N=128$, $M= 64\times 16$, $K = 512\times 128$, and $1\mathord{,}024\times 256$ for validation. The conclusions drawn from this 2D high-frequency case (see \Cref{fig:experiment5ten}) are similar to that of the 1D case in \Cref{section5.3}. We again observe plateau behavior in the first training iterations when using Adam alone due to spectral bias. When using LS/Adam, we obtain better convergence.

\pagebreak

\section{Conclusions and future work}\label[section]{section7}

We explored a strategy to speed up the convergence of gradient-descent-based algorithms for Variational Physics-Informed Neural Networks (VPINNs). We proposed an enhanced optimization scheme founded on the minimum-residual spirit of Petrov-Galerkin methods that rely on the use of Least Squares (LS). Appropriately combining this with conventional gradient-descent optimization approaches, we boost the performance of VPINNs.

Our work also focuses on automatic differentiation (AD). We explained how replacing the traditional backward mode with forward mode is critical during the LS setup. Moreover, we proposed an ultraweak-type scheme, called UltraPINNs, that avoids applying AD at this step. UltraPINNs are presumably cheaper than VPINNs because they avoid applying AD to use LS at each training iteration; however, we have observed that this tends to produce higher numerical integration errors compared to VPINNs when employing the same number of integration points (especially in singular problems). We leave a more in-depth analysis of this issue for future work. 

In our experiments, we observed substantial improvements when applying LS followed by Adam versus Adam alone. In all cases, we achieved much better accuracy thresholds and improved the slow spectral bias behavior of neural networks during training. Moreover, the proposed approach can be extended to other loss functions, including those derived from the Deep Ritz and traditional PINN methods. Additionally, integrating LS with the (L-)BFGS optimizer, which leverages second-order gradient descent approximations, might be of high interest to further improve accuracy in a limited number of iterations. Furthermore, we leave as future work studying the impact of stochastic numerical integration when using LS to train neural networks. Finally, we believe that replacing the discretization of the test space with a neural network could offer a promising alternative.

\section*{Acknowledgments}

This work has received funding from the following research grants of projects: European Union’s Horizon Europe research and innovation programme under the Marie Sklodowska-Curie grant agreement No 101119556. TED2021-132783B-I00 funded by MICIU / AEI / 10.13039 / 501100011033 and by the European Union NextGenerationEU / PRTR; PID2023-146678OB-I00 funded by MICIU / AEI  / 10.13039 / 501100011033 and by the European Union NextGenerationEU / PRTR; “BCAM Severo Ochoa” accreditation of excellence CEX2021-001142-S funded by MICIU / AEI / 10.13039 / 501100011033; Misiones Project IA4TES MIA.2021.M04.0008 funded by the Ministry of Economic Affairs and Digital Transformation; Basque Government through the BERC 2022-2025 program; BEREZ-IA (KK-2023/00012), funded by the Basque Government through ELKARTEK; Consolidated Research Group MATHMODE (IT1456-22) given by the Department of Education of the Basque Government; BCAM-IKUR-UPV/EHU, funded by the Basque Government IKUR Strategy and by the European Union NextGenerationEU/PRTR; and the Chilean grant ANID FONDECYT No. 1240643.

\appendix

\section{Numerical integration}\label{appendixIntegration}

For the 1D scenario, we consider a two-point \emph{composite} quadrature $Q=(Q_1,Q_2,\ldots,Q_P)$ of the form:

\begin{equation}
\int_\Omega I(x) \; dx = \sum_{p=1}^P \int_{\Omega_p} I(x)\; dx \approx \sum_{p=1}^P \underbrace{\sum_{k=1}^{2} \gamma_{k,p} \; I(x_{k,p})}_{Q_p(I)} =: Q(I),
\end{equation} where $I:\Omega=(0,\pi)\longrightarrow\mathbb{R}$ is a given integrand function, $\{\Omega_p=(a_p, a_{p+1})\}_{p=1}^P$ is a partition of $\Omega=(0,\pi)$, and $x_{1,p}$ and $x_{2,p}$ are two interior integration points with associated integration weights $\gamma_{1,p}$ and $\gamma_{2,p}$, respectively.

Integration points $x_{1,p}$ and $x_{2,p}$ are obtained via an affine transformation of a \emph{random} variable $X_p$ following a uniform distribution in $(-1,1)$. The resulting quadrature point and weight formulas are shown below:
\begin{equation}
\begin{cases}
x_{k,p} = (-1)^{k} \; \frac{a_{p+1}-a_p}{2} \; X_{p} + \frac{a_{p+1}+a_p}{2},\\
\gamma_{k,p} = \frac{a_{p+1}-a_p}{2},
\end{cases}\qquad k\in\{1,2\}, \quad 1\leq p\leq P.
\end{equation} By design, $x_{2,p}$ is the symmetric of $x_{1,p}$ with respect to the center of $\Omega_p$, and weights $\gamma_{1,p}$ and $\gamma_{2,p}$ are designed so that $Q$ is \emph{exact for piecewise-linear functions}.

In general, the resulting random variable $Q$ is \emph{unbiased}, i.e., the expectancy of $Q(I)$ coincides with the corresponding exact integral value,
\begin{equation}
    \mathbb{E}[Q(I)] = \int_{\Omega} I(x)\; dx.
\end{equation} 

We extend the above scheme to the two-dimensional domain $\Omega=(0,\pi)^2$ by following a Cartesian-product construction.

\bibliography{bibliography}

\end{document}